\documentclass{article}

\pdfoutput=1

\usepackage{amsmath,amssymb}
\usepackage{url}
\usepackage{ulem}
\usepackage{dsfont}
\usepackage{xcolor}
\usepackage{enumerate}
\usepackage{caption,subcaption}
\usepackage{graphicx}
\usepackage{placeins}

\usepackage{titlesec}
\titleformat{\chapter}[display]{\bfseries}{\huge\chaptertitlename~\thechapter}{10pt}{\LARGE}

\begin{document}

\title{Estimation of a discrete probability under constraint of
  $k-$monotony.
}

\author{Jade Giguelay}

\maketitle

jade.giguelay@ens-cachan.fr\\
http://maiage.jouy.inra.fr/jgiguelay\\

Laboratoire de Math\'ematiques d'Orsay, Universit\'e Paris-Sud, CNRS, Universit\'e Paris-Saclay, 91405 Orsay, FRANCE.\\
MaIAGE INRA, Universit\'e Paris-Saclay, 78350 Jouy-en-Josas, FRANCE.\\

\tableofcontents

\paragraph{Abstract} We propose two least-squares estimators of a discrete probability under the constraint of $k-$monotony and study their statistical properties. We give a characterization of these estimators based on the decomposition on a spline basis of $k$-monotone sequences. We develop an algorithm derived from the Support Reduction Algorithm and we finally  present a simulation study to illustrate their properties.


\renewcommand{\leq}{\leqslant}
\renewcommand{\geq}{\geqslant}
\renewcommand{\P}{\mathbb{P}}
\newcommand{\Z}{\mathcal{Z}}

\newtheorem{Theo}{Theorem}
\newtheorem{Prop}{Property}
\newtheorem{cor}{Corollary}
\newtheorem{Lemme}{Lemma}
\newtheorem{df}{Definition}
\newtheorem{Rem}{Remark}
\newtheorem{Exemple}{Example}

\newcommand{\cqfd}{\hfill{$\Box$}}
\newcommand{\ds}{\displaystyle}
\newcommand{\N}{\mathbb{N}}
\newcommand{\R}{\mathbb{R}}
\newcommand{\eref}[1]{(\ref{#1})}
\newcommand{\E}{\mathbb{E}}
\newcommand{\V}{\mathbb{V}}
\newcommand{\C}{\mathcal{C}}
\newcommand{\cI}{\mathcal{I}}
\newcommand{\IC}{\mathrm{CI}}
\newcommand{\SE}{\mathrm{SE}}

%%%%%%%%%%%%%%%%%%%%%%%%%%%%%%%%%%%%%%%%%%%%%%%%%%%%%%%%%%%

\section{Introduction}

The estimation of a density under shape constraint is a statistical problem that was first raised by Grenander \cite{grenander1956theory} in 1956 in the case of a density under monotony constraint.
Over the past 30 years, there has been several studies of estimators under shape constraint, most of them
being maximum likehood estimators or least squares estimators. In these cases, the authors characterize the estimators, study the asymptotic law and the rate of convergence and discuss the implementation. For such studies, the constraints are, for example, the monotony, the convexity or the log-concavity (if $\log(f)$ is concave, $f$ is log-concave) and the $k-$monotony.\\

The $k-$monotony notion was introduced by Knopp \cite{knopp1925mehrfach} in 1929 for  discrete functions: it generalizes to $k^{th}$ order the notion of convex series (or $2-$monotone series) and corresponds to the positivity of a $k$-th derivative fonction. In 1941, Feller \cite{feller1939completely} extended that definition to $k-$monotone continuous functions and Williamson \cite{williamson1955multiply} enabled characterizing these $k-$monotone functions 
with their decomposition in spline basis :
\begin{Prop}
(Williamson, 1955)\\
Let $g$ be a continuous function. Let $k\geqslant 2$. The function $g$ is $k-$monotone if and only if there exists a nonnegative mesure $\mu$ on $\mathbb{R}^{*}$ such that:
\begin{eqnarray*} g(x)=\int_{0}^{\infty}(t-x)_{+}^{k-1}d\mu (t).
\end{eqnarray*}
\end{Prop}
Consequently $k-$monotone functions can be described with an integral form. The estimation of a $k-$monotone distribution has been studied by Balabdaoui et al. (\cite{balabdaoui2004nonparametric},\cite{balabdaoui2007estimation},\cite{balabdaoui2009limit}): they proposed the maximum likehood and the least-squares estimators under $k-$monotony constraint for the continuous space and studied their theoretical properties (consistency and rate of convergence) as well as their limit distribution. They also discussed the adaptation of an algorithm proposed by Groeneboom et al. \cite{groeneboom2008support}.

Most of the work on estimation under shape constraint was focused on densities with a support on $\mathbb{R}$ or on an interval, but recently, discrete probabilities have gained  interest because of their numerous applications in ecology or financial mathematics (see \cite{durot_nonparametric} or \cite{lefevre2013multiply}). 
Jankowski and Wellner \cite{jankowski2009estimation} recently studied the estimation under monotony constraint and Balabdaoui et al. \cite{balabdaoui2014asymptotics} investigated the log-concave discrete densities. More recently, the estimation of a convex discrete distribution was treated by Durot et al. (\cite{durot_least}, \cite{durot_nonparametric}) and Balabdaoui et al. \cite{balabdaoui2014asymptotics}.
\\

In this article we propose two least-squares estimators of a $k$-monotone discrete probability with $k\geqslant 2$. The first one is the projection of the empirical estimator on the set of $k-$monotone sequences, the second one is the projection of the empirical estimator on the set of $k-$monotone probabilities. We show the existence of these estimators and give a characterization for each one of them which is based on the decomposition on a spline basis of $k-$monotone sequences showed by Lefevre and Loisel \cite{lefevre2013multiply}. Thanks to this characterization we generalize some results for the convex case ($k=2$, see \cite{durot_least}) to $k>2$, as for example the comparison with the empirical estimator (Theorem \ref{limnorme2.prop}).

However differences between the convex case and the case $k>2$ arised. First the projection of a discrete probability on the set of $k-$monotone sequences is not a probability in general when $k>2$. This structural property of the set of $k-$monotone functions, $k\geq 3$, justifies the definition of two different estimators while they are equivalent in the convex case.\\
Secondly the proofs of some other properties require new tools. In fact the results about the support of our estimator require control of the decreasing of the tail of $k-$monotone probabilities while truncation is sufficient in the convex case.
\\

Although the construction of our estimators is inspired by the work of Balabdaoui \cite{balabdaoui2004nonparametric} our results are not deduced from the continuous case. In fact, for $k\geq 3$, unlike for the convex case, we could neither construct a $k-$monotone density that goes through the points of a $k-$monotone sequence nor approach a $k-$monotone sequence with a $k-$monotone density. It is however interesting to note that connecting the points of a convex sequence can provide a convex continuous function because no differentiability assumption is required in the definition of the convexity. \\
Moreover the practical implementation of the estimator is structurally different from the continuous case. For the discrete case we implement the estimators using exact iterative algorithms inspired by the Support Reduction Algorithm described in Groeneboom et al. \cite{groeneboom2008support} and we discuss a practical stopping criterion (see Section \ref{algorithm.section}).\\
Differences with the continuous case also emerged when we consider the rate of convergence in terms of $l_2$-error since our estimators are consistent with typical parametric $\sqrt{n}$-rate of convergence (see Theorem \ref{limnorme2.prop}).
\\

The paper is organized as follows: the definition of the $k-$monotony
and some properties about $k-$monotone discrete functions are reminded
in Section \ref{characterization.section}, and a characterization of the estimator is given in
Section \ref{characterization2.section}. Statistical properties
about this estimator are then presented in Section
\ref{properties.section}. In
Section \ref{algorithm.section}, a method to implement the estimator
in practice using the Support Reduction Algorithm of Groeneboom et
al. \cite{groeneboom2008support} is presented. The stopping criterion for this algorithm,
which differs from the convex  case ($k=2$) is also discussed. In Section \ref{autreestim} we discuss the possibility to choose an estimator on the set of $k-$monotone sequences instead of the set of $k-$monotone probabilities. Finally a simulation study is given in Section \ref{simulation.section}.
\\

\section{Characterizing $k$-monotone sequences\label{characterization.section}}

Let us begin with a list of notation and definitions that will be used
throughout the paper.

The same notation is used to denote a discrete function $f : \mathbb{N}
\rightarrow \mathbb{R}^{+}$ and the corresponding sequence of real numbers
$(f(j), j \in \mathbb{N})$. For all $r \in \mathbb{N} \setminus
\left\{0 \right\}$, the classical $L^{r}$-norm of $f$ is defined as
follows: 
\begin{eqnarray*}
\|f\|_{r}   =\left(\sum_{j\geq 0} |f(i)|^{r}\right)^{1/r}, \text{  }
\|f\|_{\infty}   =\sup_{i\geq 0} |f(i)|,
\end{eqnarray*}
and we denote by $L^{r}(\mathbb{N})$ the set of functions $f$ such
that 
$\|f\|_{r}$ is finite. In particular $L^{2}$ is an Hilbert space and the associated scalar product is denoted $<,>$.\\ 

For any integer  $k\geq 1$, let $\Delta^{k} f$ be the $k^{\mbox{th}}$
differential operator of $f$ defined  for all $i\geq 0$ by the
following recurrence equation:
\begin{eqnarray*}
 \Delta^{1}f(i) & = & f(i+1) - f(i) \\
\Delta^{k}f(i) & = & \Delta^{k-1}f(i+1) - \Delta^{k-1}f(i).
\end{eqnarray*}

It is easy to see that the operator $\Delta^{k}$ satisfies the
following equation:
\begin{equation*}
\forall  i \in \mathbb{N}, \; \Delta^{k}
f(i)=\sum_{h=0}^{k}\binom{k}{h}(-1)^{k-h}f(h+i). 
\end{equation*}

\subsection*{Definitions}

\begin{itemize}
\item A sequence $f$ on $\mathbb{N}$ is $k$-monotone if
\begin{equation*}
 \left(-1\right)^{k} \Delta^{k}f(i) \geq 0 \mbox{ for all } i \in
 \mathbb{N}.  
\end{equation*}
\item Let $f$ be a $k$-monotone sequence. The integers $i$ such that
$\left(-1\right)^{k} \Delta^{k}f(i) > 0$ are called the $k$-knots of
$f$. If for all integer $i$ in the support of $f$, the quantities
$\left(-1\right)^{k} \Delta^{k}f(i)$ are strictly positive, $f$ is
said to be strictly $k$-monotone.
\item The maximum $s_f$ of the support of $f$ is  defined as 
\begin{equation*}
s_{f}  = \min_{j \geq 0} \left\{ \forall i > j, p(i) = 0\right\}
\end{equation*}
and may be infinite. 
\end{itemize}
Let us remark that if the support of a $k$-monotone sequence $f$ is
finite, then $s_{f}$ is a $k$-knot.
\\

A $k-$monotone function on $L^{1}(\mathbb{N})$ is for example a
non-negative and 
non-increasing polynomial function of degree $k-1$, such as $f(i) =
\max(0,m-i)^{k-1}$ for  some positive constant $m$.

It should be noticed that if a sequence $f$ is
$k$-monotone for some $k\geq 2$, then for all $j < k$, it
is strictly $j-$monotone on its support. This property, shown in
Section~\ref{monotony.sec}, is not true in general in the continuous
case (see Balabdaoui \cite{balabdaoui2004nonparametric}
for example).

Finally we will denote by ${\mathcal S}^{k}$ the set of $k$-monotone
sequences that are in $L^{1}(\mathbb{N})$, and by  ${\mathcal P}^{k}$
the set of 
$k$-monotone probabilities on $\mathbb{N}$. We denote by ${\mathcal
  P}$ the set of probabilities on $\mathbb{N}$.

\subsection*{Decomposition on a spline basis}

The characterization of $k$-monotone functions defined on
$\R$ as a mixture of polynomial functions has been established by
L\'evy~\cite{levy1962extensions}, and the inversion formula that specifies the mixture
function follows from the results of Williamson~\cite{williamson1955multiply}  (see Lemma
1. in~\cite{balabdaoui2007estimation} for example). In the case of
$k$-monotone sequences a similar decomposition has been simultaneously
established for convex sequences by Durot et al.~\cite{durot_least},
and in the more general case of $k$-monotony by Lefevre and Loisel~\cite{lefevre2013multiply}. Many of our proofs will rely on this decomposition.

For any integer $k$, let us define a basis of spline sequences 
$(\bar{Q}_{j}^{k})_{j\in\mathbb{N}\setminus \left\{0\right\}}$ as follows:
\begin{eqnarray}
\label{splinedef}
\forall i \in \mathbb{N}, \;
\bar{Q}_{j}^k(i)=
\binom{j-i+k-1}{k-1} \mathbb{I}_{\{j\geqslant i\}}
={(j-i+k-1)\ldots (j-i+1)\over (k-1)!} \mathbb{I}_{\{j\geqslant i\}}.
\end{eqnarray}
Let $m_{j}^{k}$ be the mass of $\bar{Q}_{j}^k$:
$m^k_j=\sum_{i=0}^{j}\bar{Q}^{k}_{j}(i)$ and $Q^k_j=\bar{Q}^k_j/m_j^k$ the 
normalized spline. We can now formulate the mixture representation of
$k$-monotone sequences.

\begin{Prop}\label{decomposition2}
Let $f\in L^1(\mathbb{N})$. 
\begin{itemize}
\item The sequence $f$ is $k-$monotone if and only if there exists a
  positive measure $\pi$ on $\mathbb{N}$, such that for all
  $i\in\mathbb{N}$,  $f(i)$ satisfies:
\begin{eqnarray}
\label{decomposition}
f(i)=\sum_{j\geqslant 0} \pi(j)Q^{k}_{j}(i)=\sum_{j\geqslant i}\pi(j)Q^{k}_{j}(i).
\end{eqnarray}
\item If $f$ is $k-$monotone, the measure $\pi$ is unique and defined as follows:
\begin{eqnarray}
\label{pi}
\forall j\geqslant 0, \; \pi(j)=(-1)^{k}\Delta^{k}f(j)m^k_j.
\end{eqnarray}
\item If $f$ is $k-$monotone, $\sum_{i=0}^{\infty}f(i)=\sum_{j=0}^{\infty}\pi_j.$
\end{itemize}
In particular $Q_{j}^{k}$ is $k-$monotone, and the set of $k$-knots of $f$
is the set of integers $j$ such that $\pi(j)$ is strictly positive.
\end{Prop}

These properties are shown in Lefevre
and Loisel~\cite{lefevre2013multiply}.

From this property, it appears that monotone discrete probabilities are
mixture of uniform distributions, convex probabilities are mixture of
triangular distribution, \ldots, $k$-monotone probabilities are
mixture of splines with degree $k-1$. 

\section{Constrained least-squares estimation on the set of
  $k$-monotone discrete probabilities \label{lse.st}}

Suppose that we observe $n$ i.i.d random variables, $X_1,\ldots,X_n$  with
distribution $p$ defined on ${\mathbb N}$, such that for all $i=1,
\ldots, n$, and $j \in {\mathbb N}$,  $p(j) = P(X_{i}=j)$. We propose to
build an estimator of $p$  that satisfies the $k$-monotony
constraint. Since the projection on the set of $k-$monotone sequences $\mathcal{S}_k$ is not a probability in general for $k\geqslant 3$ (see Section \ref{autreestim}) we consider the least-squares estimator
$\widehat{p}$ defined as follows:
\begin{equation}
 \hat{p}=\text{argmin}\left\{\| f-\tilde{p}\|^{2}_{2} , f\in {\mathcal
   P}^k\right\}
\label{lse.eq}
\end{equation}
 where $\tilde{p}$ is the empirical estimator of $p$: 
\begin{equation*}
\forall j \in {\mathbb N}, \;  \tilde{p}(j)=\frac{1}{n}\sum_{i=1}^{n}\mathbb{I}_{\{X_i =j\}}.
\end{equation*}

Since the set ${\mathcal P}^k$ of $k-$monotone discrete
probabilities  is convex and closed in the Hilbert space
$L^{2}(\mathbb{N})$, it follows, from the projection theorem for
Hilbert spaces,  that $\widehat{p}$ exists and is unique.

\subsection{Characterizing the constrained least-squares estimator\label{characterization2.section}}

A connerstone for deriving some statistical properties of our estimator is the following characterization of $\hat{p}$.
Let us begin with a few notation. 
For any positive sequence $f$ in $L^{1}({\mathbb N})$ we define the $j$-th primitive of
$f$ as follows: for all $l\in\mathbb{N}$
\begin{eqnarray*}
 F^{1}_{f}(l) & = & \sum_{i=0}^{l}f(i), \\
F^{j}_{f}(l) & = & \sum^{l}_{i=0}F^{j-1}_{f}(i) \mbox{ for all } j\geq 2.
\end{eqnarray*}
Moreover, we define the quantity $\beta(f)$
\begin{eqnarray}
\beta(f)=\sum_{i=0}^{\infty}f(i)(f(i)-\tilde{p}(i)).
\label{beta.eq}
\end{eqnarray}

\begin{Theo}\label{theo2}
Let $f \in {\mathcal P}$. The projection $\hat{p}$ defined as Equation (\ref{lse.eq}) is the unique $k-$monotone probability $f$ satisfying:
\begin{enumerate}
\item For all $l\in\mathbb{N}$,
\begin{eqnarray}
\label{criterearret2}
F^{k}_{f}(l)- F^{k}_{\tilde{p}}(l)\geqslant \beta(f)m^k_l.
\end{eqnarray}
\item If $l$ is a $k-$knot of $f$, then the previous inequality is an equality. 
\end{enumerate}
\end{Theo}
The proof of this theorem is given in Section~\ref{theo2.st}. It uses
the connections between successive primitives of the spline sequences
$(Q^k_j, j \in {\mathbb N})$.

 In the particular case of convexity
 the same result can be established with 0 in place of
$\beta(f)$, see Lemma 2 in~\cite{durot_least}. 
Let us recall that in that case, we have the nice property that
the least-squares estimator over convex sequences is a
convex probability distribution. This property is no longer satisfied
when $k\geq 3$. 
We will come back to this point in Section~\ref{autreestim}.

\subsection{Support of $\hat{p}$}
 
A key feature of the estimator $\hat{p}$ is that its support is finite. Let us denote by $\widehat{s}=s_{\widehat{p}}$, respectively 
$\widetilde{s}=s_{\widetilde{p}}$, the maximum of the support of
$\widehat{p}$, respectively $\widetilde{p}$. 

\begin{Theo}\label{supportfini.prop}
Let  $\widehat{p}$ be the
least-squares estimator defined by Equation~(\ref{lse.eq}).
\begin{enumerate}
\item The support of 
$\hat{p}$ is finite.  
\item If $\widehat{s} \geq \widetilde{s}+1$, then $\Delta^{k}\widehat{p}(i) = 0 $
\begin{itemize}
\item for all $i \in \left[ \tilde{s}-k+2,
    \hat{s}-1\right]$ if $k$ is even,
\item for all $i \in \left[ \tilde{s}-k+2,
    \hat{s}-2\right]$ if $k$ is odd.
\end{itemize}
\end{enumerate}
\end{Theo} 

In the particular case of convexity, when $k=2$, it is shown that
$\widehat{s} \geq \widetilde{s}$.
The question whether such a property still holds for $k\geqslant 3$ remains open.

\subsection{Statistical Properties\label{properties.section} of $\widehat{p}$ when $p$ is $k$-monotone}

Let us now evaluate the behaviour of $\widehat{p}$ for estimating
a  probability, in particular, how does it compare with the empirical
estimator $\widetilde{p}$. It is proved in the following theorem that the
constrained least-squares estimator is closer (with respect to the
$L^{2}$-norm) to any $k$-monotone
probability than is $\widetilde{p}$. 

\begin{Theo}\label{limnorme2.prop}
For any  $k-$monotone probability $f$, the following inequality is satisfied:
\begin{eqnarray}
\vert\vert f-\hat{p}\vert\vert_{2}\leqslant \vert\vert
f-\tilde{p}\vert\vert_{2}.
\label{majNorm.eq}
\end{eqnarray}
If $\tilde{p}$ is not $k-$monotone, then the  inequality is strict.

Moreover if $p$ is $k-$monotone  and 
if there exists $i\in\mathbb{N}$ such as $\Delta^{k}p(i)=0$, then for all
$k-$monotone probability $f$, we have :
\begin{eqnarray*}
\underset{n\rightarrow \infty}{\liminf} \; \mathbb{P}(\vert\vert f-\hat{p}\vert\vert_{2}<\vert\vert f-\tilde{p}\vert\vert_{2})\geqslant 1/2.
\end{eqnarray*}
\end{Theo}

In particular, if $p$ is $k$-monotone and not
strictly $k$-monotone, the
estimator $\widehat{p}$ is strictly closer to $p$ than is
$\widetilde{p}$ with probability at least 1/2. This theorem is a
straightforward generalization of Theorem~4 in~\cite{durot_least}
for the convex case and its proof is omitted.

In the following theorem   the  moments of the
distributions $\hat{p}$ and $\tilde{p}$ are compared. 
\begin{Theo}\label{moments}
For all $u\geqslant \max(1,k-3)$ and $0\leqslant a\leqslant \hat{s}$
the following inequality is satisfied:
\begin{eqnarray}
\sum_{i\geqslant 0}\vert i-a\vert^{u} \left(\widehat{p}(i)-\tilde{p}(i)
\right) \geq \beta\left(\widehat{p}\right)m(a,u),
\label{majmom.eq}
\end{eqnarray}
where $\beta$ is defined at Equation~(\ref{beta.eq}) and $m(a,u)=\sum_{i=0}^a (a-i)^u$.

Moreover $\hat{p}(0)-\tilde{p}(0)\geqslant \beta(\hat{p})$.
\end{Theo}

If $\widehat{p}$ satisfies
$\beta\left(\widehat{p}\right)=0$, the result is the same as the one
obtained in the convex case. In fact, it will be
stated in Section~\ref{autreestim}, that $\beta(\hat{p}) \leq 0$, and
that 
$\beta(\hat{p})=0$ if  the mimimizer of $\left\|f -
  \widetilde{p} \right\|^{2}$ over $f \in {\mathcal S}^{k}$ equals the
mimimizer of $\left\|f - \widetilde{p} \right\|^{2}$ over $f \in
{\mathcal P}^{k}$. This is  the case if $\widetilde{p}$ is
$k$-monotone, or if $k=2$.

\subsection{Asymptotic properties of $\widehat{p}$ \label{propertiesAs.section}}

In this section we consider the asymptotic properties of $\widehat{p}$
when the sample size $n$ tends to infinity. We first establish the
consistency of $\widehat{p}$ both in the case of a well-specified model or a
misspecified model.

\begin{Theo}\label{vitesse}Let $p_{\mathcal{S}_k}$ be the orthogonal projection
  of $p$ on the set $\mathcal{P}_k$. Then, for all $r\in[2,+\infty]$, the random
  variable $\sqrt{n} \left\| p_{\mathcal{S}_k}-\hat{p}\right\|_{r}$
is bounded in probability.
\end{Theo}
In particular, this theorem states that if the distribution $p$ is
$k$-monotone, then the convergence of $\widehat{p}$ to $p$ is of the
order $\sqrt{n}$ with respect to the $L^{r}$-norm. 

\paragraph{The case of a finite support}

In the particular case where the distribution $p$ is $k$-monotone and
has 
a finite support, we  characterize the asymptotic behaviour
of the $k$-knots of $\widehat{p}$, and  give an upper
bound for $\widehat{s}$, the maximum of the support of $\widehat{p}$.

\begin{Theo}\label{estimnoeud.prop}\label{maxsupport.prop}
Let $p$ be a $k-$monotone probability 
  with finite support. 
\begin{enumerate}
\item Let $j\in\mathbb{N}$ be a $k-$knot of $p$. 
Then with probability one there exists $n_0\in\mathbb{N}$ such that for all $n\geqslant n_0$, $j$ is a $k$-knot of $\hat{p}$.
\item  Let $s$, respectively $\hat{s}$, be the maximum of the support
  of $p$, respectively $\hat{p}$. Then, with probability one, there exists $n_0\in\mathbb{N}$ such that for all $n\geqslant n_0$ we have
\begin{itemize}
\item  $\hat{s}\leqslant s$ if $k$ is even.
\item $\hat{s}\leqslant s+1$ if $k$ is odd.
\end{itemize} 
\end{enumerate}
\end{Theo}

The proof of the first part of the theorem (see
Section~\ref{estimnoeud.proof}) is based on the fact 
that for all $j \in {\mathbb N}$, 
\begin{equation*}
\P \left(\lim_{n\rightarrow \infty} (-1)^{k}
  \Delta^{k}\widehat{p}(j) = (-1)^{k} \Delta^{k}p(j)\right) = 1.
\end{equation*}
It follows that if $j$ is a $k$-knot of $p$, then $(-1)^{k}
\Delta^{k}\widehat{p}(j)$ will be strictly positive for $n$ large
enough. Conversely, if $j$ is not a $k$-knot of $p$,
which means that $\Delta^{k}p(j)=0$, then
$(-1)^{k}\Delta^{k}\widehat{p}(j)$ may be strictly 
positive for all $n$. Therefore, the set of $k$-knots of $\widehat{p}$
does not estimate consistently the set of $k$-knots of $p$.

Concerning the second part of the theorem, we can notice that  the
result we get concerning $\widehat{s}$ is weaker than what was obtained
in the convex case. Indeed, when $k=2$, we know that $\widehat{s} \geq
\widetilde{s}$, and consequently
that,  $\widehat{s}=s$ for $n$ large enough if $p$
has a finite support.

\section{Implementing the estimator $\hat{p}$\label{algorithm.section}}

The practical implementation of $\widehat{p}$ requires the use of a
specific algorithm that is composed of two parts. The first part
consists in solving the problem  
defined at Equation~(\ref{lse.eq}) for sequences $f$  whose
support is finite. More precisely, for a chosen positive integer $L$, we compute
$\widehat{p}_{L}$, the minimizer of $\|f - \widetilde{p}\|^{2}$ over
sequences $f \in {\mathcal P}^{k}$ whose support is included in
$\left\{0, \ldots, L\right\}$. The
second part consists   in checking if $\widehat{p}_{L} = \widehat{p}$. For that
purpose, starting from Theorem~\ref{theo2}, we propose a stopping
criterion that can be calculated in practice.

\subsection{Constrained least-squares estimation on a given finite
  support \label{Part1algo.st}}

We know from Property~\ref{decomposition} that if $f \in {\mathcal
  P}^{k}$, there exists a unique probability $\pi$ on ${\mathbb N}$,
such that $f$ and $\pi$ satisfy Equation~(\ref{decomposition}). 
Therefore, solving~(\ref{lse.eq}) is equivalent to minimizing on the
set of probabilities $\pi$ on ${\mathbb N}$, the 
criterion $\Psi(\pi)$ defined as follows:
\begin{equation*}
\Psi(\pi) = \sum_{i\geq 0} \left(
\sum_{j\geq i}\pi(j)Q_{j}^{k}(i) - \widetilde{p}(i)
\right)^{2}.
\end{equation*}

The first part of our algorithm computes
\begin{equation*}
 \widehat{p}_{L} = {\mathrm{ argmin}} \left\{ \| f -
   \widetilde{p}\|^{2}_{2}, \; f \in {\mathcal P}^{k}, s_{f} \leq
   L\right\}.
\end{equation*}
 The solution is given by $\widehat{p}_{L} = \sum_{j \geq 0} \widehat{\pi}_{L}(j)
 Q_{j}^{k}$ where $\widehat{\pi}_{L}$ is the minimizer of $\Psi(\pi)$
 over probabilities $\pi$ whose support is included in $\left\{0,
   \dots, L\right\}$:
\begin{equation}
\widehat{\pi}_{L} = {\mathrm{argmin}} \left\{ \Psi(\pi), \;
  \pi \in {\mathcal P}, s_{\pi} \leq L\right\}.
\label{piHat.eq}
\end{equation}

The algorithm we use to compute $\widehat{\pi}_{L}$ is based on the
support reduction algorithm introduced by Groeneboom et
al.~\cite{groeneboom2008support}. In the particular case where $k=2$ it
was developped by Durot and al.~\cite{durot_least}. Nevertheless when
$k\geq 3$ an adaptation is needed to guarantee that
$\widehat{\pi}_{L}$ is a probability. Let us underline that this
algorithm is proven to give the exact solution in a finite number of
steps.

For all $\nu \in \mathcal{P}$, let  $D_{\nu}\Psi$ be
the directionnal derivative function of $\Psi$ in the direction $\nu$
defined as follows: 
\begin{eqnarray*}
\forall \mu\in\mathcal{P}, D_{\nu}\Psi (\mu)=\underset{\epsilon\searrow 0^+}{\lim}\frac{1}{\epsilon}\big(\Psi((1-\epsilon)\mu+\epsilon\nu)-\Psi (\mu)\big).
\end{eqnarray*}
The Support Reduction Algorithm is based on the property that
$\widehat{\pi}_{L}$ is solution of~(\ref{piHat.eq}) if and only if the
directionnal derivative functions calculated in $\widehat{\pi}_{L}$ in
the directions $\nu=\delta_{j}$, where $\delta_{j}$ denotes the Dirac
probability in $\{j\}$, are non negative for all $0 \leq j \leq
L$. Moreover, these derivatives 
are exactly  0 for all $j$ in the support of $\widehat{\pi}_{L}$.

Starting from this property, the Support Reduction Algorithm is
composed of two steps. In the first step, the support of the current
probability $\mu$ is augmented by a point $j$ where
$D_{\delta_{j}}\Psi (\mu)$ is strictly negative (if any). In the
second step the minimisation of $\Psi(\mu)$ over sequences $\mu$ such
that $\sum_{j\geq 0}|\mu(j)|=1$ and 
whose support is the current support, is performed. The current
support  is reduced to obtain a positive sequence. 

Let us introduce notation used in the second step of the algorithm. For a set
$S=\{j_1,\ldots,j_s\}\subset\{0,\ldots,L\}$ we note 
\begin{itemize}
\item $Q_{S}$ the matrix whose
component $(Q_{S})_{i+1, \ell} = Q^{k}_{j_\ell}(i)$ for $0\leq i\leq
L$ and $j_\ell \in S$, $l=1,\hdots,s$,
\item $H_{S}$ the projection matrix $H_{S} =
Q_S(Q_S^TQ_S)^{-1}Q_S^T$, and
\item $\lambda_{S}$ the Lagrange multiplier
\begin{equation*}
 \lambda_S=\frac{\langle H_{S} \tilde{p}, \mathbb{I} \rangle-1}
{\langle H_{S} \mathbb{I},\mathbb{I}\rangle},
\end{equation*}
where $\mathbb{I}$ is the vector with $L+1$ components all equal to
1. 
\end{itemize}

The algorithm for computing $\widehat{\pi}_{L}$ for a fixed $L$ is
given at Table~\ref{algo.tb}. It is shown in Section~\ref{algo2.st} that this
algorithm returns $\widehat{\pi}_{L}$ in a finite number of steps. 

\begin{table}[h]
\fbox{\begin{minipage}{0.8\textwidth}
\begin{itemize}
\item
\textbf{Initialisation :}\\
S $\leftarrow$ $\{L\}$\\
$\pi \leftarrow \delta_{L}$
\item
\textbf{Step 1 :}\\
For all $j\in\{0,\ldots,L\}$ compute  $D_{\delta_j}\Psi(\pi)$.
\begin{itemize}
\item If for all $j\in\{0,\ldots,L\}$  $D_{\delta_j}\Psi(\pi)\geqslant
  0$, stop.
\item Else choose $j\in\{0,\ldots,L\}$ such as $D_{\delta_j}\Psi(\pi)<0$.
\\
$S' \leftarrow S+\{j\}$.\\
Go to step 2.
\end{itemize}
\item
\textbf{Step 2 :}\\
$\lambda$ $\leftarrow$ $\lambda_{S'}$\\
$\pi_{S'}$ $\leftarrow$ argmin$\{\Psi(\mu)+\lambda_{S'}(\sum_{j\in
  S'}\mu(j)-1)$,  supp$(\mu)\subset S'\}$.
\begin{itemize}
\item If for all $l\in S', \pi_{S'}(l)\geqslant 0$,\\
$\pi \leftarrow \pi_{S'}$\\
$S \leftarrow S'$\\
Return to step 1.
\item  Else 
\begin{eqnarray*}
l \leftarrow %\underset{j'\in S'}
{\text{argmin}}\{\epsilon_{j'}=\frac{\pi_{j'}}{\pi_{j'}-\pi_{S'}(j')},
\; j'\in S', \pi_{S'}(j')<\pi_{j'}\} 
\end{eqnarray*}
$S' \leftarrow S'-\{l\}$\\
Return to step 2.
\end{itemize}
\end{itemize}
\end{minipage}}
\caption{\label{algo.tb}} Algorithm for computing $\widehat{\pi}_{L}$
for a fixed $L$.
\end{table}

\subsection{Stopping criterion}

The second step of the algorithm is to find $L$ such that $\widehat{p}_{L} =
\widehat{p}$. 

The characterization of $\widehat{p}$ given by
Theorem~\ref{theo2} cannot help for practical
implementation because the necessary condition in that theorem
requires an infinite number of calculations. Nevertheless, if $f$ is a
$k$-monotone probability, with maximum support $s_{f}$, it is
possible to find an integer $M$ such that if $f$ satisfies
Inequality~(\ref{criterearret2}) for all $l \leq M$, then $f$
satisfies 
Inequality~(\ref{criterearret2}) for all $l > M$. Such a property results from
the writting of 
\begin{equation*}
P_{f} (l) = F_{f}^{k}(l) - F_{\widetilde{p}}^{k}(l) - \beta(f) m_{l}^{k}
\end{equation*}
as a polynomial function in the variable $l$. On the one hand, Property~\ref{theoclé},
shown in Section~\ref{theoclé.st},
states that $F_{f}^{k}(l) - F_{\widetilde{p}}^{k}(l)$ is a polynomial
function in $l$ of degree $k-1$ as soon as $l$ is greater than the maxima of the
support of $f$ and $\widetilde{p}$.

\begin{Prop}\label{theoclé}
Let $f$ be a discrete sequence with finite support and $s_{f}$ be the
maximum of its support. Let $\tau = \max(s_{f}, \widetilde{s})$, then for all $l \geq \tau+1$,  we have the following equalities:
\begin{eqnarray}F^{k}_{f}(l)-F^{k}_{\tilde{p}}(l)
&=&\sum_{j=1}^{k}\bar{Q}_{l-1}^{k-j+1}(\tau)\big(F^j_{f}(\tau)-F^{j}_{\tilde{p}}(\tau)\big) \nonumber\\
&= &\sum_{j=1}^{k} \frac{F^j_{f}(\tau)-F^{j}_{\tilde{p}}(\tau)}{(k-j)!}
\left((l-\tau+k-j-1) \ldots (l-\tau)\right) \label{pol1.eq}
\end{eqnarray}
\end{Prop}

On the other hand starting from Pascal's rule, it is easy to see that 
for all $k\geq 2$ and $l\geq 0$, 
\begin{equation}
 m_{l}^{k} = \bar{Q}^{k+1}_{l}(0) = \frac{(l+k)(l+k-1)\ldots(l+1)}{k!}.
\label{massespline}
\end{equation}
Putting Equations~\eref{pol1.eq} and~\eref{massespline} together, it
appears that for all $l \geq \tau = \max(s_{f}; \widetilde{s})$, there
exist coefficients $(a_{0}, a_{1}, \ldots, a_{k-1})$ such that 
\begin{equation*}
 P_{f} (l) = \sum_{j=0}^{k-1} a_{j} l^{j}.
\end{equation*}
Let $d$ be the degree of this polynomial (the smallest $j$ such that
$a_{j}=0$ for all $j \geq d+1$) and let $M$ be defined by 
\begin{equation*}
M = \max\left( 1+\frac{a_{d-1}}{a_d},\ldots,1+\frac{a_0}{a_d}\right),
\end{equation*}
By Cauchy's Theorem for  localization of polynomial's roots,
the largest root of $P_{f}(l)$ is bounded by $M$. Therefore if $a_{d}$
is positive, $P_{f}(l)$ is positive beyond $M$. This leads to the
following characterization of $\widehat{p}$ which is a corollary of
Theorem~\ref{theo2}. Its proof is omitted.

\begin{Theo}
\label{critere_k}
Let $f$ be a sequence in ${\mathcal P}$ with a finite support. Let $M$
and $a_{d}$ be defined as above. The two following assertions are equivalent:
\begin{enumerate}
\item The sequence $f$ satisfies 
\begin{enumerate}
\item $a_{d}$ is positive. 
\item $\forall l\leqslant M, F^k_{f}(l)- F^k_{\tilde{p}}(l)\geqslant
  \beta(f)m^k_l$. 
\item If $l$ is a $k-$knot of $f$, the previous inequality is an
        equality.
      \end{enumerate}
\item The sequence $f$ is exactly $\hat{p}$.
\end{enumerate}
\end{Theo}

This Theorem answers our initial problem, namely checking if $\widehat{p}_{L}$ equals
$\widehat{p}$: the coefficients $a_j$ depending on $\widehat{p}_{L}$ and
$\widetilde{p}$, can be calculated in practice, as well as
$M$. Nevertheless $M$ can be very large, leading to tedious
computation.  For small values of $k$ more efficient criteria can be
 proposed.  In particular, for $k=3$ and $k=4$, it is possible to
 obtain a characterization of $\widehat{p}$ that only depends on
 $\widetilde{s}$ and $\widehat{s}$. This is the object of the following theorem shown
 in Section~\ref{critere_k=34.st}.

\begin{Theo}\label{critere_k=34}
\label{critere_k}
Let $f$ be a sequence in ${\mathcal P}$ with a finite support and 
$s'
= \max\left\{s_{f}, \widetilde{s} \right\}$.  Let
$k\in \{3,4\}$. The two following assertions are equivalent:
\begin{enumerate}
\item The sequence $f$ satisfies 
\begin{enumerate}
\item $\forall l\leqslant s'+1, F^k_{f}(l)- F^k_{\tilde{p}}(l) \geqslant 
  \beta(f)m^k_l$. 
      \item If $l$ is a $k-$knot of $f$, the previous inequality is an
        equality.
      \item for all $2\leqslant j\leqslant k-1$, $F^j_f(s'+1)-
        F^j_{\tilde{p}}(s'+1)  \geqslant \beta(f)m^j_{s'+1}$.
      \item $\beta(f)\leqslant 0$.
\end{enumerate}
\item The sequence $f$ is exactly $\hat{p}$.
\end{enumerate}
\end{Theo}

When $k>4$ we are not able
to propose a similar stopping criterion.  Indeed the proof is based on
the properties of the spline function $Q_j^k$ for $k>4$  and in
particular, requires the calculation of the number of $k'$-knots of
$Q_j^k$ for $k'\geqslant k$, which is intractable.

\section{\label{autreestim} Constrained least-squares estimation on the set of $k$-monotone sequences}

By definition, our estimator $\widehat{p}$ is a probability. We could
have proposed to estimate $p$ by minimizing the least-squares criterion
under the constraint of $k$-monotony only. Let $\widehat{p}^{*}$
be
that estimator:
\begin{equation}
  \hat{p}^{*}=\text{argmin}\left\{\| f-\tilde{p}\|^{2}_{2} , f\in {\mathcal
   S}^k\right\}
\label{lseSk.eq}
\end{equation}
The
following property, shown in Section~\ref{beta.st},  establishes the link between both estimators.
\begin{Prop}\label{beta} Let $\widehat{p}$ and $\widehat{p}^{*}$ be
  defined at Equations~\eref{lse.eq} and~\eref{lseSk.eq}, and let
  $\beta$ be defined at Equation~\eref{beta.eq}. The coefficient $\beta(\hat{p})$ is null if and only if
  $\hat{p} = \hat{p}^{*}$.
\end{Prop}

In the particular case where $k=2$, $\hat{p}^{*}$ is exactly $\widehat{p}$ (see
Theorem 1 in~\cite{durot_least}). As soon as
$k\geq 3$ this property is no longer satisfied in general, and it is
proven in Section~\ref{masse.prop.st} that the mass of
$\widehat{p}^{*}$ is greater than or equal to 1.
This result was expected because a similar property was shown by
Balabdaoui~\cite{balabdaoui2004nonparametric} in the continuous
framework. 

To
illustrate this point in the discrete framework, let us consider the  projection of $\delta_1$
(the Dirac probability in 1) on the set of $3-$monotone
sequences. Some calculation (see Section~\ref{cex_masse.st} for a
proof) leads to the following result:
\begin{equation*}
{\mathrm{ Proj}}_{{\mathcal S}^{3}}\left(\delta_{1} \right)  = \frac{3}{238}\bar{Q}^3_5+\frac{1}{238}\bar{Q}^3_6,
\end{equation*}
and its mass is close to 1.14.

Nevertheless we can show (see Section~\ref{prop5.st}) the following asymptotic
result.
\begin{Prop}\label{prop5}Let $p$ be a $k-$monotone probability, and
  $\widehat{p}^{*}$ be defined at Equation~\eref{lseSk.eq}. Then, with
  probability one the mass of $\widehat{p}^{*}$ converges to one.
\end{Prop}

The properties shown in Section~\ref{lse.st} for the estimator $\widehat{p}$ hold true for
$\widehat{p}^{*}$. More precisely, the estimator $\hat{p}^{*}$ satisfies Theorems~\ref{supportfini.prop},
\ref{vitesse}, \ref{estimnoeud.prop}. Theorem~\ref{limnorme2.prop} is
also true for $\widehat{p}^{*}$ apart from the first assertion, where
Equation~\eref{majNorm.eq} is satisfied for any $k$-monotone sequence
$f$. Finally Theorem~\ref{moments} is true with 0 in place of 
$\beta(\widehat{p})(\widehat{s}+1-a)$ in Equation~\eref{majmom.eq}.

The implementation of $\widehat{p}^{*}$ is similar to that of
$\widehat{p}$ except that in the first part of the algorithm (see
Section~\ref{Part1algo.st})  the Support Reduction Algorithm can be
used without any modification at Step 2 (where the estimator of $\pi$
is constraint to have a sum of one). 
The stopping criterion used in the second part is obtained in the same
way as for $\hat{p}$. 
The proofs of these last results are omitted. They are based on Property~\ref{beta}. 

\section{Simulation\label{simulation.section}}

We designed a simulation study to assess the quality of the
 least-squares  estimator $\widehat{p}$  on the set of $k$-monotone
 probabilities,  as compared to the empirical estimator $\widetilde{p}$,
and to the least-squares  estimator $\widehat{p}^{*}$  on the set of
$k$-monotone
 sequences for $k\in\{2,3,4\}$.

\subsection{Simulation Design\label{simulationdesign}}

We considered two shapes for the distribution $p$: the spline
distribution $Q_{j}^{\ell}$ with $j=10$ and  $\ell \in\{2,3,4,10\}$,
and the Poisson distribution $\mathcal{P}(\lambda)$ for $\lambda \in
\{0.3,0.35,0.45,2-\sqrt{2},0.7,1\}$. Those two families of
distribution differ by the finiteness of their support, and by the
number of knots  in their decomposition on the spline
basis. Precisely, the distribution $Q_{j}^{\ell}$ has one $\ell$-knot
in $j$
while a $\ell$-monotone Poisson distribution has an infinite
number of $\ell$-knots. The following proposition, shown in
Section~\ref{Poiss.st}, gives the property of $k$-monotony for Poisson distributions.
\begin{Prop}\label{Poiss.prop}
Let $\mathcal{P}(\lambda)$ be the Poisson distribution with parameter
$\lambda$. For each $\ell \geqslant 1$, let $\lambda_\ell$ be defined
as the smallest root of the following polynomial function:
\begin{eqnarray*}
P_\ell(\lambda)=\sum_{h=0}^\ell (-1)^h \frac{(\ell!)^2}{h!((\ell-h)!)^2}\lambda^h.
\end{eqnarray*} 
Then $\mathcal{P}(\lambda)$  is $\ell$-monotone
if and only if $\lambda\leqslant \lambda_\ell$
\end{Prop}
Some simple calculation gives the following values: $\lambda_{1}=1$,
$\lambda_{2}=2-\sqrt{2}\simeq0.585$, $\lambda_3\simeq 0.415$,
$\lambda_4\simeq 0.322$, $\lambda_5\simeq 0.264$. 
Therefore the considered Poisson distributions $\mathcal{P}(\lambda)$
are $\left\{4, 3, 2, 2, 1\right\}$-monotone when
$\lambda$ belongs   to
$\{0.3,0.35,0.45,2-\sqrt{2},0.7\}$. When $\lambda=1$, the Poisson
distribution is not strictly decreasing.

For each distribution $p$, we considered several values for the sample
size $n$: $n \in \left\{20, 50, 100, 250, 500, 1000\right\}$. In some
 cases we also considered  very large values of $n$ in order to
 illustrate the asymptotic framework. We
denote by $\widetilde{p}_{n}$ the empirical estimator and by
$\widehat{p}_{n}^{k}$, respectively $\widehat{p}_{n}^{* k}$, the
least-squares estimator of $p$  on the set of $k$-monotone
 probabilities, respectively sequences. 
For each simulation
configuration, 1000 random samples were generated. 
 The simulations were
carried out with  R ({\tt www.r-project.org}), the code being available at
the following web address: http://www.maiage.jouy.inra.fr/jgiguelay

\subsection{Global fit}

To assess the quality of the estimators for estimating the
distribution $p$ we consider the $l_{2}$-loss and the Hellinger
loss. We have also considered the total variation loss, but the
results are not shown because they are very similar to those obtained
for the $l_{2}$-loss,  

\subsubsection{Estimators comparison based on the $l_{2}$-loss}

The $l_{2}$-loss between $p$  and any estimator of $p$, say $\widehat{q}$, is
defined as the expectation of the $l_{2}$-error, $l_{2}(p,\widehat{q})
= E(\|p-\widehat{q}\|_{2}^{2})$.

We first compared the quality of the fit of the estimators
$\widehat{p}_{n}^{k}$ 
and $\widetilde{p}_{n}$  by computing for each
simulated sample  $\|p-\widehat{p}_{n}^{k}\|^{2}_{2}$
and $\|p-\widetilde{p}_{n}\|^{2}_{2}$. The  $l_{2}$-losses were
estimated by the 
mean of 1000  independant replications of the $l_{2}$-errors. In all
simulation configurations, the  $l_{2}$-losses are decreasing
towards 0 
when $n$ increases. In what follows we will consider  the
ratios  $l_{2}(p,\widehat{p}_{n}^{k})/l_{2}(p,\widetilde{p}_{n})$ to
compare the estimators.

\begin{figure}
\centering
\begin{subfigure}[b]{0.45\textwidth}
\includegraphics[width=0.95\textwidth]{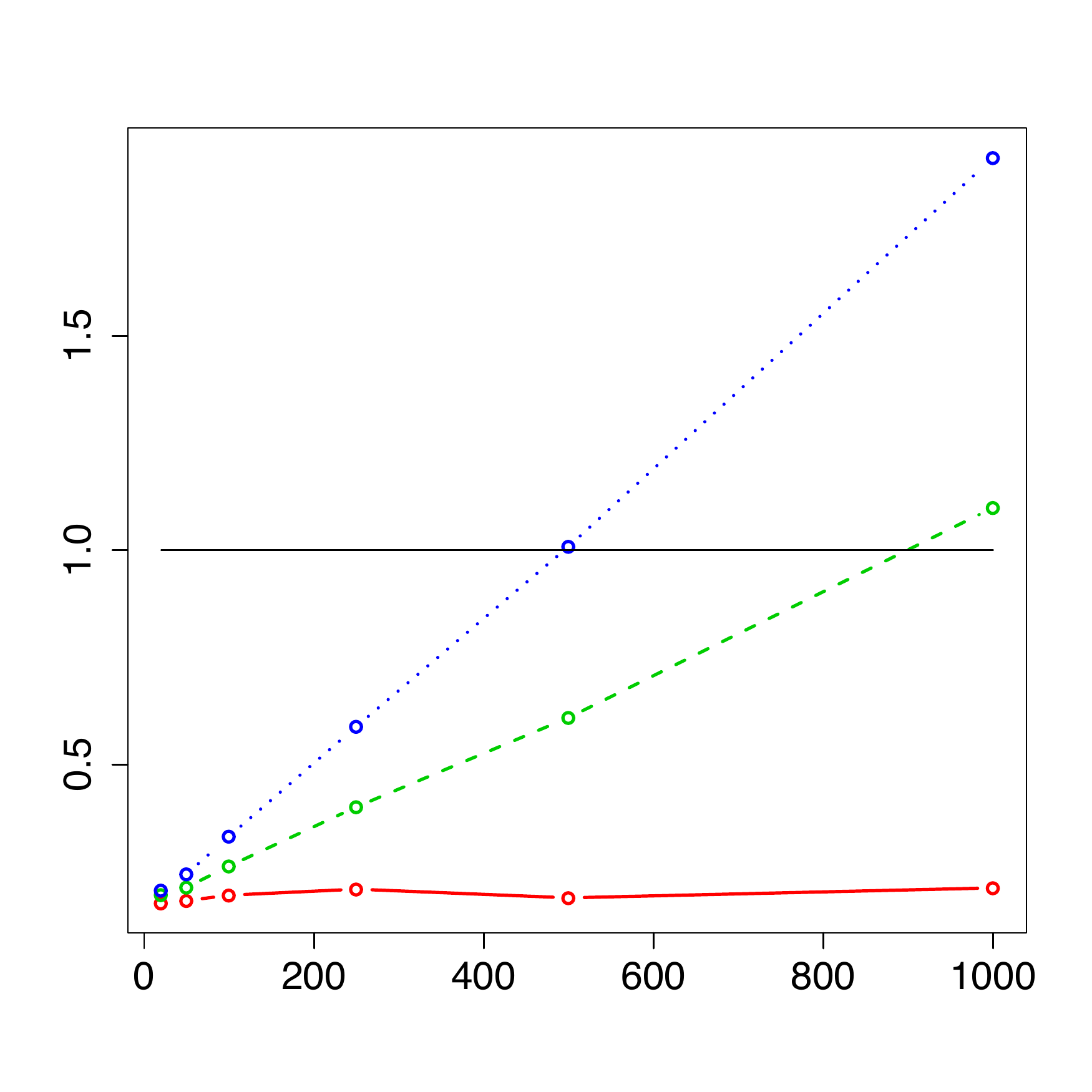}
\caption{$p=Q^2_{10}$.}
\label{fig1a}
\end{subfigure}
\hfill
\begin{subfigure}[b]{0.45\textwidth}
\includegraphics[width=0.95\textwidth]{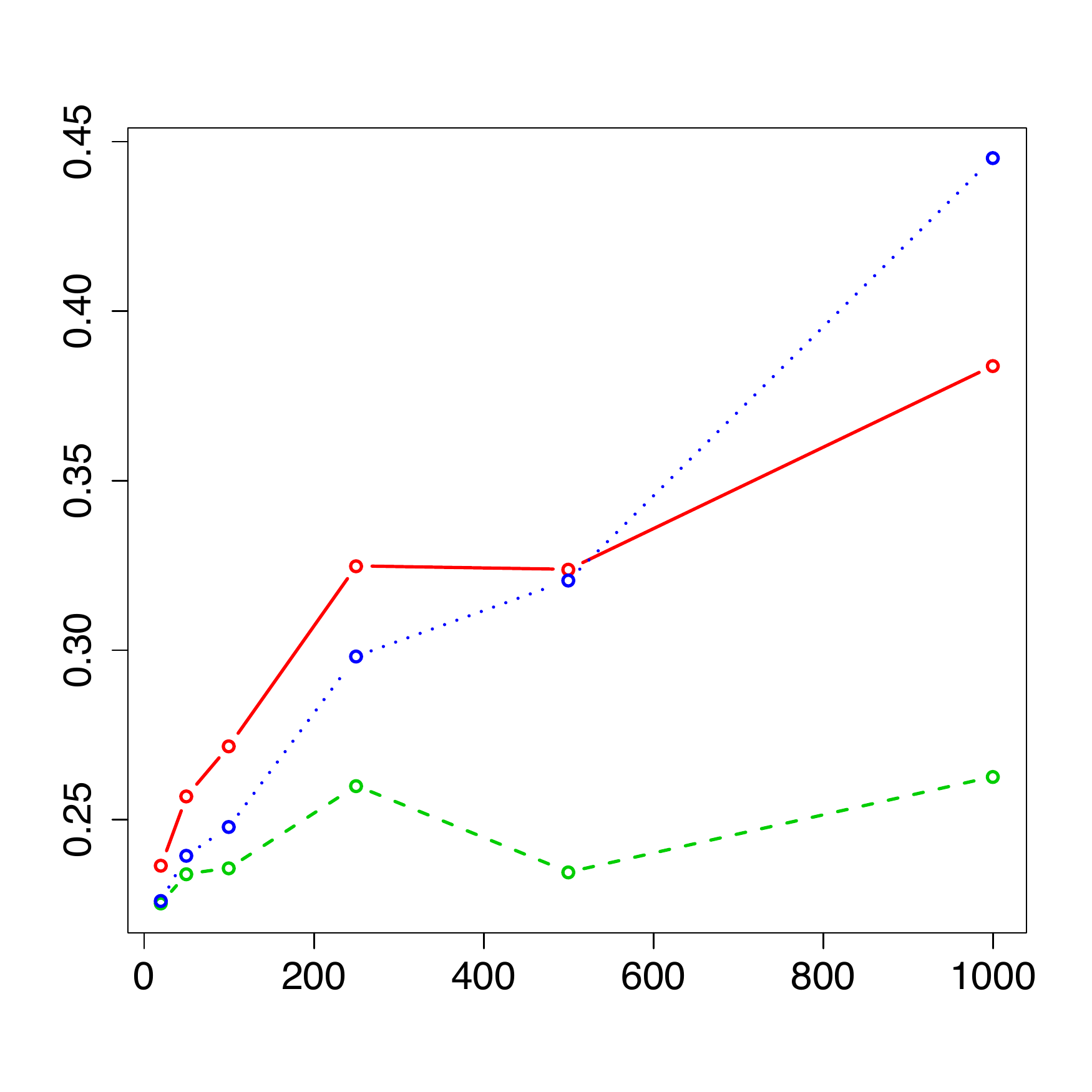}
\caption{$p=Q^3_{10}$.}
\label{fig1b}
\end{subfigure}
\begin{subfigure}[b]{0.45\textwidth}
\includegraphics[width=0.95\textwidth]{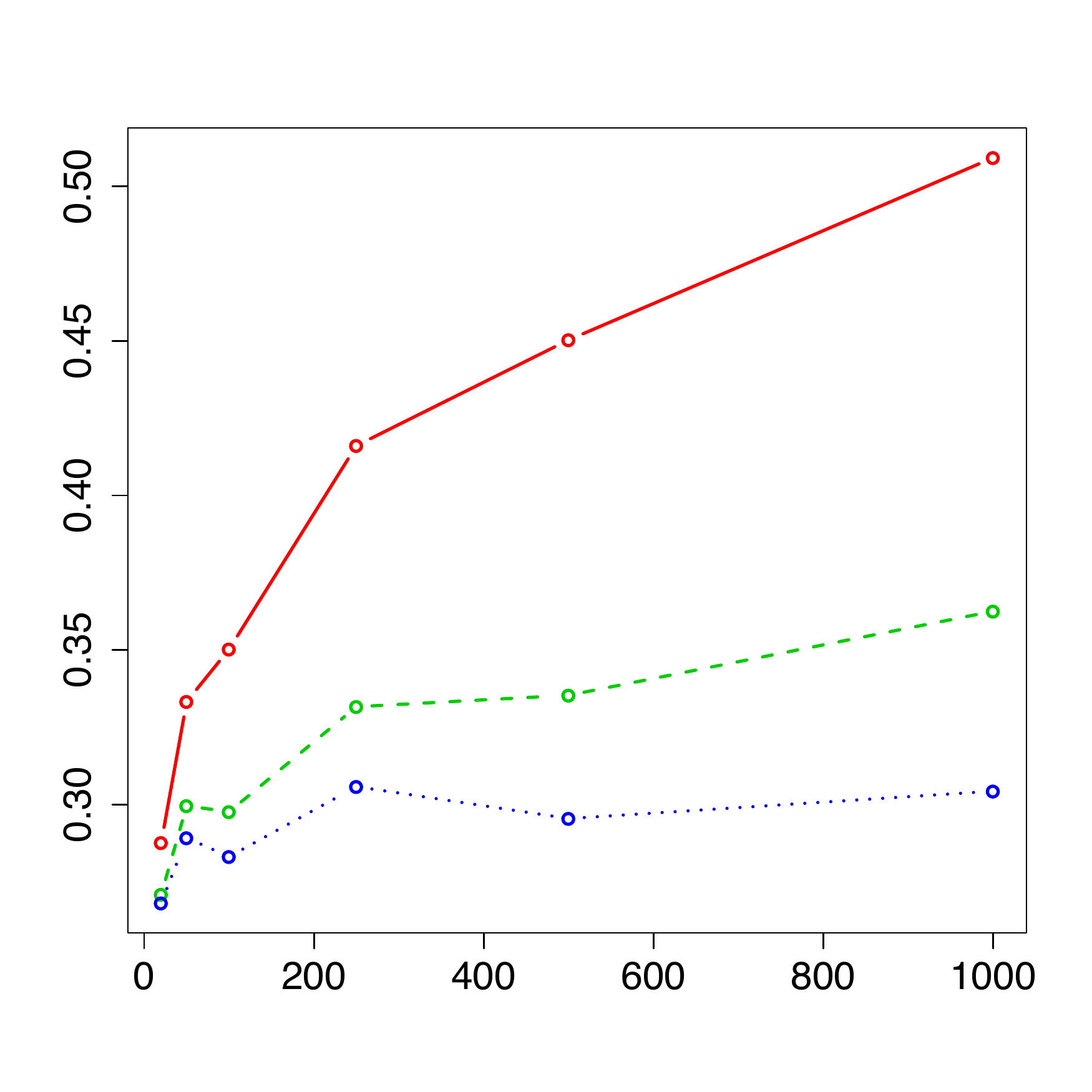}
\caption{$p=Q^4_{10}$.}
\label{fig1c}
\end{subfigure}
\hfill
\begin{subfigure}[b]{0.45\textwidth}
\includegraphics[width=0.95\textwidth]{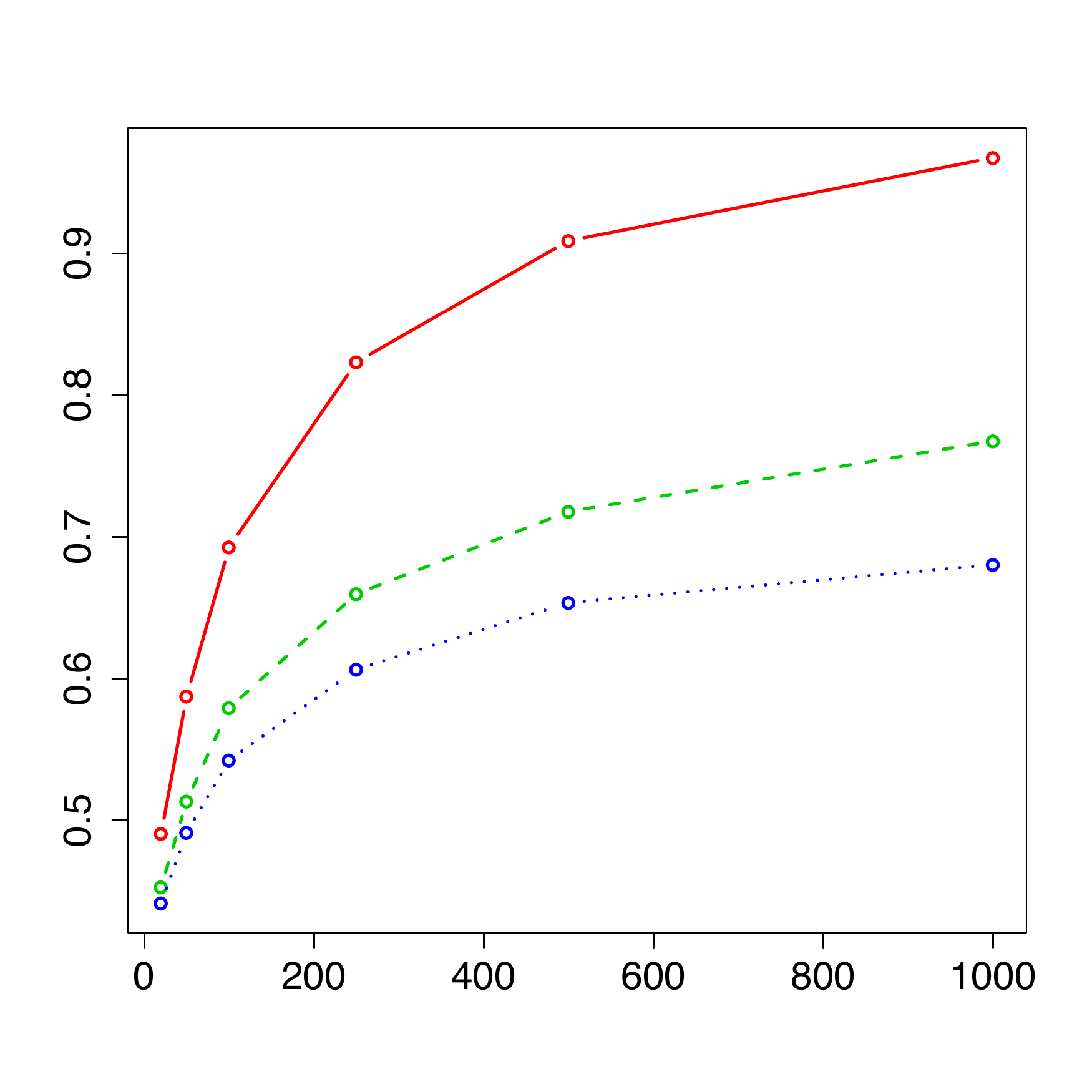}
\caption{$p=Q^{10}_{10}$.}
\label{fig1d}
\end{subfigure}
\caption{Spline distributions: ratio between the $l_2$-loss of 
  $\widehat{p}_{n}^{k}$ and the $l_2$-loss of
  $\widetilde{p}_{n}$ versus the sample size $n$:  for $k=2$ in "\textbf{\textcolor{red}{---}}",
  $k=3$ in "\textbf{\textcolor{green}{- -}}", $k=4$ in
  "\textbf{\textcolor{blue}{...}}". Each subfigure
  corresponds to the results  obtained with  $p = Q_{j}^{\ell}$, for
  $\ell\in\{2, 3, 4, 10\}$.}
\label{l2.fig}
\end{figure}

The results for the spline distributions $Q_{j}^{\ell}$ are presented on
Figure~\ref{l2.fig}. When $n$ is small, $\widehat{p}^{k}_{n}$ has
smaller $l_{2}$-loss than $\widetilde{p}_{n}$ whatever the value of
$k$. When $n$ tends to infinity, we have to consider two cases
according to the discrepancy between $k$ which defines the degree of
monotonty of the estimator, and $\ell$ which is the degree of monotony
of $p$. As it was expected considering Theorem~\ref{limnorme2.prop}, when $k \leq
\ell$, then the ratio is smaller than 1. Moreover we note that
the smaller the deviation $\ell-k$ is, the smaller the ratio. In
particular when $k=\ell$, the ratio tends to a constant strictly smaller than
1, while when  $k<\ell$, the ratio tends to 1. For
    example, when 
$\ell=4$, $k=3$, the ratio of the $l_{2}$-losses
equals $0.45$ for $n=10000$ and $0.80$ for $n=100000$. When $k>\ell$, the ratio tends to infinity.  For example, when $\ell=2$, $k=3$, the ratio of the $l_{2}$-losses
equals $9.93$ for $n=10000$ and $259$ for $n=100000$. This result was expected because the empirical estimator
$\widetilde{p}_{n}$ 
is consistent while our estimator is not. Indeed,
following Theorem~\ref{vitesse}, the ratio of $l_{2}$-losses is of
order $n \|p-p_{S_{k}}\|^{2}_{2}$ which is zero if $p$ is
$\ell$-monotone and $k\leq \ell$.

\begin{figure}[h]
\centering
\includegraphics[width=0.45\textwidth]{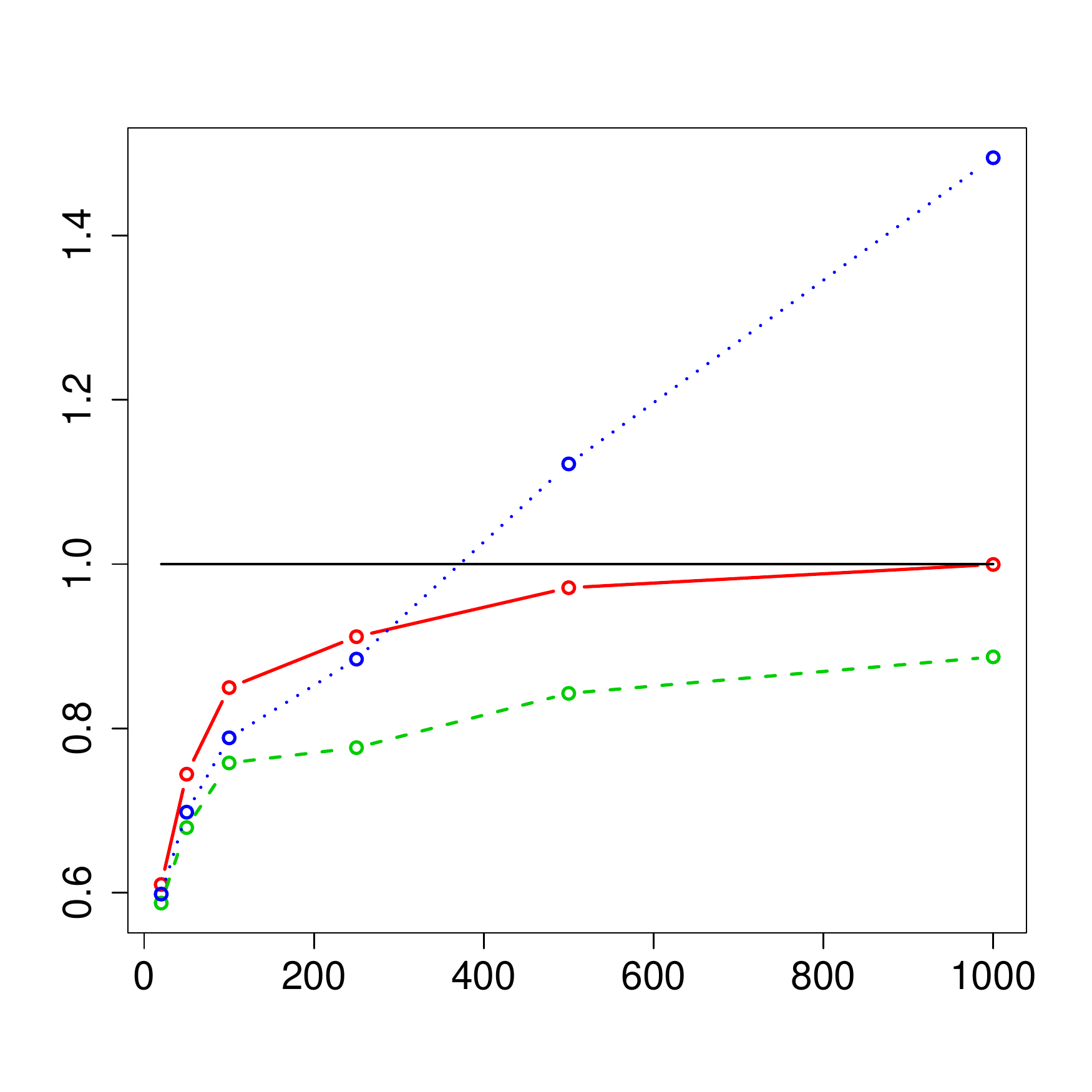}
\caption{Poisson distribution with parameter $\lambda=0.35$: ratio between the $l_2$ loss 
  $\widehat{p}^{k}_{n}$ and the $l_2$ loss of 
  $\widetilde{p}_{n}$ versus  the sample size $n$: for $k=2$ in
  "\textbf{\textcolor{red}{---}}, $k=3$ in
  "\textbf{\textcolor{green}{- -}}", and $k=4$ in
  "\textbf{\textcolor{blue}{...}}.} 
\label{fig7}
\end{figure}

The results for the Poisson distribution are similar to those obtained
for the spline distributions except that the asymptotic is
achieved for smaller values of the sample size $n$. Only the case
$\lambda=0.35$, where the corresponding Poisson distribution is 3-monotone, is
presented in Figure~\ref{fig7}. It appears that when $k=2$ the ratio
of $l_{2}$-losses tends to one, when $k=3$ it tends to a value close
to 0.9, and when $k=4$ it tends to infinity.

\FloatBarrier
\FloatBarrier

Finally we compare the $l_{2}$-losses for the estimators
$\widehat{p}^{k}_{n}$, $\widehat{p}^{* k}_{n}$ and
$\widetilde{p}_{n}$ for $k=3$ and $k=4$ (recall that for $k=2$, $\widehat{p}^{* k}_{n}=\widehat{p}^{k}_{n}$). The ratios 
$l_{2}(p,\widehat{p}_{n}^{* k})/l_{2}(p,\widetilde{p}_{n})$
behave similarly to the ratios
$l_{2}(p,\widehat{p}_{n}^{k})/l_{2}(p,\widetilde{p}_{n})$ (not shown).

Next we compare the values of the $l_2$ losses for $\widehat{p}^{* k}_{n}$ and $\widehat{p}^{k}_{n}$. When we consider the
spline distributions $Q_{j}^{\ell}$ with $l=2$ and $l=3$, the difference between the $l_2$ losses are not significant (they are smaller than $2$-times their empirical standard-error calculated on the basis of $1000$ simulations). 
When $l$ increases, the distribution $p$ is more hollow and it appears that $l_2(p,\widehat{p}^{* k}_{n})$ is greater than $l_2(p,\widehat{p}^{k}_{n})$, see Table \ref{l2.tab}. 
  
%\FloatBarrier
\begin{table}
\centering
\begin{tabular}{|c|c|c|c|}
\hline
 & $n=20$ & $n=100$ & $n=1000$ \\
\hline
$p=Q^4_{10}$  & \textcolor{green}{-} & \textcolor{green}{-} & \textcolor{green}{-}\\
             & \textcolor{blue}{-} & \textcolor{blue}{-}& \textcolor{blue}{0.06} \\
\hline
$p=Q^{10}_{10}$ & \textcolor{green}{0.89} & \textcolor{green}{0.13} & \textcolor{green}{0.02}\\
             & \textcolor{blue}{0.92} & \textcolor{blue}{0.24}&\textcolor{blue}{0.01}\\
\hline
\end{tabular}
\caption{Spline distributions : difference ($\times 1000$)
between the $l_2$-loss of $\widehat{p}_{n}^{*k}$ and the $l_2$-loss of $\widehat{p}_{n}^{k}$, 
  for different values of $n$, for  $k=3$  in green and $k=4$ in
  blue. The symbol "-" is for non-significant result. }
\label{l2.tab}
\end{table}
%\FloatBarrier

\subsubsection{Estimators comparison based on the Hellinger loss}

Let us now consider the Hellinger loss defined, for any
estimator $\widehat{q}$, as $H(p,\hat{q}) =
E\left(\| \sqrt{p} - \sqrt{\widehat{q}}\|_{2}^{2}\right)$. 

\begin{figure}[h]
\centering
\includegraphics[width=0.45\textwidth]{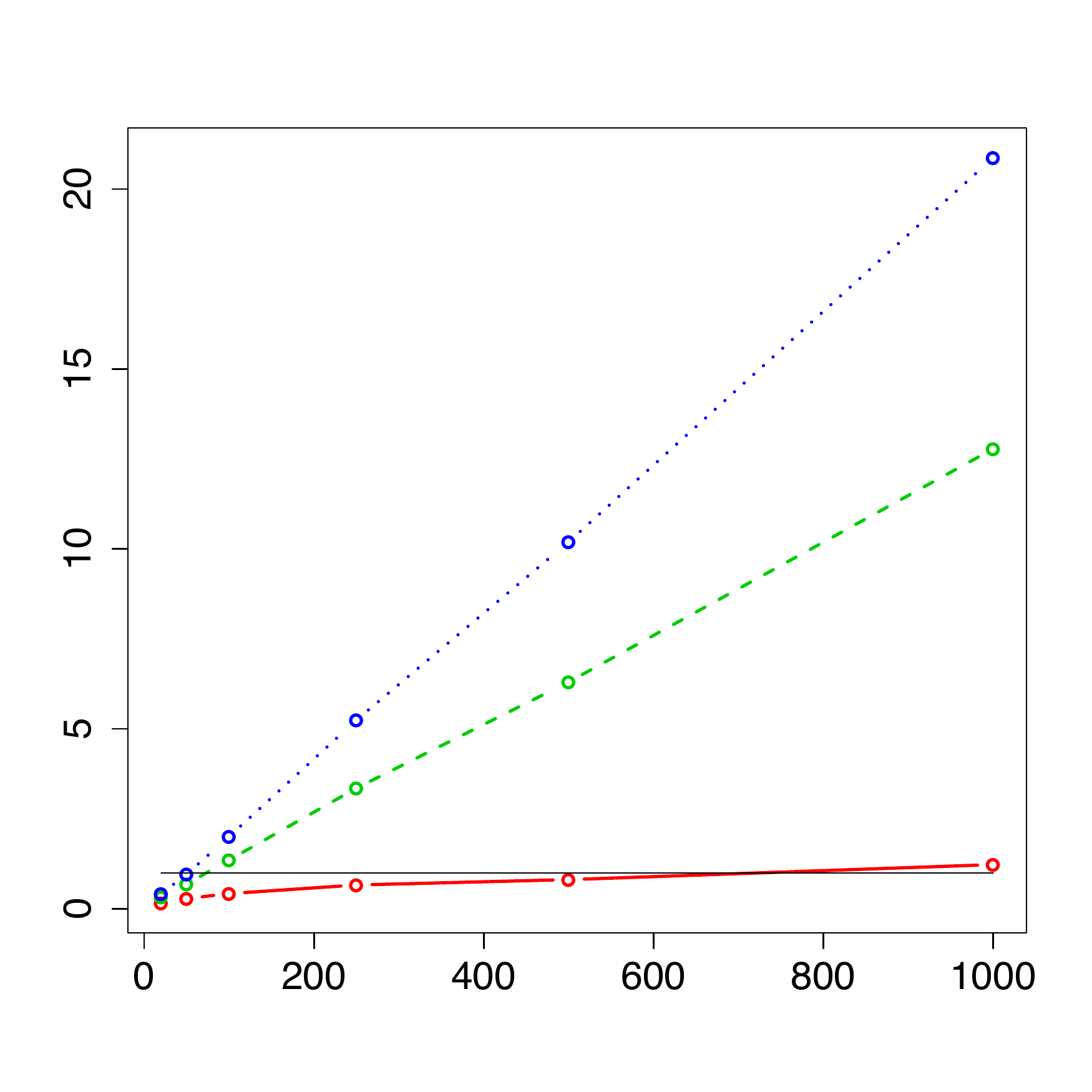}
\label{fig2a}
\caption{Triangular distribution $Q_{j}^{2}$ : ratio between the
  Hellinger loss of $\widehat{p}^{k}_{n}$ and the Hellinger loss of
  $\widetilde{p}_{n}$ versus the sample size $n$: for $k=2$ in
  "\textbf{\textcolor{red}{---}}, $k=3$ in
  "\textbf{\textcolor{green}{- -}}", and $k=4$ in "\textbf{\textcolor{blue}{...}}".}
\label{Hell.Qj2.fig}
\end{figure}

The results for the spline distributions $Q_{j}^{\ell}$ are similar to
those obtained for the $l_{2}$-loss, except that the ratios
$H(p,\widehat{p}_{n}^{k})/H(p,\widetilde{p}_{n})$ are not necessary
smaller than 1 when $k \leq \ell$, see Figure~\ref{Hell.Qj2.fig} for the
Triangular distribution $Q_{j}^{2}$.

\FloatBarrier
\FloatBarrier

\begin{figure}[!h]
\centering
\begin{subfigure}[b]{0.45\textwidth}
\includegraphics[width=0.95\textwidth]{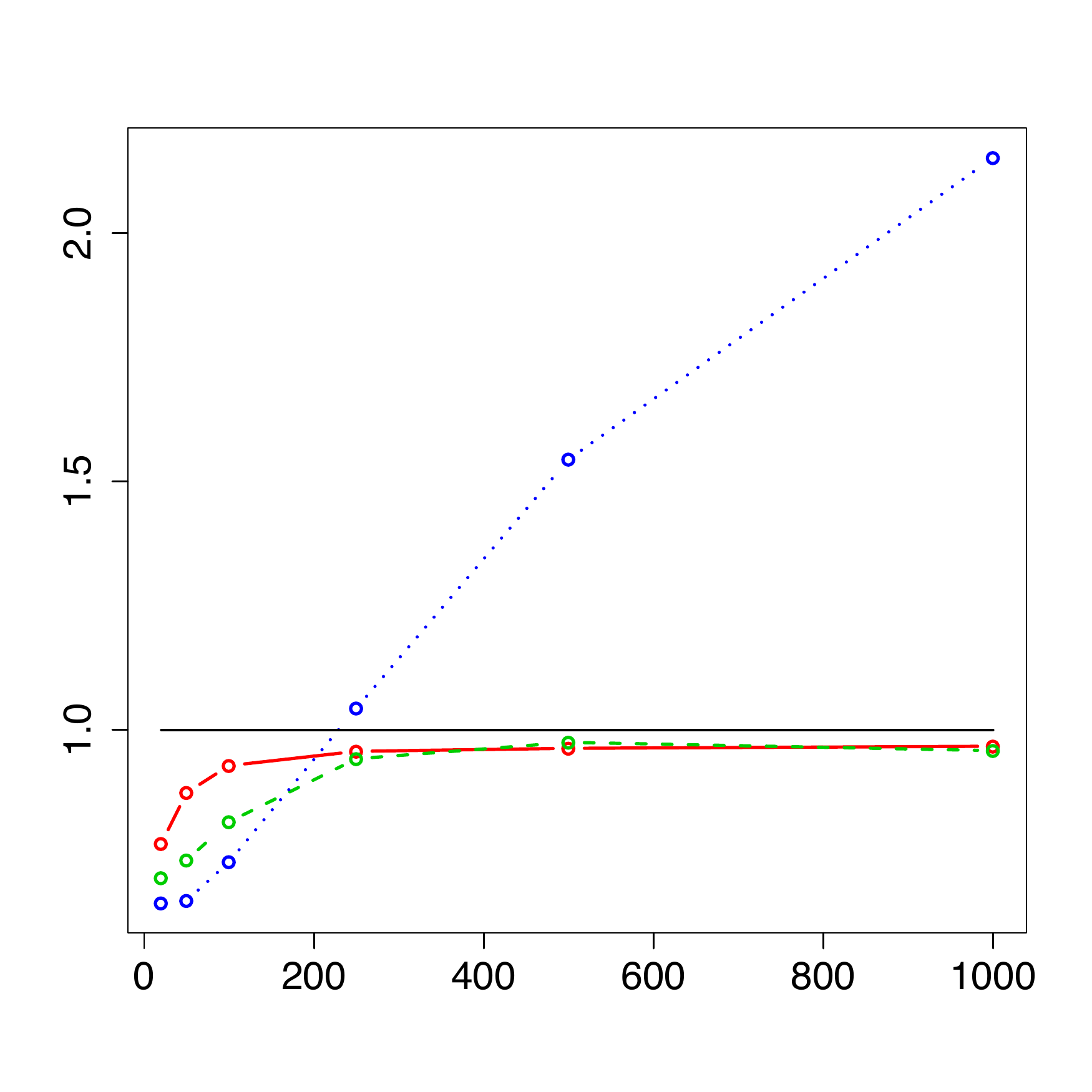}
\caption{$\lambda=0.3; \ell=4$.}
\label{fig8a}
\end{subfigure}
\hfill
\begin{subfigure}[b]{0.45\textwidth}
\includegraphics[width=0.95\textwidth]{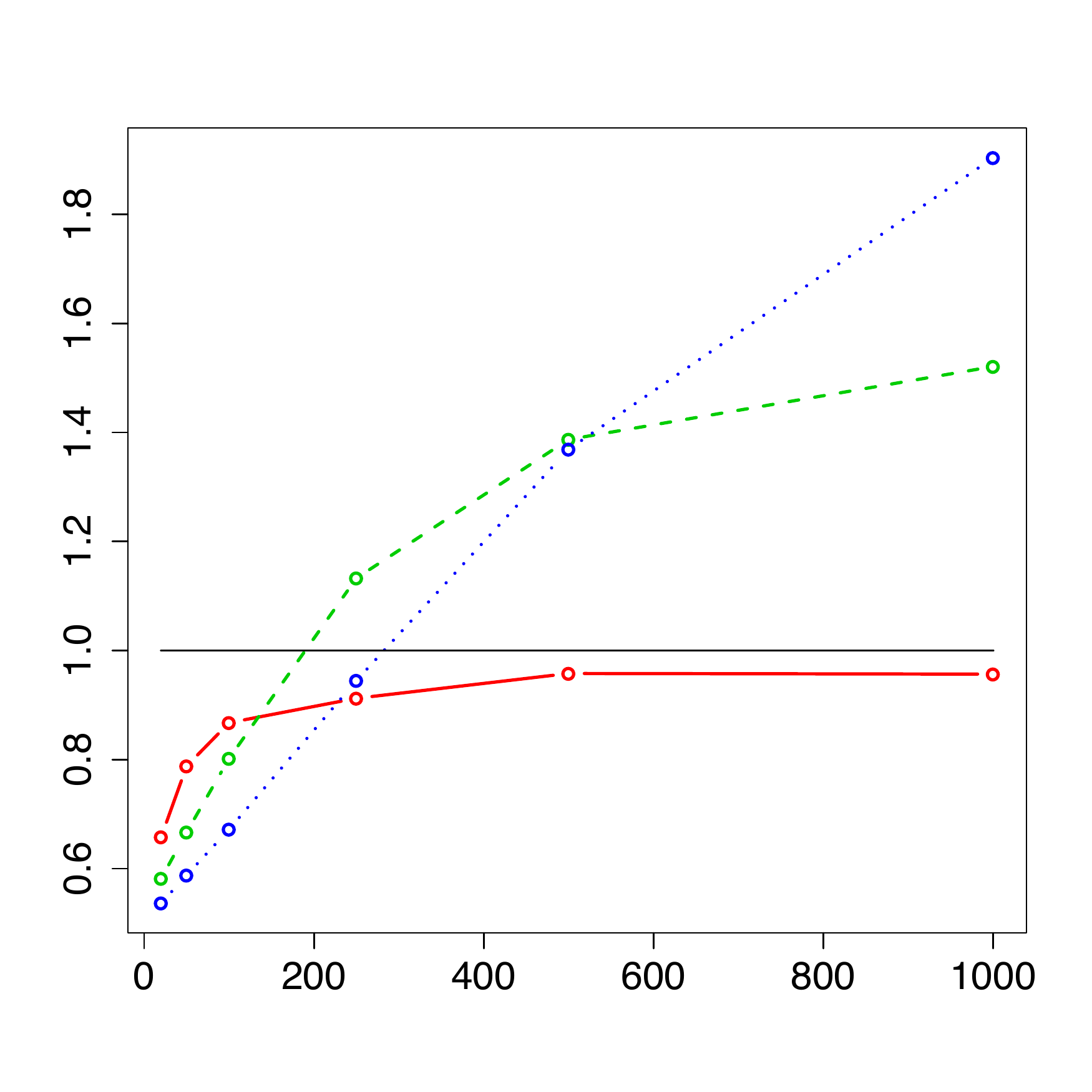}
\caption{$\lambda=0.35; \ell=3$.}
\label{fig8b}
\end{subfigure}
\begin{subfigure}[b]{0.45\textwidth}
\includegraphics[width=0.95\textwidth]{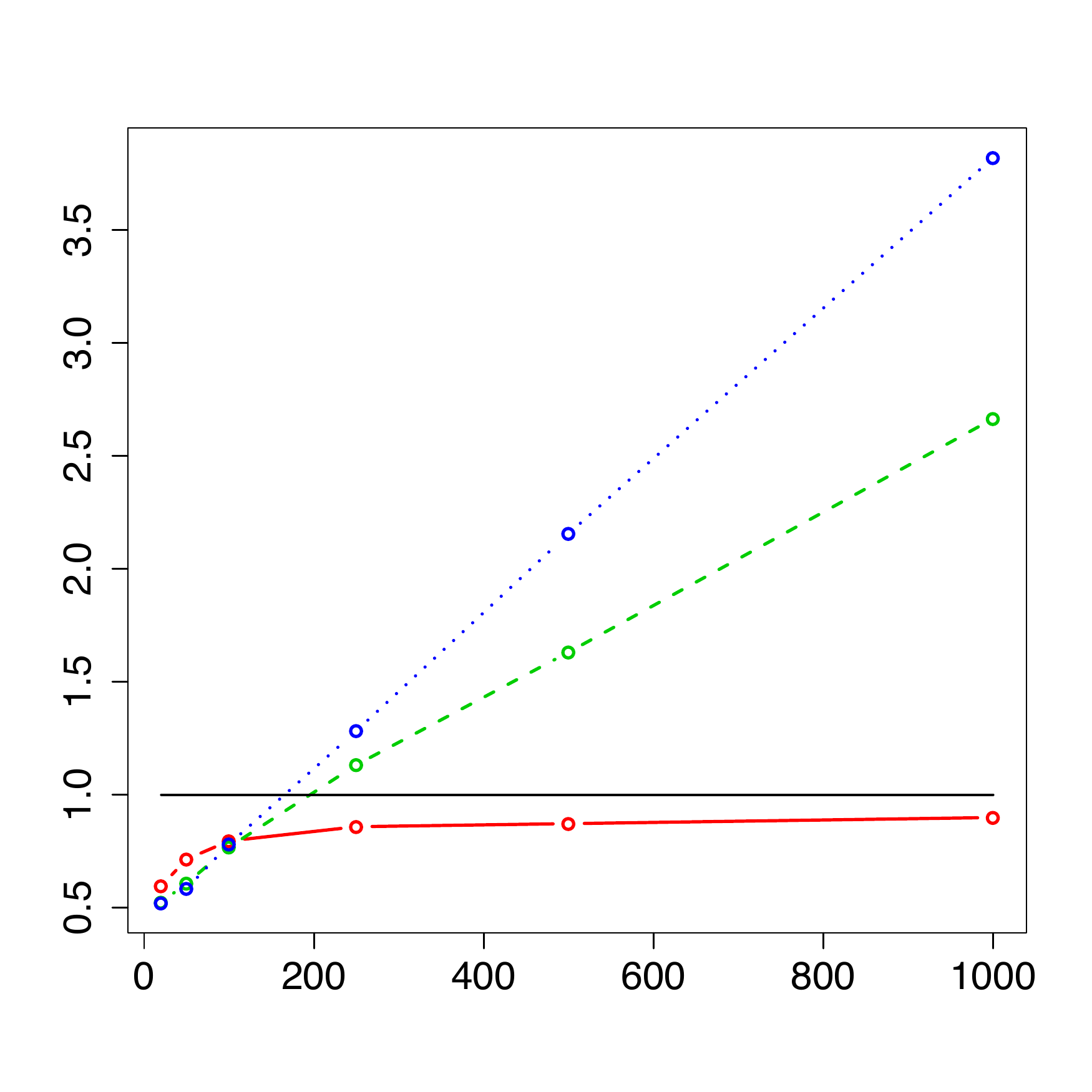}
\caption{$\lambda=0.45; \ell=2$.}
\label{fig8c}
\end{subfigure}
\hfill
\begin{subfigure}[b]{0.45\textwidth}
\includegraphics[width=0.95\textwidth]{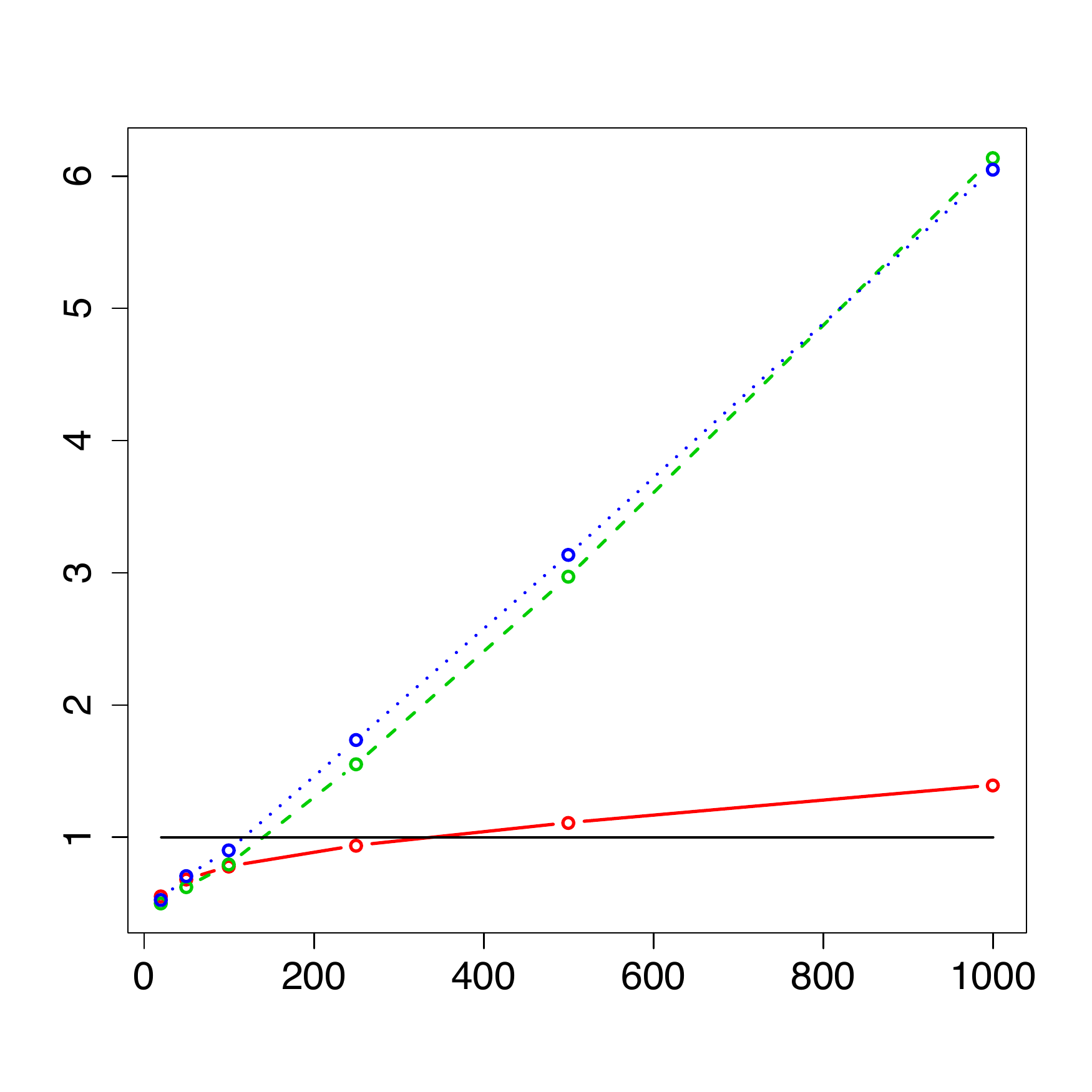}
\caption{$\lambda=2-\sqrt{2}; \ell=2$.}
\label{fig8d}
\end{subfigure}
\caption{Poisson distributions: ratio between the Hellinger loss of
  $\widehat{p}_{n}^{k}$ and the Hellinger loss of $\widetilde{p}_{n}$
  versus the sample size $n$: for $k=2$ in
  "\textbf{\textcolor{red}{---}}, $k=3$ in
  "\textbf{\textcolor{green}{- -}}", and $k=4$ in
  "\textbf{\textcolor{blue}{...}}". Each subfigure corresponds to the
  results obtained with $p={\mathcal P}(\lambda)$ for $\lambda \in
  \left\{0.3, 0.35, 0.45, 2-\sqrt{2}\right\}$. The degree of monotony of these
  distributions is given by $\ell$.}
\label{fig8}
\end{figure}

In the case of the Poisson distributions the differences between the
$l_{2}$-loss and the Hellinger loss are more obvious. As it is
illustrated by Figure~\ref{fig8}, if $\ell$ the degree of monotony of
$p$ is strictly greater than $k$, then the ratio is smaller than 1
(see case (a) with $k=2,3$ and case (b) with $k=2$). If $k=\ell$,
then $H(p, \widetilde{p}^{k}_{n})$ is smaller than $H(p, \widehat{p}_{n})$ 
if the distribution $p$ is ``$\ell$-monotone enough'', that is to say if
the parameter $\lambda$ of the Poisson distribution is such that
$\lambda_{\ell} - \lambda$  is large enough, where $\lambda_{\ell}$
has been defined in Property~\ref{Poiss.prop}, see for example cases
(c) and (d) with $k=2$, where $\lambda_{2} = 2-\sqrt{2}$.

\FloatBarrier
\FloatBarrier

\subsection{Some characteristics of interest}

We consider the estimation of some charactéeristics that may be of
interest as the entropy, the variance and the probability at 0. For
each of these characteristics denoted $L(p)$, we measure the
performance in terms of the root mean squared error of prediction
calculated as follows:
\begin{equation*}
 \mathrm{ RMSEP} = \sqrt{\mathrm{ BIAS}^{2} + \mathrm{ SE}^{2}},
\end{equation*}
where $\mathrm{ BIAS}$ and $\mathrm{SE}$ are the estimated bias and
standard-error of the estimator based on the simulations. Let
$\widehat{L}$ be an estimator of $L(p)$, then $\mathrm{ BIAS}=\widehat{L}_{\cdot} - L$, where $\widehat{L}_{\cdot} =
\sum_{s}\widehat{L}_{s}/1000$ with $\widehat{L}_{s}$ being the
estimate of $L(p)$ at simulation $s$, and  $\mathrm{ SE}^{2} =
\sum_{s}(\widehat{L}_{s}- \widehat{L}_{\cdot}
)^{2}/1000$. 

\subsubsection{Entropy}

The entropy is defined as 
\begin{eqnarray*}
\text{Ent}(f)=\sum_{i=0}^{\infty}f(i)\log(f(i)).
\end{eqnarray*}

\begin{figure}[!h]
\centering
\begin{subfigure}[b]{0.45\textwidth}
\includegraphics[width=0.95\textwidth]{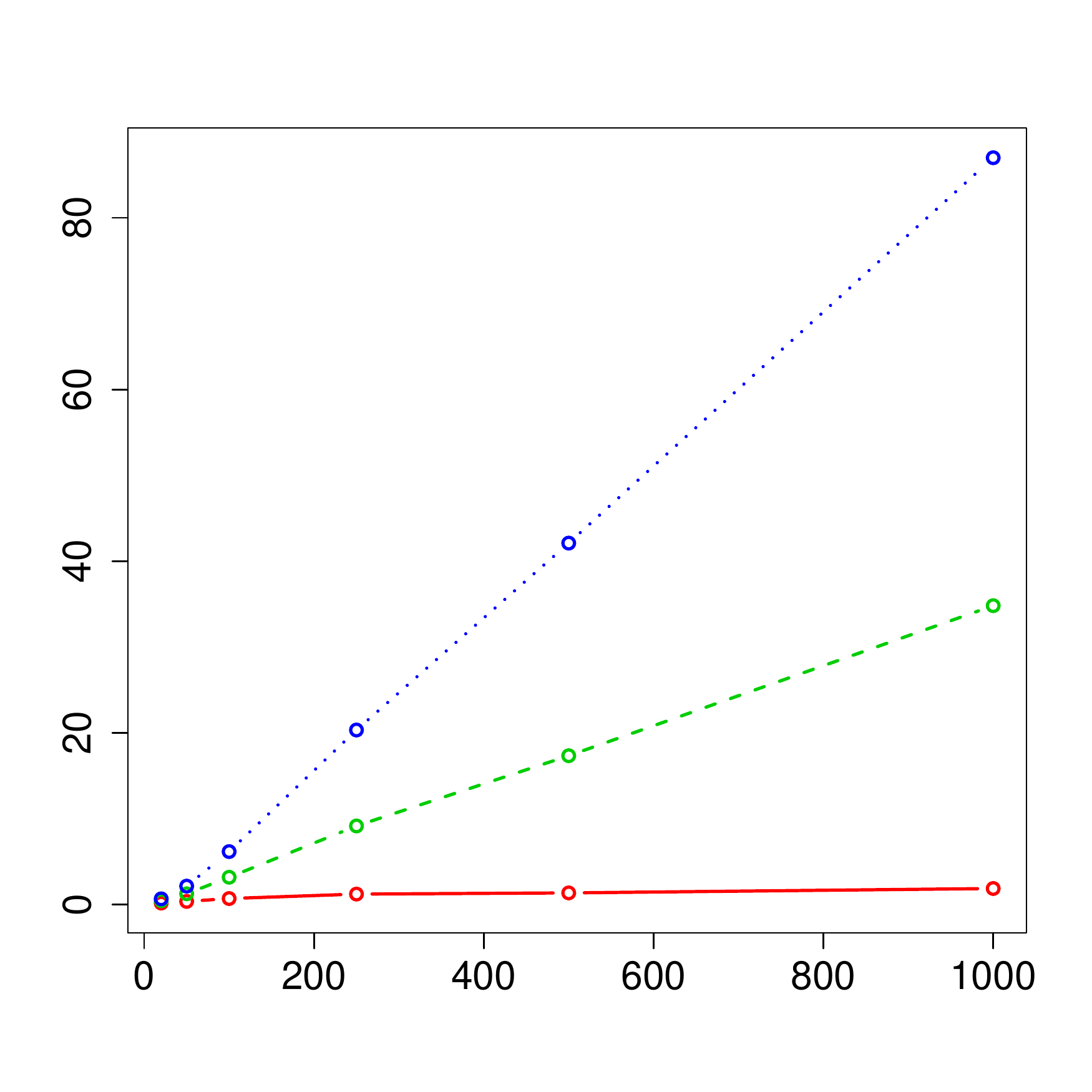}
\caption{$p=Q^2_{10}$.}
\label{fig4a}
\end{subfigure}
\hfill
\begin{subfigure}[b]{0.45\textwidth}
\includegraphics[width=0.95\textwidth]{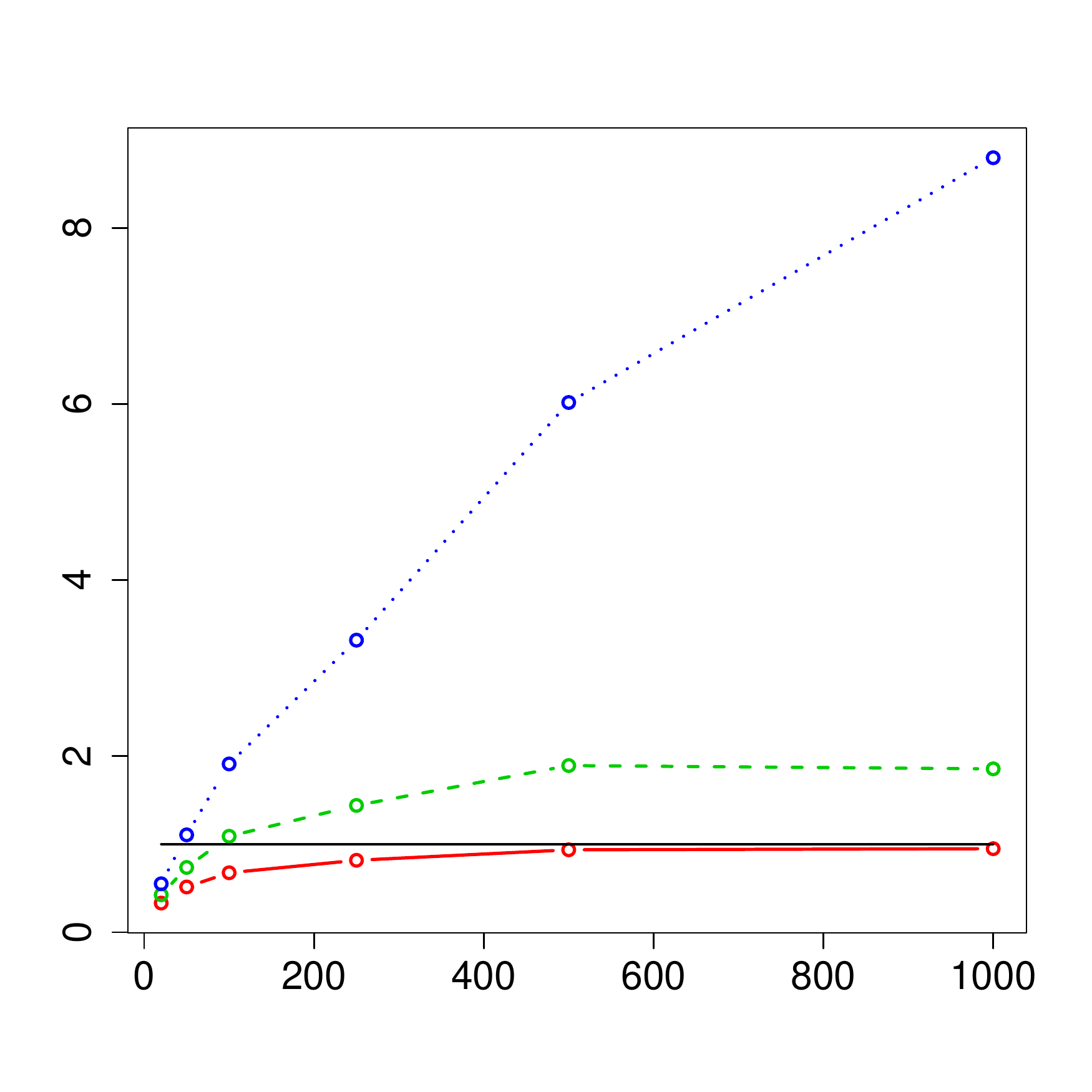}
\caption{$p=Q^3_{10}$.}
\label{fig4b}
\end{subfigure}
\begin{subfigure}[b]{0.45\textwidth}
\includegraphics[width=0.95\textwidth]{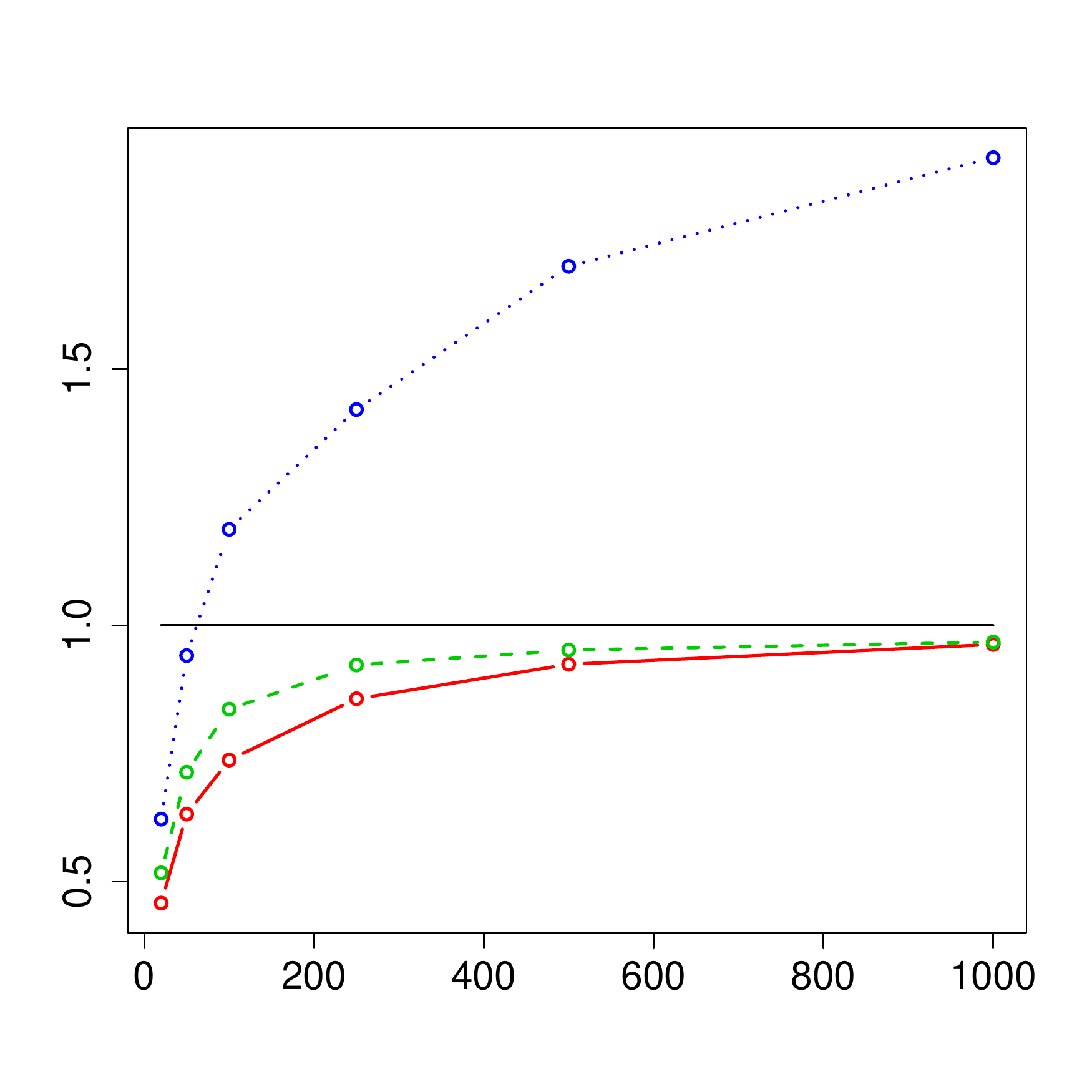}
\caption{$p=Q^4_{10}$.}
\label{fig4c}
\end{subfigure}
\hfill
\begin{subfigure}[b]{0.45\textwidth}
\includegraphics[width=0.95\textwidth]{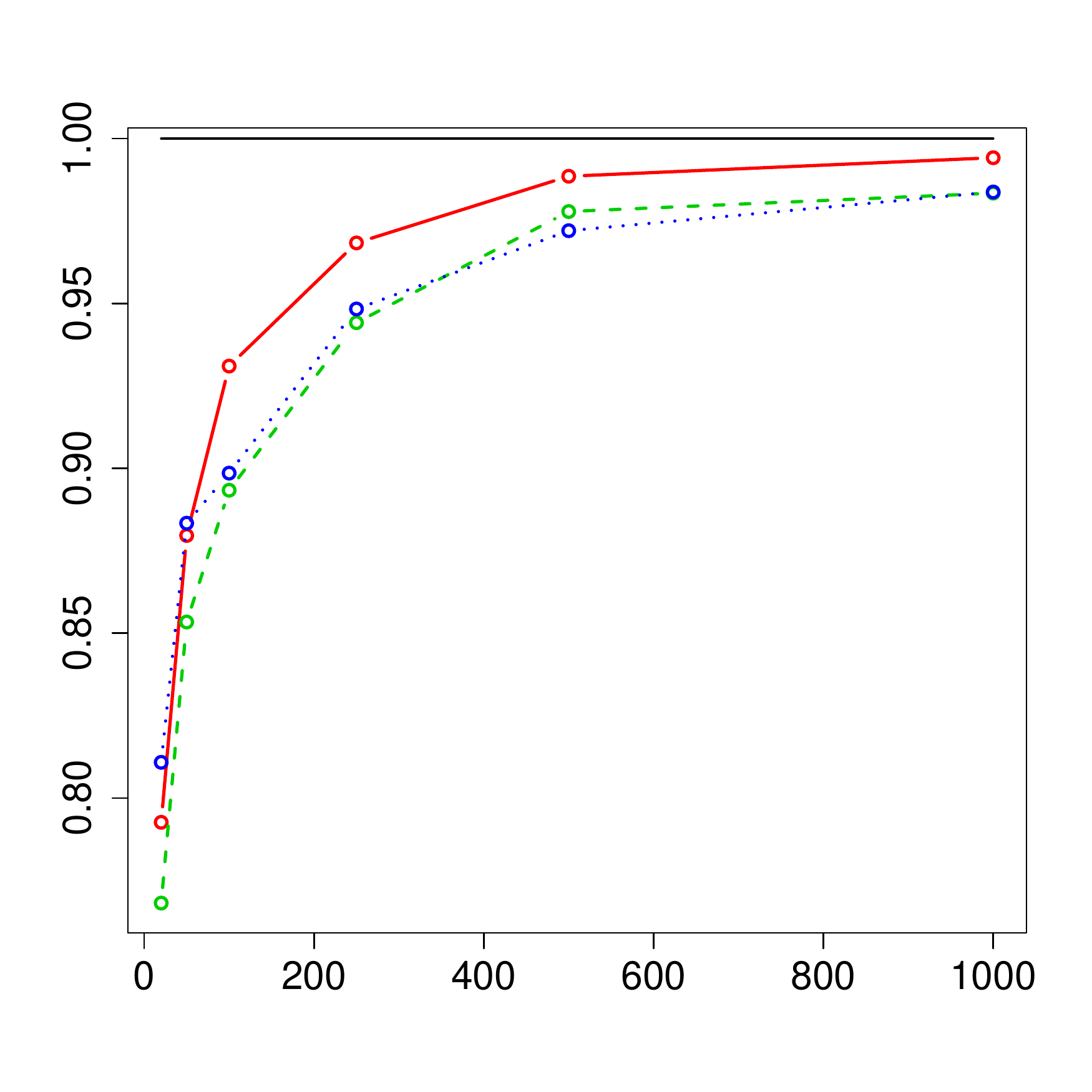}
\caption{$p=Q^{10}_{10}$.}
\label{fig4d}
\end{subfigure}
\caption{Spline distributions: ratio between the $\mathrm{RMSEP}$ of
  $\text{Ent}(\widehat{p}^{k}_{n})$ and the $\mathrm{RMSEP}$ of
  $\text{Ent}(\widetilde{p}_{n})$ versus the sample size $n$: for $k=2$ in "\textbf{\textcolor{red}{---}}",
  $k=3$ in "\textbf{\textcolor{green}{- -}}", $k=4$ in
  "\textbf{\textcolor{blue}{...}}". Each subfigure
  corresponds to the results  obtained with  $p = Q_{j}^{\ell}$, for
  $\ell\in\{2, 3, 4, 10\}$.}
\label{fig4}
\end{figure}

We compare the estimators $\text{Ent}(\widehat{p}^{k}_{n})$ and
$\text{Ent}(\widetilde{p}_{n})$ by the ratio of their $\mathrm{
  RMSEP}$. 
The results differ according to the family of distributions. 
For the spline distributions $Q_{j}^{\ell}$, see Figure~\ref{fig4}, it appears that if $k <
\ell$, then $\text{Ent}(\widehat{p}^{k}_{n})$ has smaller $\mathrm{
  RMSEP}$ than $\text{Ent}(\widetilde{p}_{n})$. However, when $k=l$, the ratio of the $\mathrm{
  RMSEP}$'s increases and  reaches an asymptote greater than
1. For example, in Figure~\ref{fig4}, case (b) with $k=3$, the ratio
tends to $0.96$, in case (c) with $k=4$, the ratio tends
to $1.93$. In fact, if we consider the space of $\ell$-monotone distributions
with maximum support $j$, the distribution $Q_{j}^{\ell}$ may appear
as a ``limiting  case'' in this space, in  that it admits only one $\ell$-knot in
$j$. It seems that for these $Q_{j}^{\ell}$  distributions,
the projection on the space of 
$\ell-1$-monotone discrete probabilities give better results than on
the space of $\ell$-monotone discrete probabilities.

%\FloatBarrier
\FloatBarrier

For the Poisson distributions, see Figure~\ref{fig15}, when $n$ is
small, the estimator
based on the emprirical distribution, $\text{Ent}(\widetilde{p}_{n})$,
has a smaller RMSEP than $\text{Ent}(\widehat{p}^{k}_{n})$. When $n$
is large the RMSEP ratio tend to one if $k \leq \ell$, and tend to
infinity if $k>\ell$.

\FloatBarrier
\begin{figure}[!h]
\centering
\begin{subfigure}[b]{0.45\textwidth}
\includegraphics[width=0.95\textwidth]{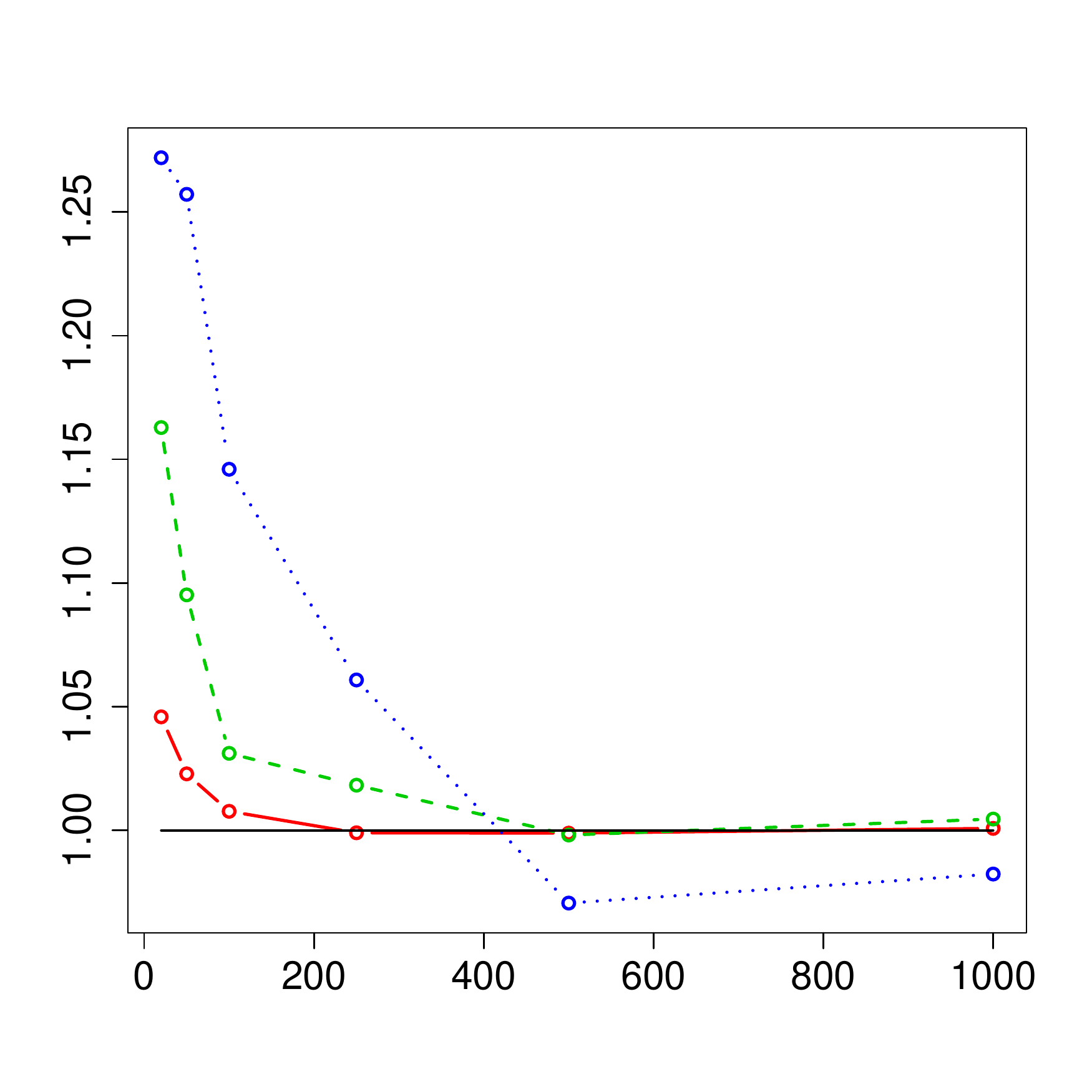}
\caption{$\lambda=0.3; \ell=4$.}
\label{fig15a}
\end{subfigure}
\hfill
\begin{subfigure}[b]{0.45\textwidth}
\includegraphics[width=0.95\textwidth]{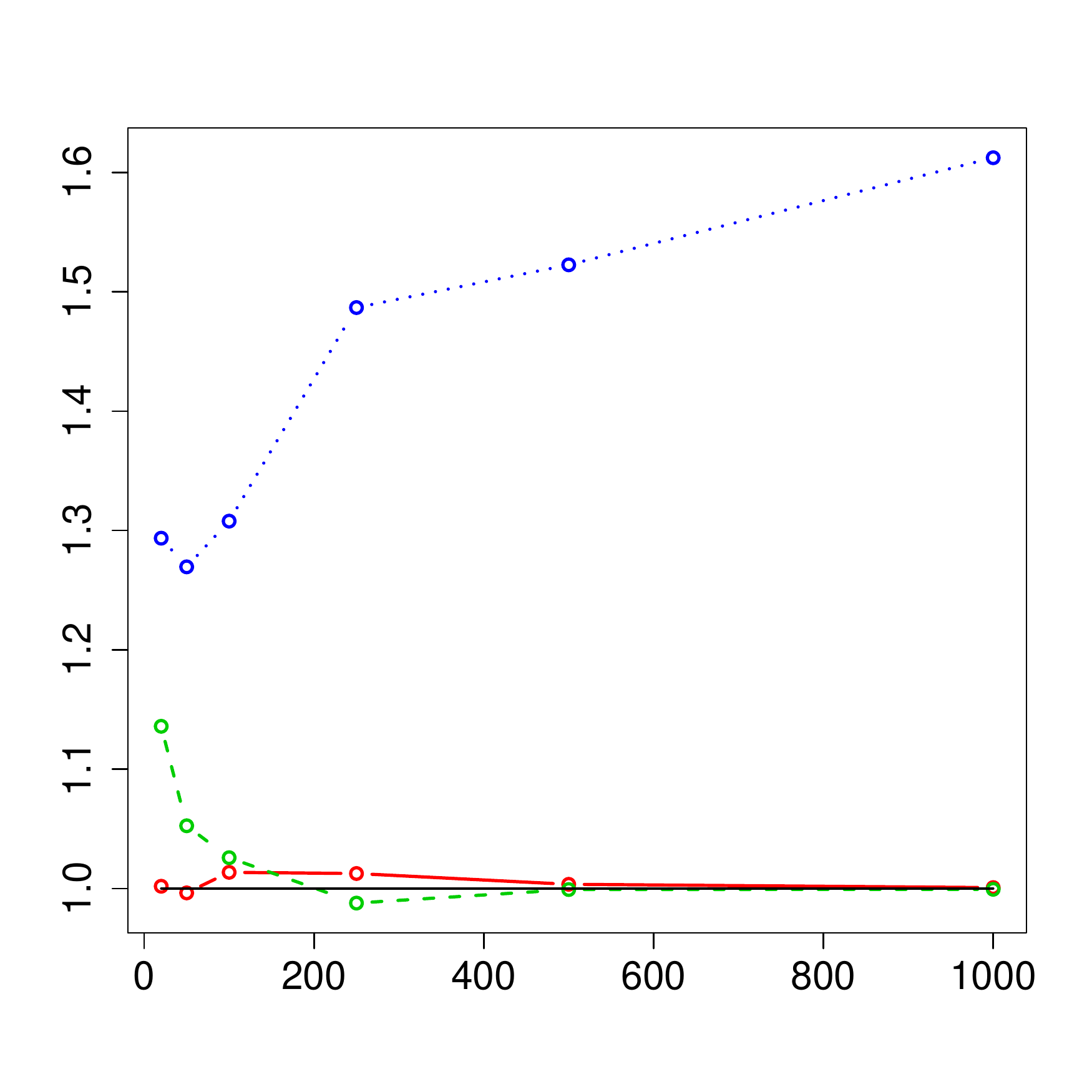}
\caption{$\lambda=0.35; \ell=3$.}
\label{fig15b}
\end{subfigure}
\caption{Poisson distributions: ratio between the $\mathrm{RMSEP}$ of
  $\text{Ent}(\widehat{p}^{k}_{n})$ and the $\mathrm{RMSEP}$ of
  $\text{Ent}(\widetilde{p}_{n})$ versus the sample size $n$: for $k=2$ in "\textbf{\textcolor{red}{---}}",
  $k=3$ in "\textbf{\textcolor{green}{- -}}", $k=4$ in
  "\textbf{\textcolor{blue}{...}}". Each subfigure
  corresponds to the results  obtained with $p={\mathcal P}(\lambda)$
  with $\lambda \in \{0.3, 0.35\}$.}
\label{fig15}
\end{figure}
\FloatBarrier

\subsubsection{Probability mass in 0.}

We compare the performances of $\widehat{p}^{k}_{n}(0)$
and $\widetilde{p}_{n}(0)$  by comparing the corresponding renormalized SE and BIAS. 

The results for the spline distributions are presented in Table \ref{p0eerror.tab}.

When $k\leqslant l$, $\widehat{p}^{k}_{n}(0)$ has smaller SE than $\widetilde{p}_{n}(0)$. Its bias is greater in absolute value and always negative, but the RMSEP stays smaller.
For each $k$, the variations of $\sqrt{n}SE/p(0)$ versus $n$ are very small and tend to stabilize around a value that increases with $l-k$.

When $k>l$, $\widehat{p}^{k}_{n}(0)$ keeps a smaller RMSEP than $\widetilde{p}_{n}(0)$ for small $n$. But, when $n$ increases the absolute bias as well as the standard error increase.

The results for the Poisson distributions are similar and omitted.

\FloatBarrier
\begin{table}
\centering
\begin{tabular}{|c|c|c|c|c|c|c|c|c|c|}
\hline
 & \multicolumn{3}{c|}{$n=20$} & \multicolumn{3}{c|}{$n=100$} &  \multicolumn{3}{c|}{$n=1000$} \\
\hline 
&SE & BIAS & RMSEP& SE & BIAS & RMSEP & SE & BIAS & RMSEP\\
\hline
$p=Q^2_{10}$ & 2.25 & 7e-4  & 2.25& 2.234 & 0.002 & 2.234 & 2.284 & 0.017  & 2.284\\
             & \textcolor{red}{1.800} & \textcolor{red}{0.181} & \textcolor{red}{1.809}  & \textcolor{red}{1.819} & \textcolor{red}{0.170} & \textcolor{red}{1.82}  & \textcolor{red}{1.745} & \textcolor{red}{0.162} & \textcolor{red}{1.752}\\
            & \textcolor{green}{1.757} & \textcolor{green}{0.157} & \textcolor{green}{1.764} &  \textcolor{green}{1.783} & \textcolor{green}{0.188} & \textcolor{green}{1.792}  & \textcolor{green}{2.231} & \textcolor{green}{0.334} & \textcolor{green}{2.255}\\
             & \textcolor{blue}{1.742} & \textcolor{blue}{0.155} & \textcolor{blue}{1.748} & \textcolor{blue}{1.780} & \textcolor{blue}{0.196} & \textcolor{blue}{1.790}& \textcolor{blue}{2.622} & \textcolor{blue}{0.408} & \textcolor{blue}{2.653}  \\
\hline
$p=Q^4_{10}$ & 1.634 & 0.008 & 1.634 & 1.601 & 0.013 & 1.601 & 1.626 & 0.006 & 1.626\\
             & \textcolor{red}{1.362} & \textcolor{red}{0.143} & \textcolor{red}{1.369} & \textcolor{red}{1.389} & \textcolor{red}{0.120} & \textcolor{red}{1.394}& \textcolor{red}{1.488} & \textcolor{red}{0.052} & \textcolor{red}{1.489}   \\
             & \textcolor{green}{1.354} & \textcolor{green}{0.137} & \textcolor{green}{1.361} & \textcolor{green}{1.372} & \textcolor{green}{0.132} & \textcolor{green}{1.378}& \textcolor{green}{1.439} & \textcolor{green}{0.088} & \textcolor{green}{1.442}   \\
             & \textcolor{blue}{1.340} & \textcolor{blue}{0.135} & \textcolor{blue}{1.347}& \textcolor{blue}{1.353} & \textcolor{blue}{0.136} & \textcolor{green}{1.359}& \textcolor{blue}{1.362} & \textcolor{blue}{0.109} & \textcolor{blue}{1.366} \\
\hline
$p=Q^{10}_{10}$ & 1.010 & 2e-4  & 1.010 & 0.98 & 6e-4 & 0.98 & 0.984 & 0.006  & 0.984 \\
                & \textcolor{red}{0.884} & \textcolor{red}{0.058} & \textcolor{red}{0.886}& \textcolor{red}{0.934} & \textcolor{red}{0.022} & \textcolor{red}{0.934}& \textcolor{red}{0.982} & \textcolor{red}{0.006} & \textcolor{red}{0.982}  \\
         & \textcolor{green}{0.886} & \textcolor{green}{0.057} & \textcolor{green}{0.888}& \textcolor{green}{0.919} & \textcolor{green}{0.039} & \textcolor{green}{0.920}& \textcolor{green}{0.957} & \textcolor{green}{0.009} & \textcolor{green}{0.957}\\
               & \textcolor{blue}{0.887} & \textcolor{blue}{0.053} & \textcolor{blue}{0.889}& \textcolor{blue}{0.921} & \textcolor{blue}{0.042} & \textcolor{blue}{0.922}& \textcolor{blue}{0.940} & \textcolor{blue}{0.018} & \textcolor{blue}{0.940} \\
\hline
\end{tabular}
\caption{Spline distributions: $\sqrt{n} \mathrm{SE}/p(0)$, 
  $\sqrt{n} \vert\mathrm{BIAS}\vert/p(0)$ and $\sqrt{n} \mathrm{RMSEP}/p(0)$ for $\widetilde{p}_{n}(0)$
  in black, $\widehat{p}_{n}^{k}(0)$ for $k=2$ in red, $k=3$ in
  green and $k=4$ in blue.}
\label{p0eerror.tab}
\end{table}
\FloatBarrier

\subsubsection{Variance}

We compare the estimators of the variance of $p$, denoted $\text{var}(\widehat{p}^{k}_{n})$ and
$\text{var}(\widetilde{p}_{n})$ comparing the ratio of their $\mathrm{
  RMSEP}$. 
The results are similar for the spline distributions and the Poisson's distributions and we present only the RMSEP for the spline distributions $Q^l_j$ in Figure \ref{fig15}.\\

When $k=l$, the ratio of the RMSEP tends to a constant smaller than $1$ when $n$ tends to infinity. Conversely if we are not in a good model ($k>l$) the ratio of the RMSEPs tends to infinity when $n$ tends to infinity.

When $k<l$ and $n$ large the ratio of the RMSEPs increases with $l-k$ and goes beyond $1$. For example for $k=3$ and $l=4$ the ratio of the RMSEPs is equal to $0.68$ when $n=10000$, while if $l=10$ the ratio is greater than $1$ as soon as $k\leqslant 3$ and $n\geqslant 1000$.

When $k>l$ the ratio of the RMSEPs tends to infinity when n tends to infinity.

When $n$ is small $\text{var}(\widehat{p}^{k}_{n})$ has smaller RMSEP than $\text{var}(\widetilde{p}_{n})$ whatever the value of $k$ and $l$.\\

%\FloatBarrier
\begin{figure}[!h]
\centering
\begin{subfigure}[b]{0.45\textwidth}
\includegraphics[width=0.95\textwidth]{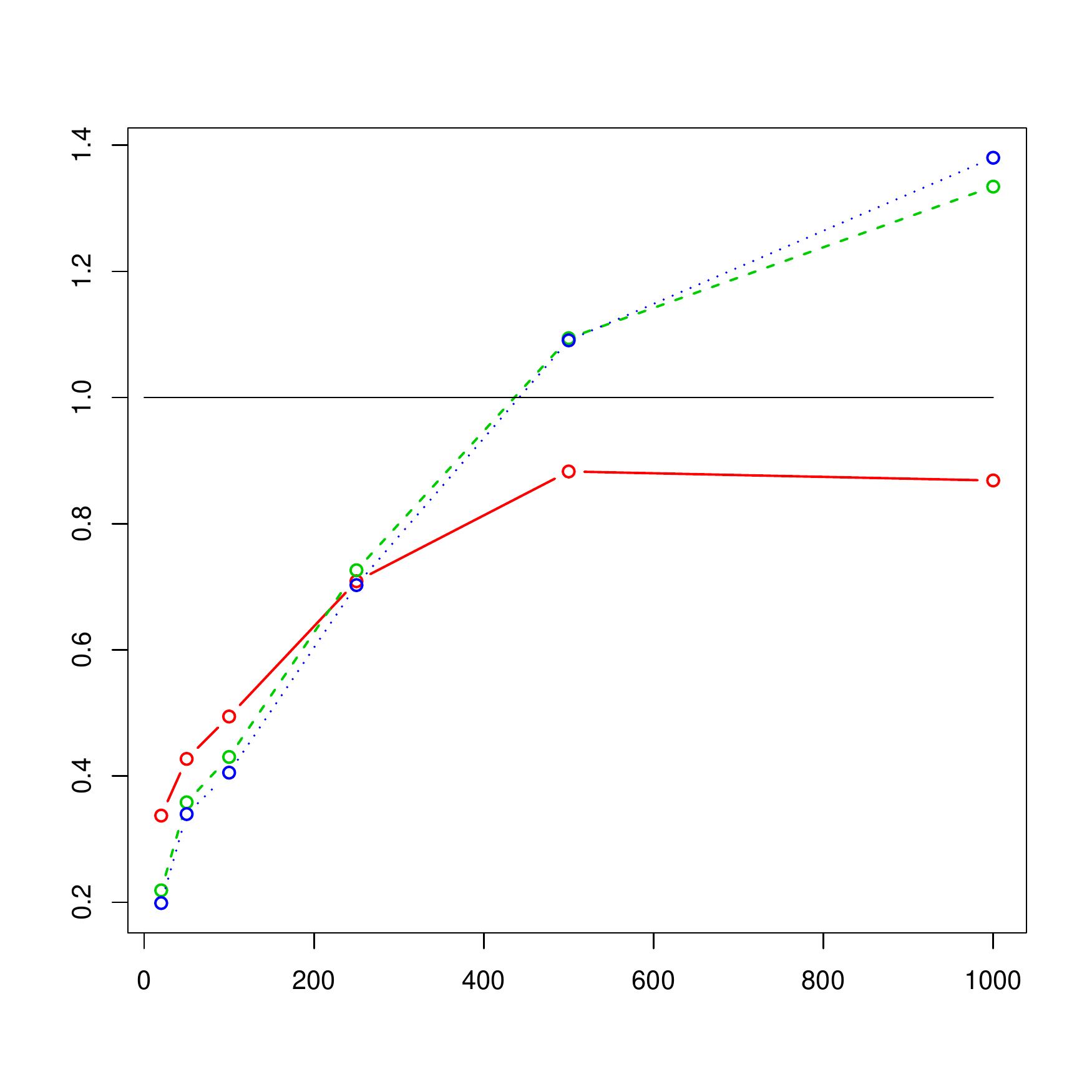}
\caption{$p=Q^2_{10}$.}
\label{fig15a}
\end{subfigure}
\hfill
\begin{subfigure}[b]{0.45\textwidth}
\includegraphics[width=0.95\textwidth]{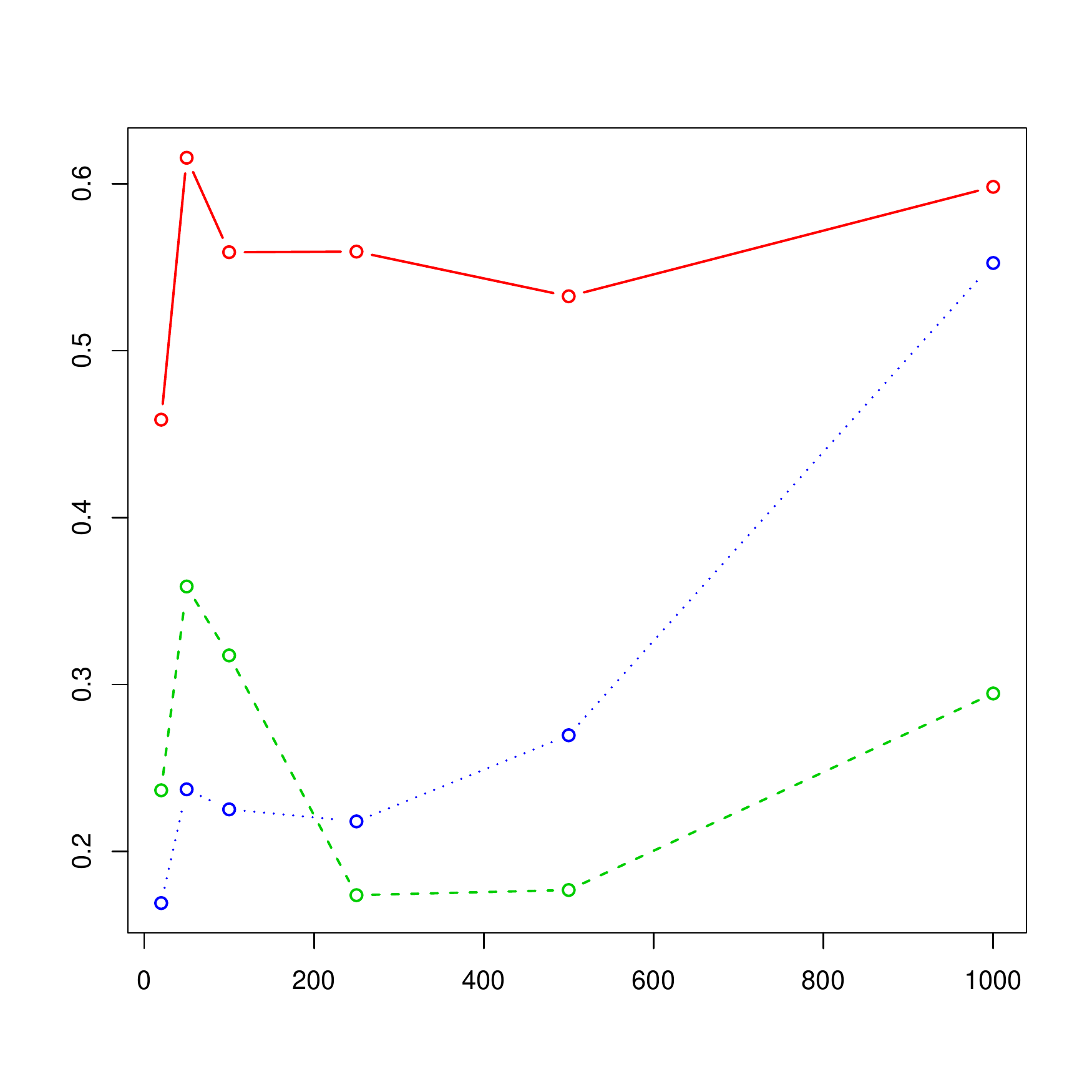}
\caption{$p=Q^3_{10}$.}
\label{fig15b}
\end{subfigure}
\begin{subfigure}[b]{0.45\textwidth}
\includegraphics[width=0.95\textwidth]{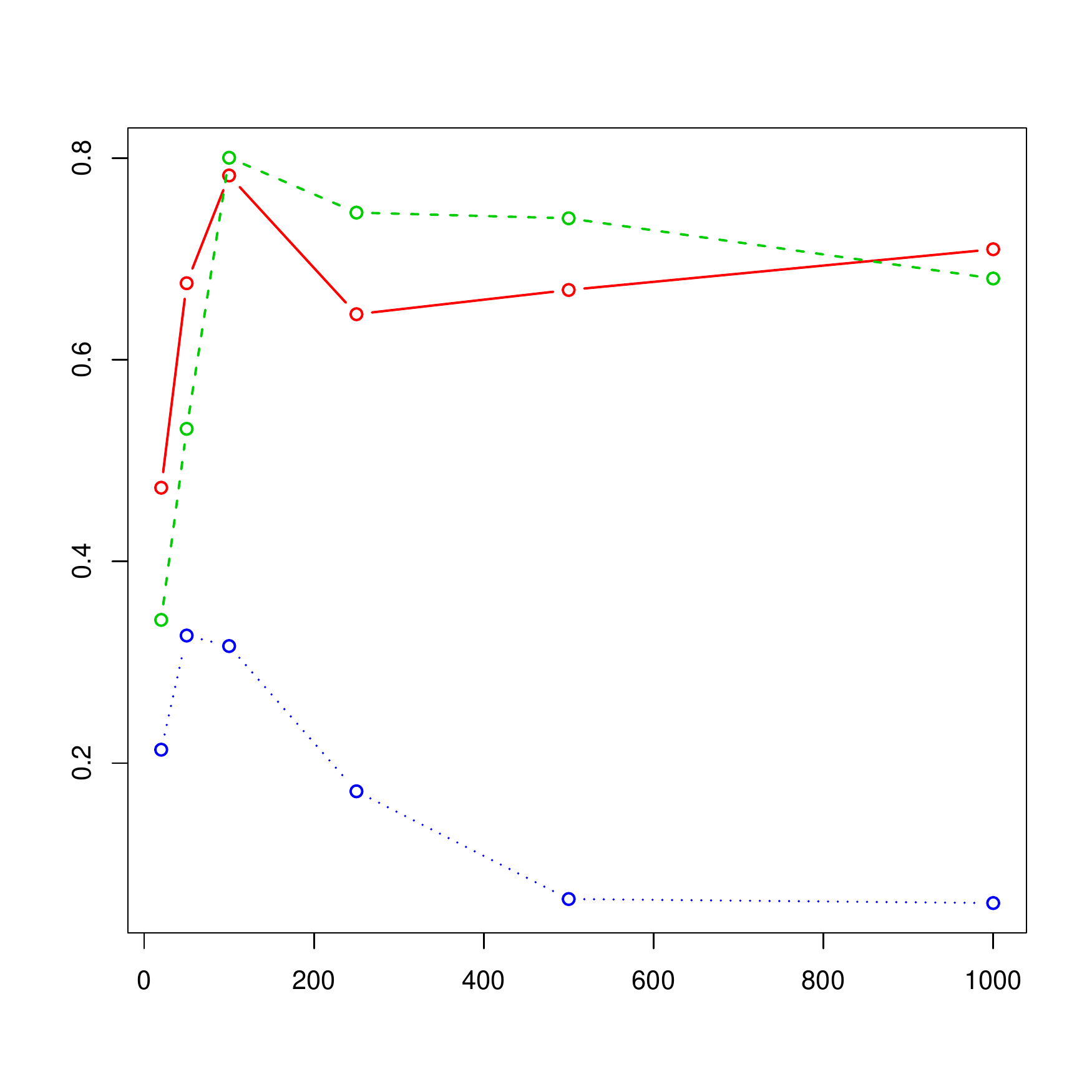}
\caption{$p=Q^4_{10}$.}
\label{fig15c}
\end{subfigure}
\hfill
\begin{subfigure}[b]{0.45\textwidth}
\includegraphics[width=0.95\textwidth]{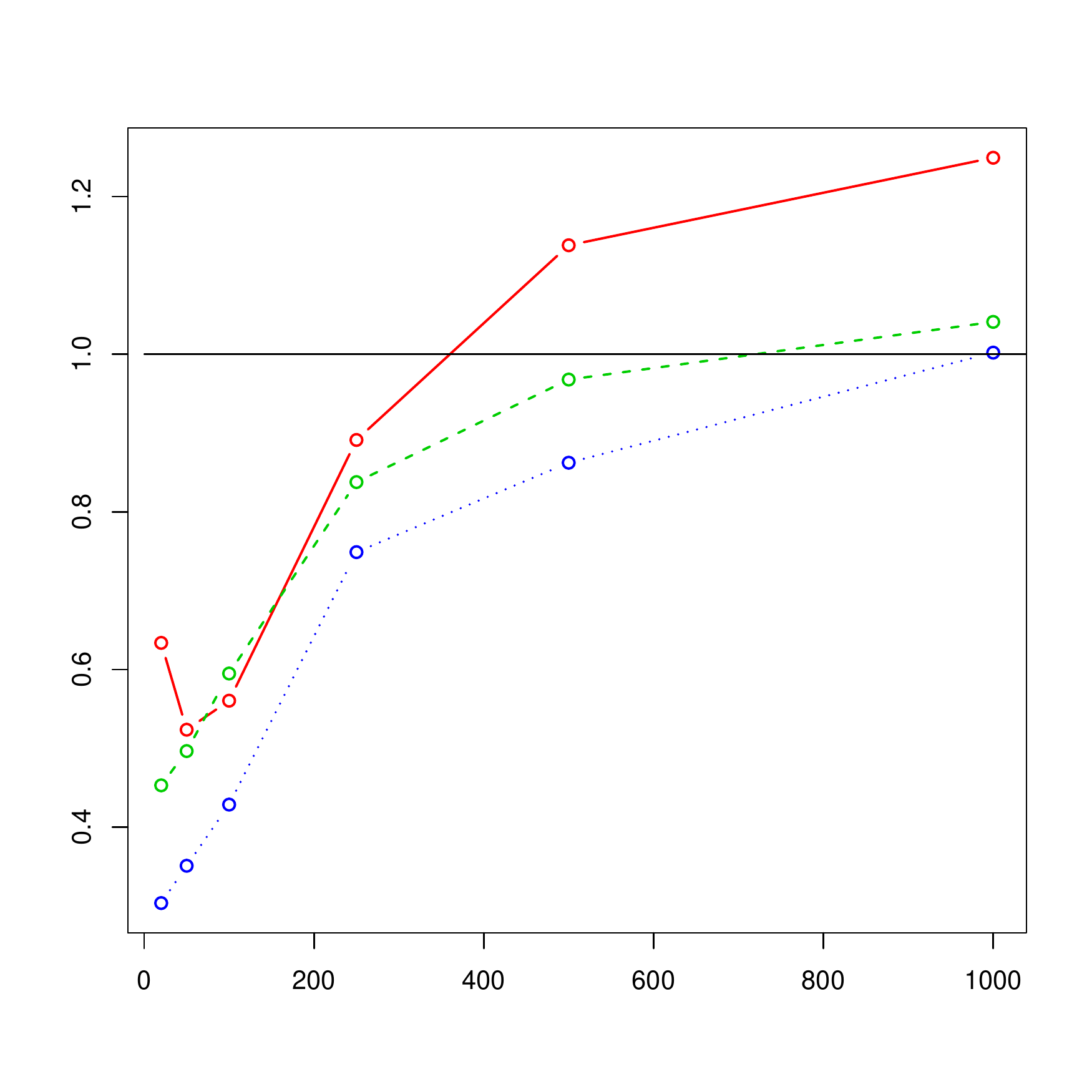}
\caption{$p=Q^{10}_{10}$.}
\label{fig15d}
\end{subfigure}
\caption{Spline distributions: ratio between the $\mathrm{RMSEP}$ of
  $\text{var}(\widehat{p}^{k}_{n})$ and the $\mathrm{RMSEP}$ of
  $\text{var}(\widetilde{p}_{n})$ versus the sample size $n$: for $k=2$ in "\textbf{\textcolor{red}{---}}",
  $k=3$ in "\textbf{\textcolor{green}{- -}}", $k=4$ in
  "\textbf{\textcolor{blue}{...}}". Each subfigure
  corresponds to the results  obtained with  $p = Q_{j}^{\ell}$, for
  $\ell\in\{2, 3, 4, 10\}$.}
\label{fig15}
\end{figure}
\FloatBarrier

\subsection{About the mass of the non-constrained estimator $\hat{p}^{*k}$}

\begin{figure}[!h]
\centering
\begin{subfigure}[b]{0.30\textwidth}
\includegraphics[width=0.95\textwidth]{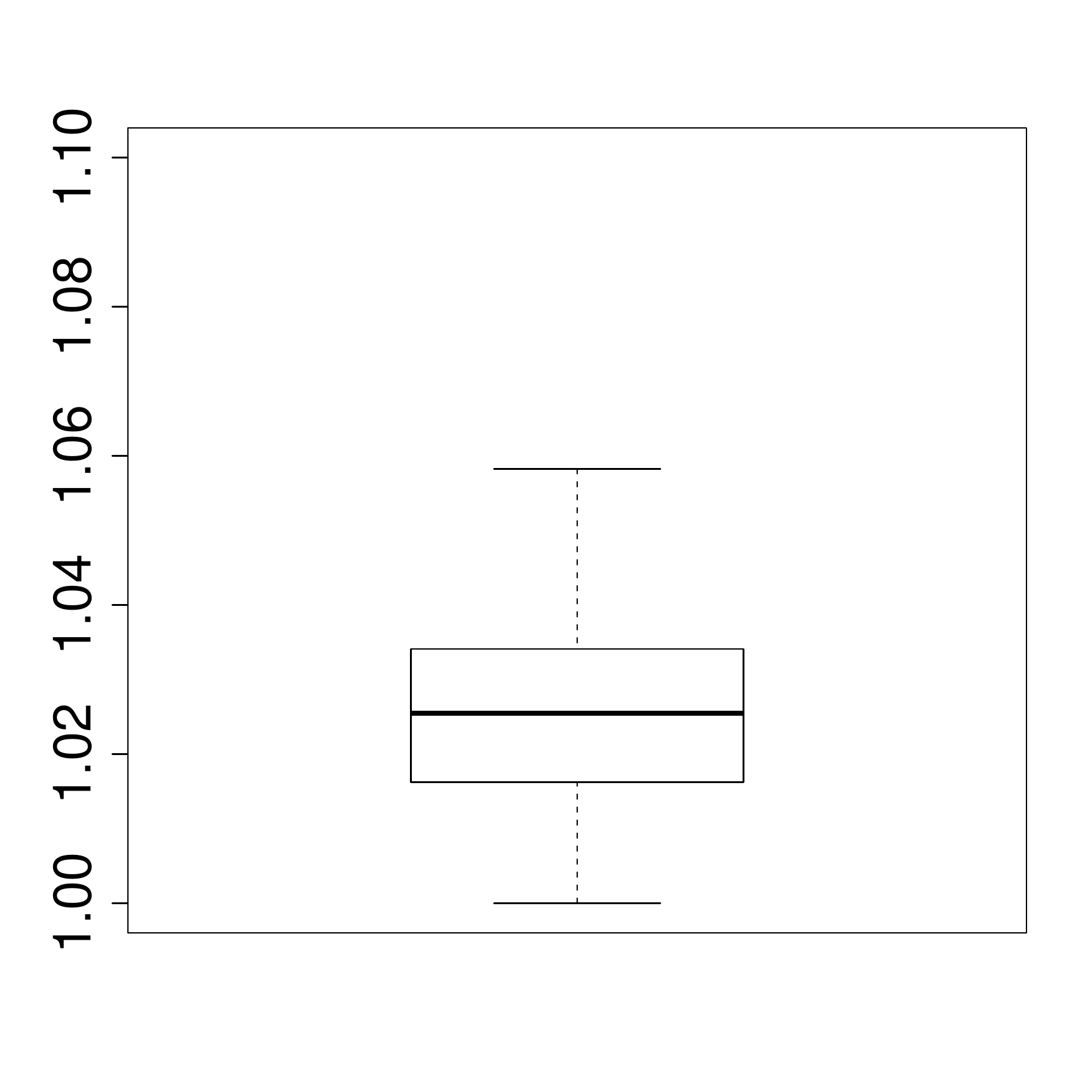}
\caption{$p=Q^2_{10}$.}
\label{fig9a}
\end{subfigure}
\hfill
\begin{subfigure}[b]{0.30\textwidth}
\includegraphics[width=0.95\textwidth]{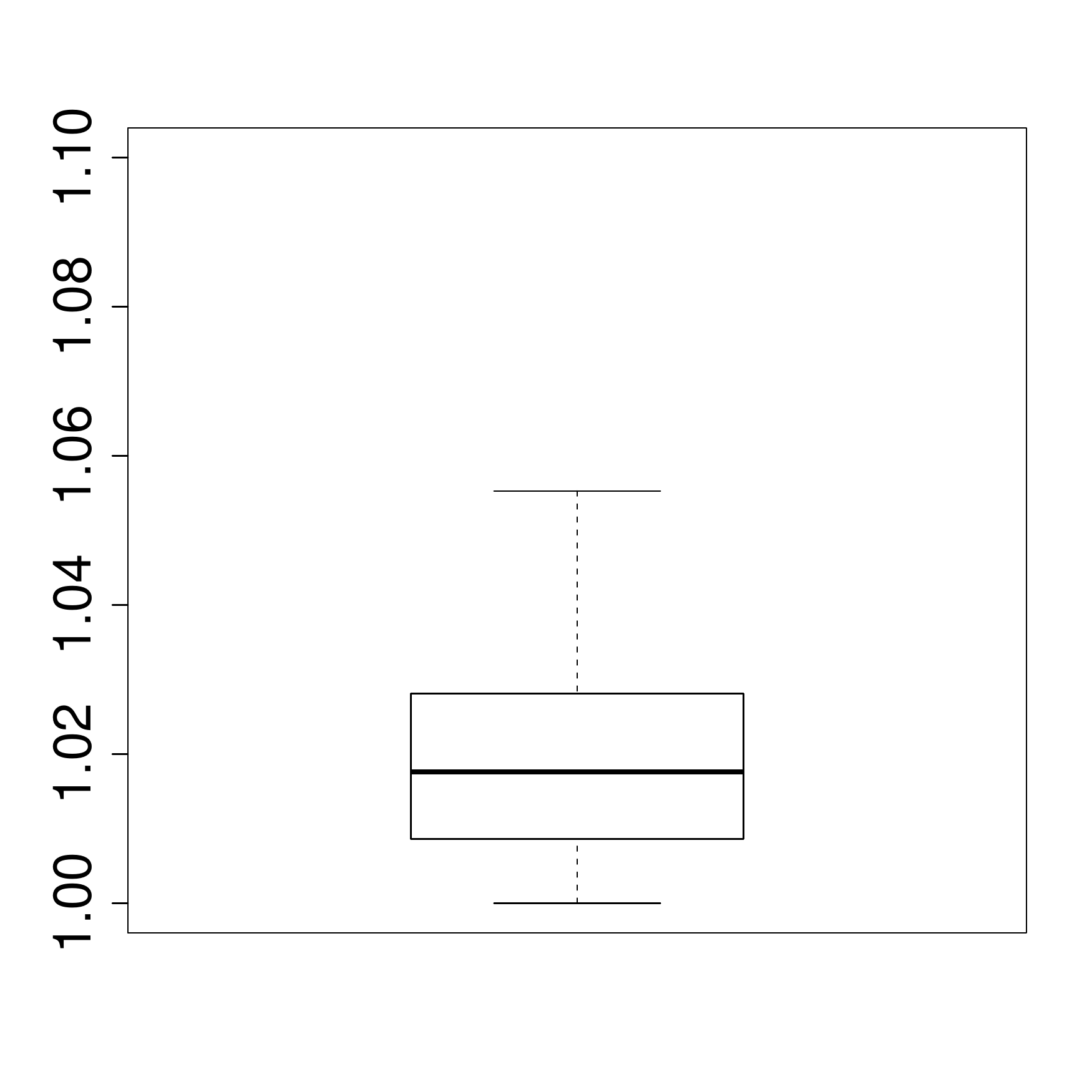}
\caption{$p=Q^4_{10}$.}
\label{fig9b}
\end{subfigure}
\hfill
\begin{subfigure}[b]{0.30\textwidth}
\includegraphics[width=0.95\textwidth]{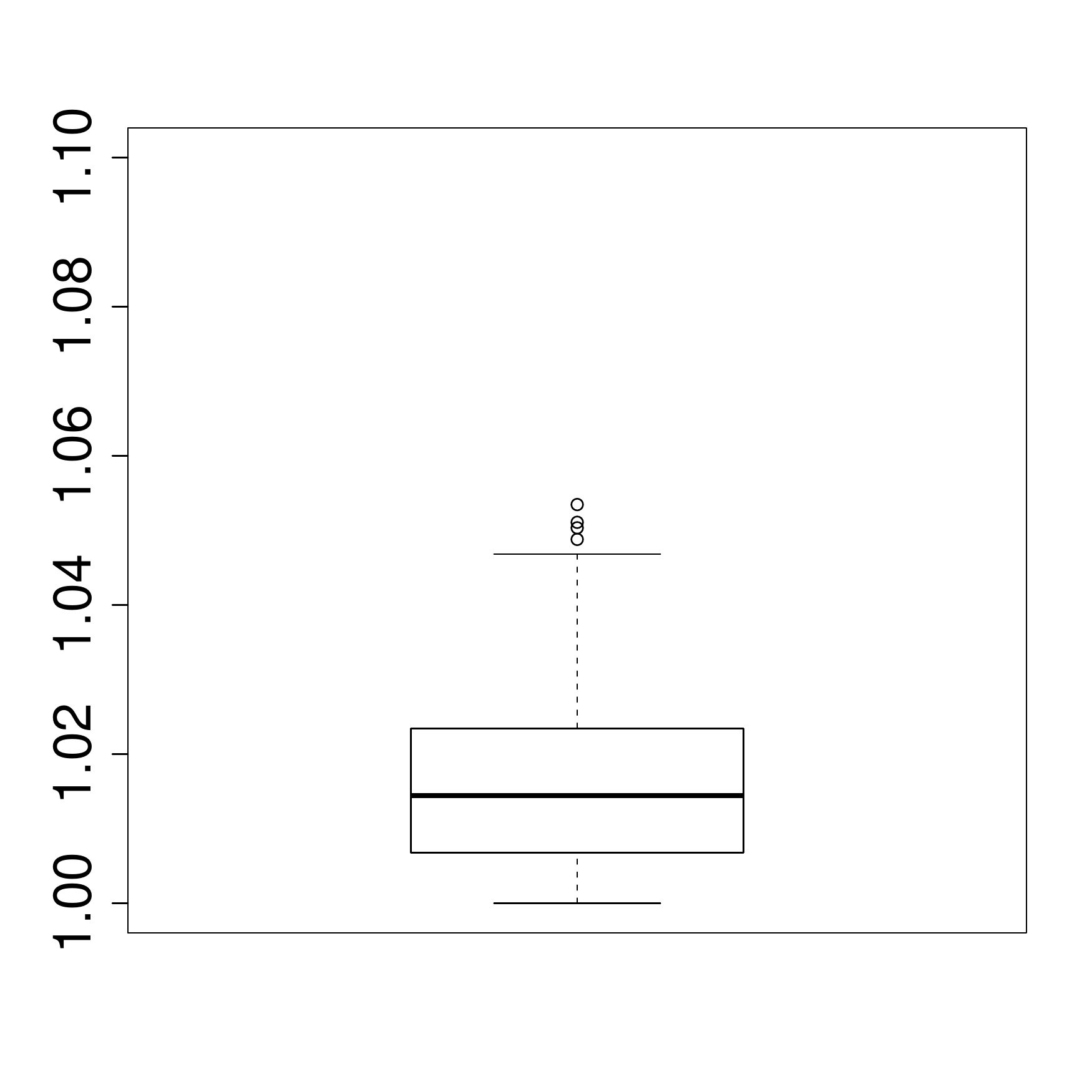}
\caption{$p=Q^{10}_{10}$.}
\label{fig9c}
\end{subfigure}
\begin{subfigure}[b]{0.30\textwidth}
\includegraphics[width=0.95\textwidth]{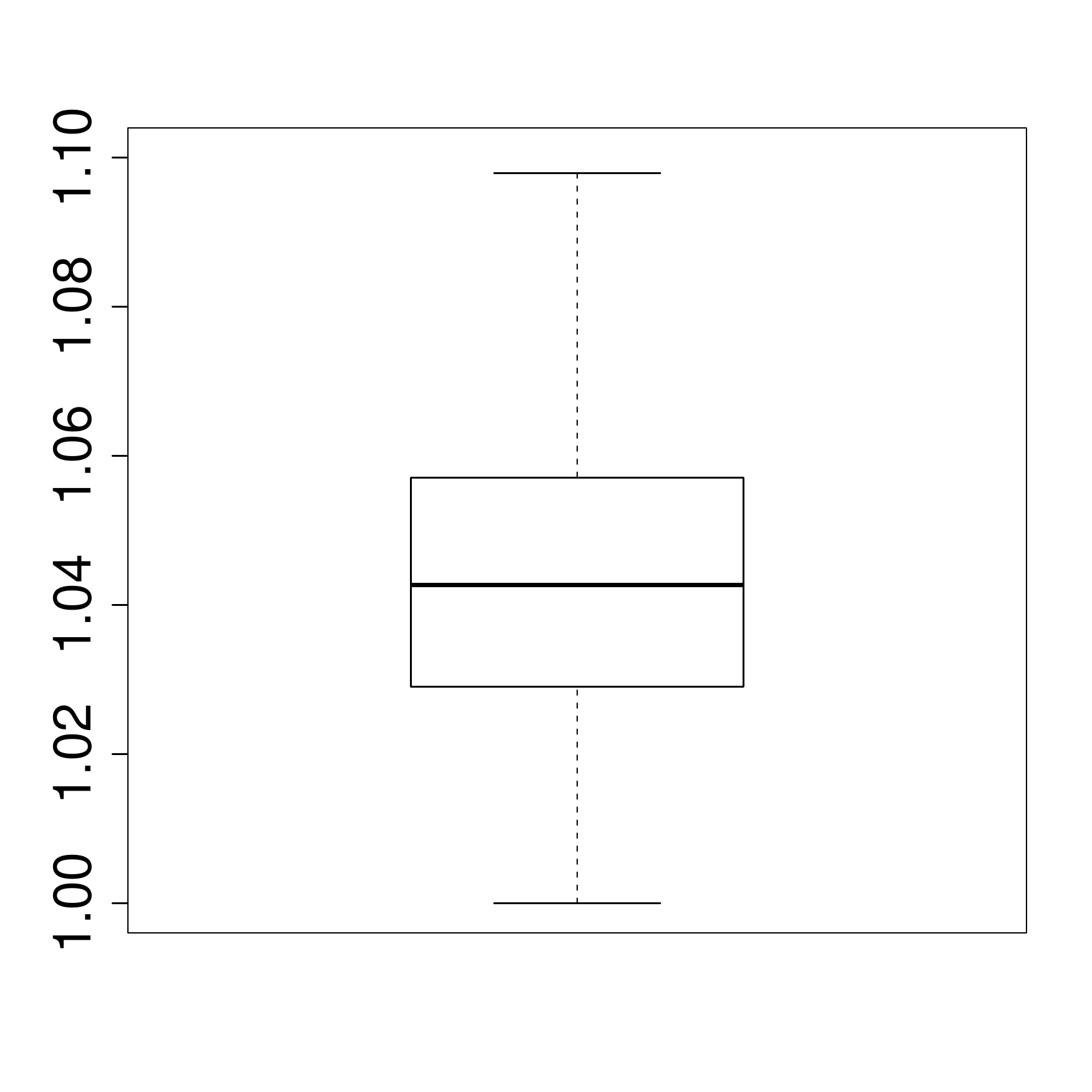}
\caption{$p=Q^2_{10}$.}
\label{fig9d}
\end{subfigure}
\hfill
\begin{subfigure}[b]{0.30\textwidth}
\includegraphics[width=0.95\textwidth]{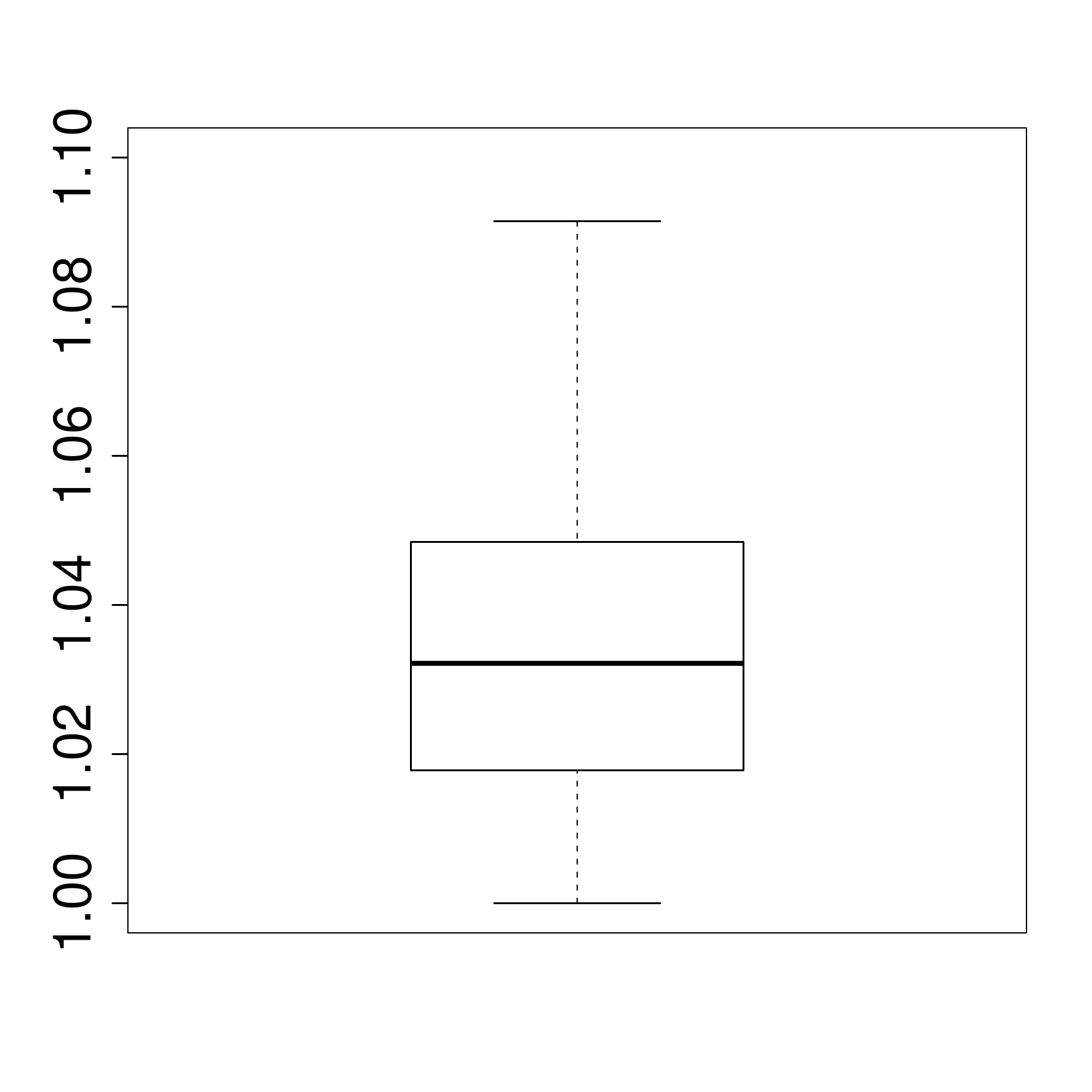}
\caption{$p=Q^4_{10}$.}
\label{fig9e}
\end{subfigure}
\hfill
\begin{subfigure}[b]{0.30\textwidth}
\includegraphics[width=0.95\textwidth]{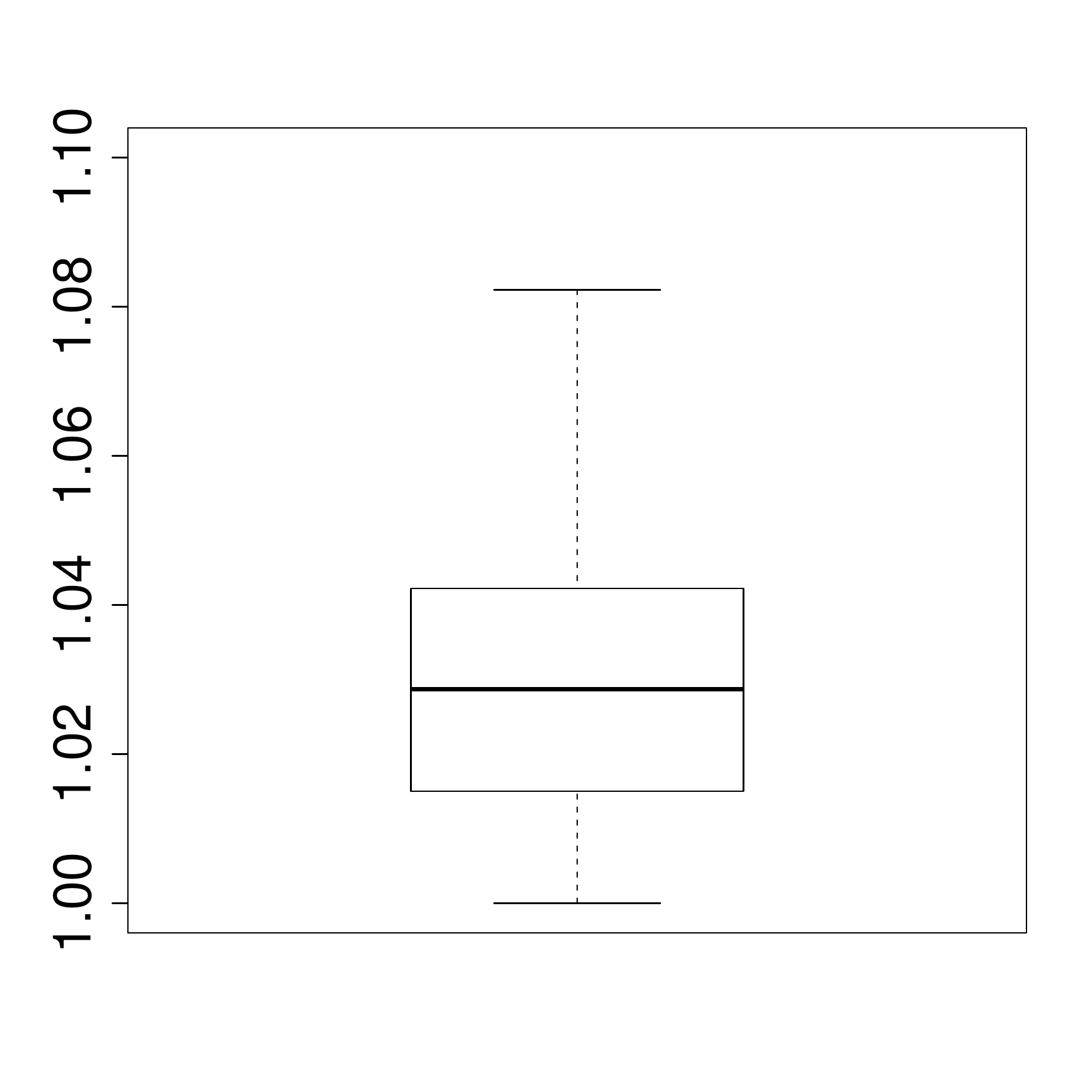}
\caption{$p=Q^{10}_{10}$.}
\label{fig9f}
\end{subfigure}
\caption{Splines distributions; $n=20$: Boxplot of the mass of $\hat{p}^{*k}$ for $k=3$ at the top and $k=2$ at bottom.  Each column corresponds to the results obtained with $p=Q_j^l$ for $l=2,4,10$.}
\label{fig9}
\end{figure}

\begin{figure}[!h]
\centering
\begin{subfigure}[b]{0.30\textwidth}
\includegraphics[width=0.95\textwidth]{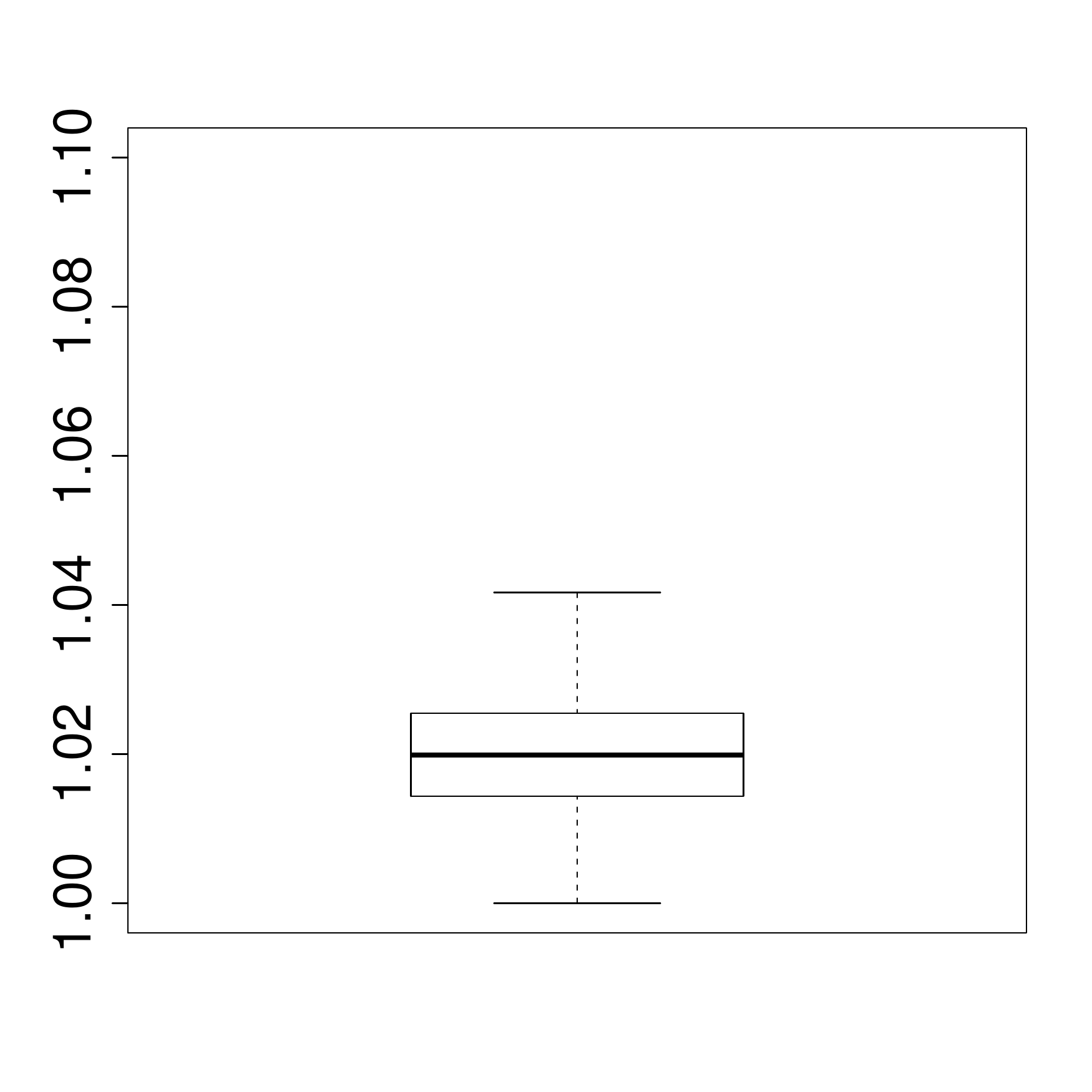}
\caption{$p=Q^2_{10}$.}
\label{fig10a}
\end{subfigure}
\hfill
\begin{subfigure}[b]{0.30\textwidth}
\includegraphics[width=0.95\textwidth]{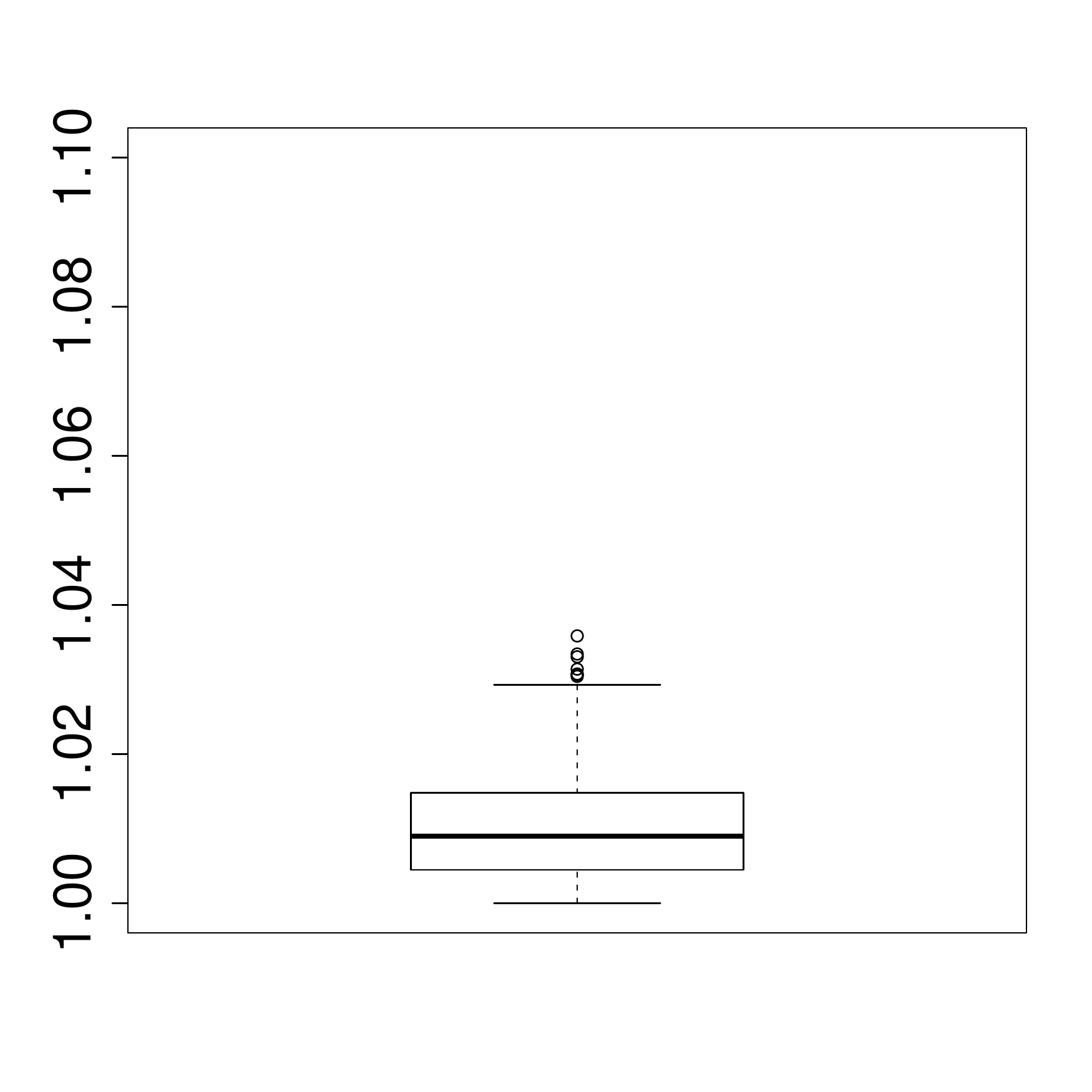}
\caption{$p=Q^4_{10}$.}
\label{fig10b}
\end{subfigure}
\hfill
\begin{subfigure}[b]{0.30\textwidth}
\includegraphics[width=0.95\textwidth]{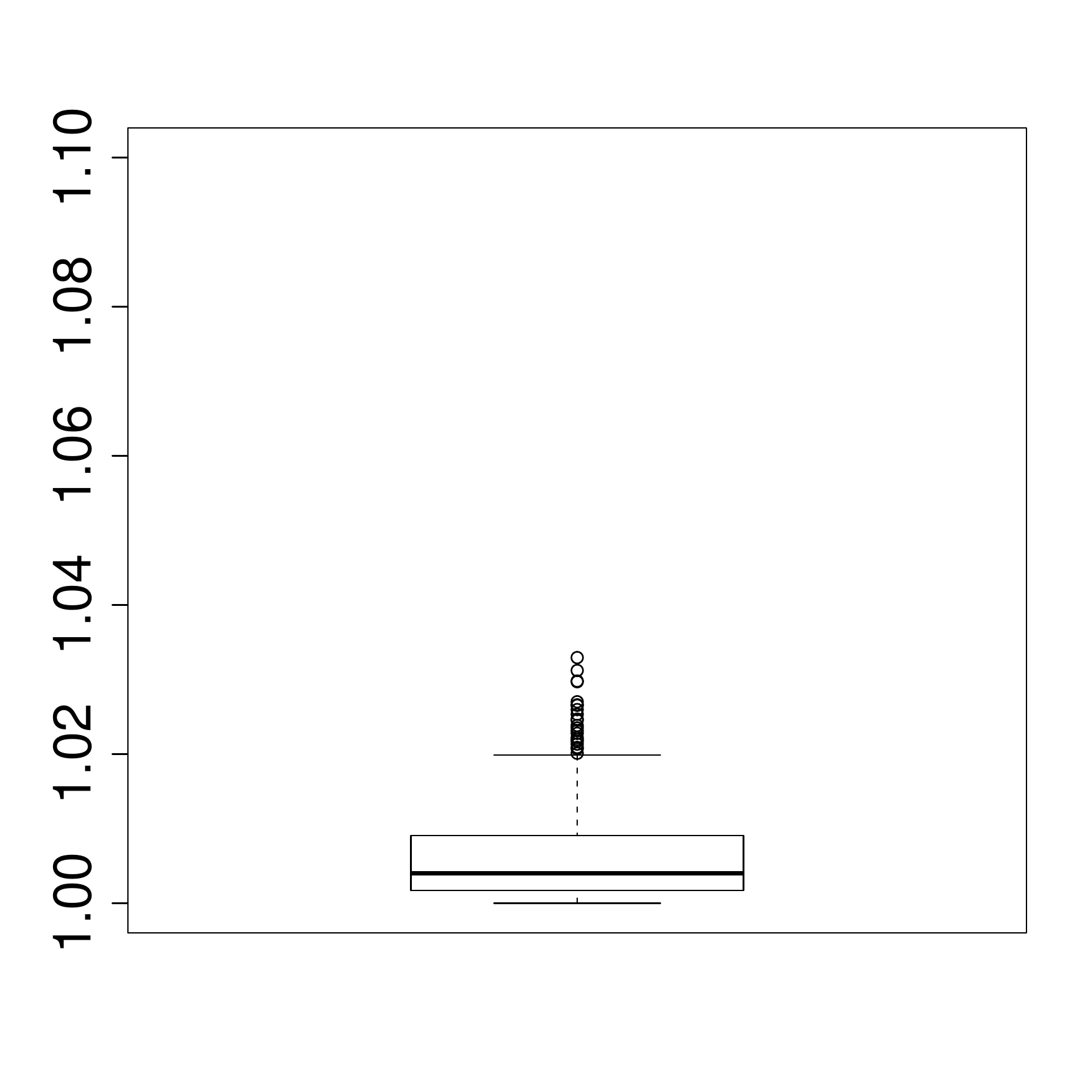}
\caption{$p=Q^{10}_{10}$.}
\label{fig10c}
\end{subfigure}
\begin{subfigure}[b]{0.30\textwidth}
\includegraphics[width=0.95\textwidth]{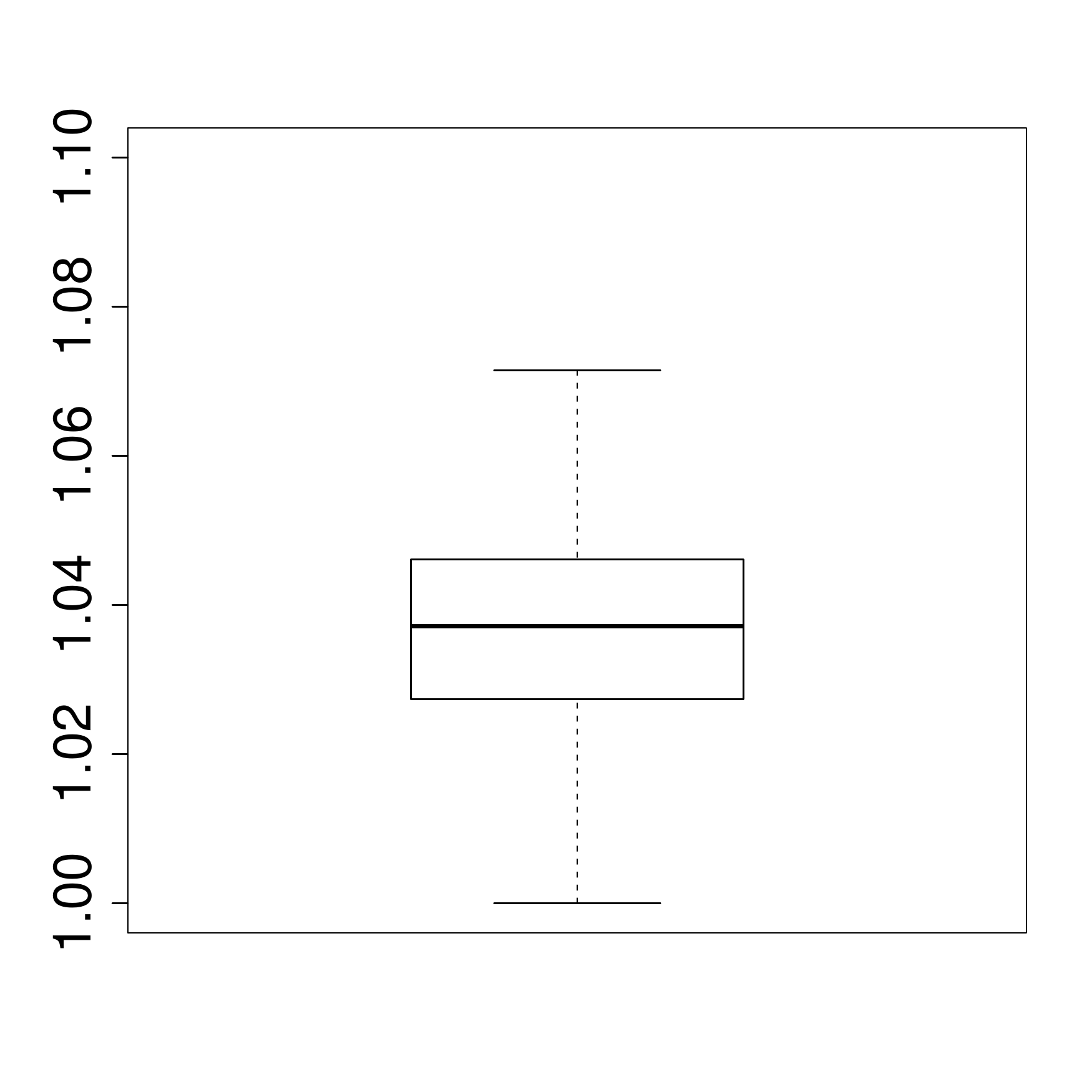}
\caption{$p=Q^2_{10}$.}
\label{fig10d}
\end{subfigure}
\hfill
\begin{subfigure}[b]{0.30\textwidth}
\includegraphics[width=0.95\textwidth]{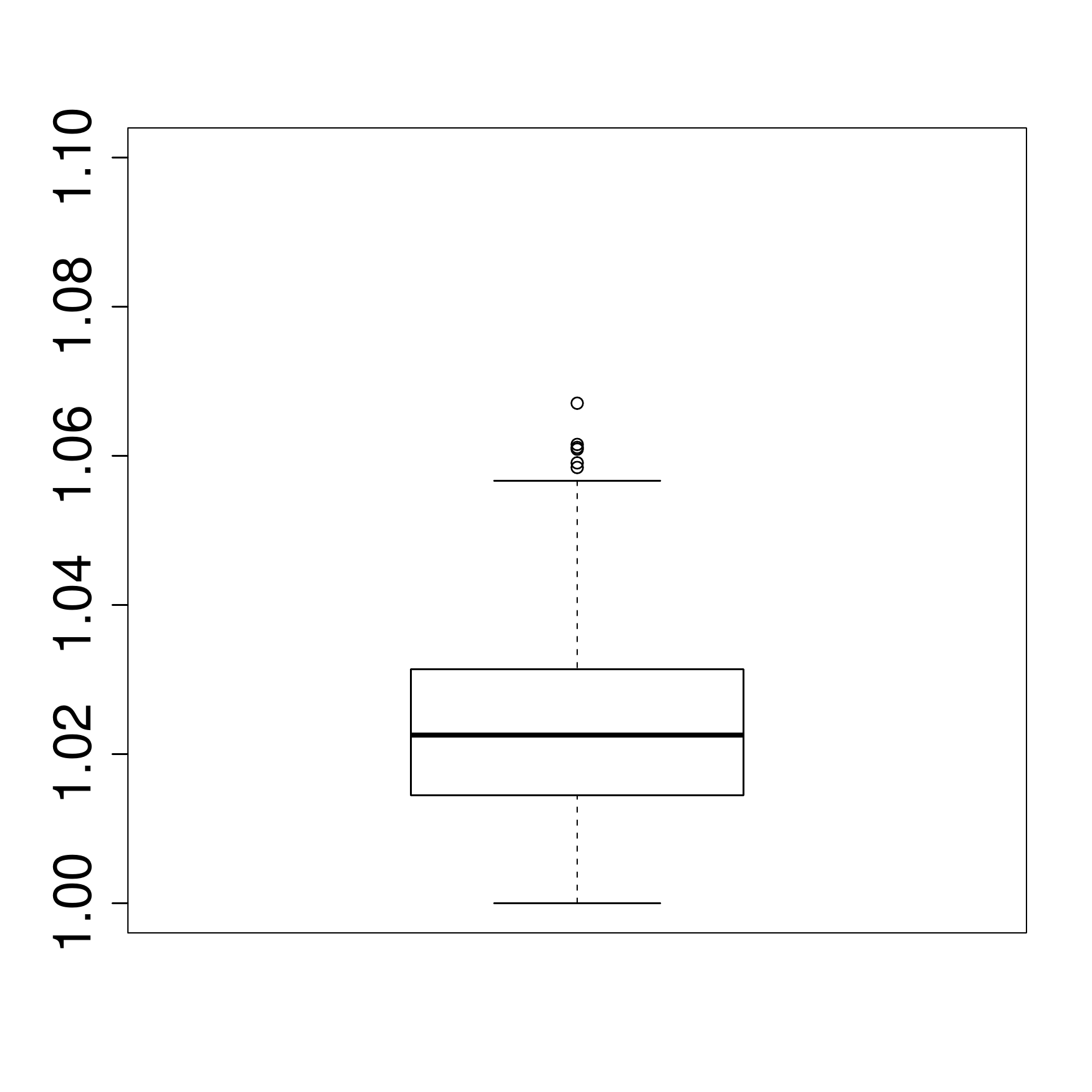}
\caption{$p=Q^4_{10}$.}
\label{fig10e}
\end{subfigure}
\hfill
\begin{subfigure}[b]{0.30\textwidth}
\includegraphics[width=0.95\textwidth]{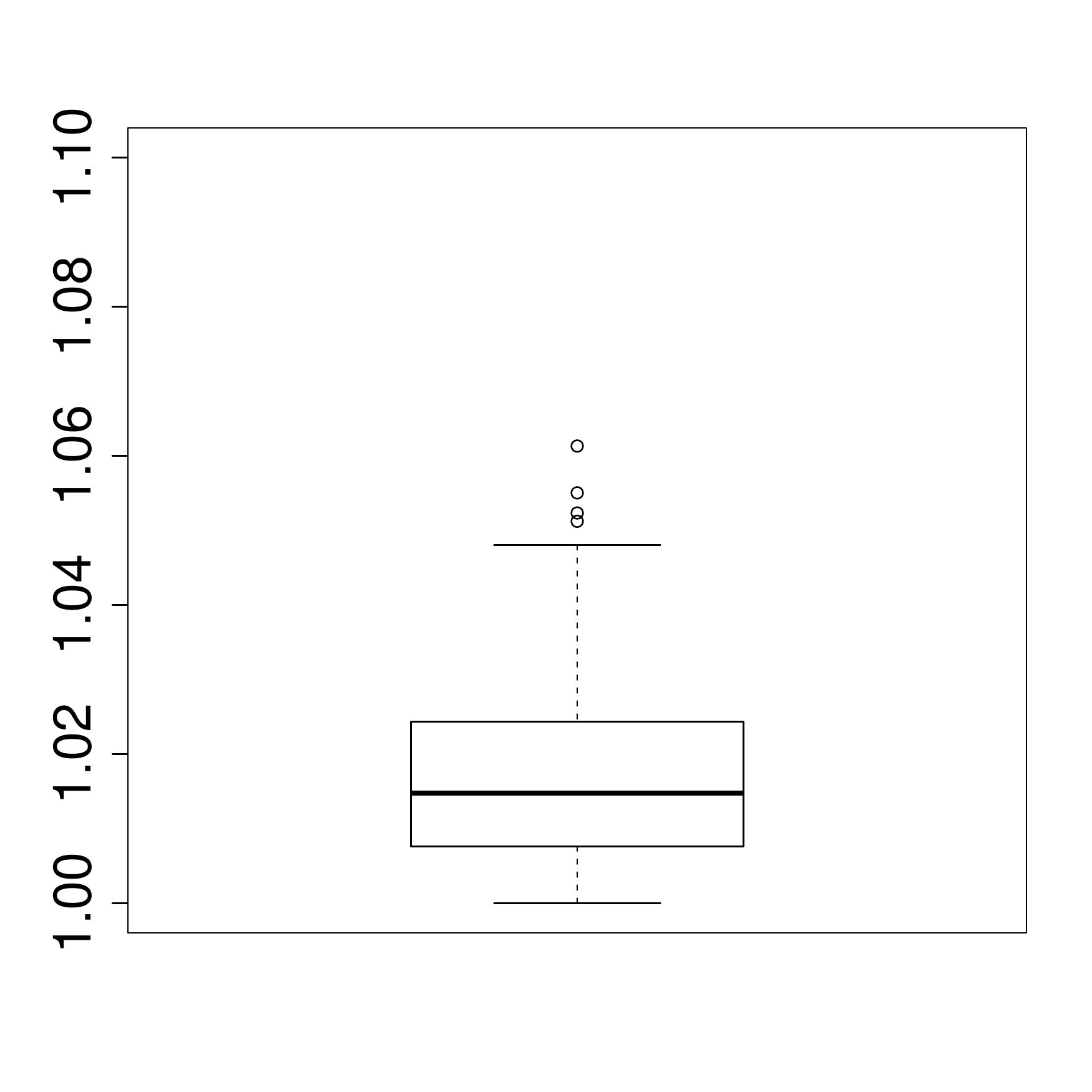}
\caption{$p=Q^{10}_{10}$.}
\label{fig10f}
\end{subfigure}
\caption{Repartition of the mass of $\hat{p}^{*k}$ for $n=100$.  Each column represents the estimation of a different probability $p$ explained in subtitle. The first line is for the mass of $\hat{p}^{*3}$ and the second line for the mass of $\hat{p}^{*4}$.}
\label{fig10}
\end{figure}

We were also interested in the estimation of the mass of the non-constrained estimator $\hat{p}^{*k}$.
Figures \ref{fig9} and \ref{fig10} illustrated the results for the spline distributions with $n=20$ and $n=100$. As expected the mass is always larger than 1 and whatever $k$, the distribution of the mass comes closer to one when $n$ increases (compare figures \ref{fig9} and \ref{fig10}). The larger $l$ is, the smaller the median and the dispersion around the median are.  On the other hand when $k$ increases the distributions are more scattered and their medians move away from $1$ (compare the lines of each figure).

\FloatBarrier
\FloatBarrier

\FloatBarrier
\FloatBarrier

\subsection{Conclusion}

Let us consider the case where $p^*$ is $l-$monotone and $\hat{p}_n^k$ is the least-squares estimator of $p^*$ on $\mathcal{S}_k$ for $k\leq l$.\\
Concerning the $l_2$-loss, the total variation loss and the estimation of $p^*(0)$, $\hat{p}^k_n$ performs better than the empirical estimator $\tilde{p}_n$. Moreover the superiority of the performance of $\hat{p}^k_n$ is larger when $n$ is small.\\
Concerning the Hellinger loss, or the estimation of the variance and the entropy, we get the following results. For small $n$, as before, the least-squares estimator is always better than the empirical estimator $\tilde{p}_n$. When $n$ is large, $\hat{p}^k_n$ and $\tilde{p}_n$ behave similarly. If $p$ is a \textit{frontier} distribution in $\mathcal{S}_l$, as for example the Poisson distribution with $\lambda=\lambda_l$ or a spline distribution $Q_j^l$, then $\hat{p}^{l-1}_n$ performs better than $\hat{p}^l_n$. If not, then $\hat{p}^l_n$ performs better than $\hat{p}^k_n$ for all $k\leq l$.\\
Finally, for all considered criteria, the estimator $\hat{p}^k_n$ performs better than $\hat{p}^{*k}_n$ when $n$ is small and both estimators perform similarly when $n$ is large.

\section{Proofs\label{proof.section}}

For all discrete function $f$, let $\mathcal{Q}(f)=\frac{1}{2}\vert\vert f-\tilde{p}\vert\vert^2$ . When no confusion is possible we write $\hat{p}$ instead of $\hat{p}^k$.

\subsection{Properties of the estimator} 
\subsubsection{Proof of Theorem \ref{theo2}\label{theo2.st} : Characterization of $\hat{p}$.}

Let us first prove that $\hat{p}$ satisfies \textit{1.} 
or equivalently, that for all integer $l$  
the following inequality is satisfied:
\begin{eqnarray}
\label{aux82}
F^{k}_{\hat{p}}(l)- F^{k}_{\tilde{p}}(l)\geqslant \beta(f)m^k_l.
\end{eqnarray}
By definition $\beta(\hat{p})=<\hat{p},\hat{p}-\tilde{p}>$, then (\ref{aux82}) is equivalent to:
\begin{eqnarray}
\label{aux83}
\frac{1}{m^k_l}(F^{k}_{\hat{p}}(l)-F^{k}_{\tilde{p}}(l))-
\sum_{i=0}^{\infty}\hat{p}(i)(\hat{p}(i)-\tilde{p}(i)) \geqslant 0.
\end{eqnarray}

Let us rewrite this equation
by considering limits of the directionnal derivatives.\\
For all $\epsilon\in]0,1]$, $l\geqslant 0$ we define a function $q_{\epsilon l}$ as follows:
\begin{equation}
\label{defqeps}
 q_{\epsilon l}(i)=(1-\epsilon)\hat{p}(i)+\epsilon \frac{\bar{Q}^{k}_{l}(i)}{m^k_l}=\left \{
    \begin{array}{l}
         (1-\epsilon)\hat{p}(i)+\frac{\epsilon}{m^k_l} \bar{Q}^{k}_{l}(i) \text{ if $i\in\{0,\ldots,l\}$}\\
         (1-\epsilon)\hat{p}(i) \text{ if $i\geqslant l+1$}.
    \end{array}
 \right.
\end{equation}

The function $q_{\epsilon l}$ is a $k-$monotone probability and because $\hat{p}$ minimizes $\mathcal{Q}$ on the set of $k$-monotone probabilities, we have $\mathcal{Q}(q_{\epsilon l})\geqslant \mathcal{Q}(\hat{p})$ and:
\begin{eqnarray*}
\liminf_{\epsilon\searrow 0}\frac{1}{\epsilon}(\mathcal{Q}(q_{\epsilon l})-\mathcal{Q}(\hat{p}))\geqslant 0,
\end{eqnarray*}
that is equivalent to:
\begin{eqnarray*}
\liminf_{\epsilon\searrow 0}\frac{1}{\epsilon}\left(\sum_{i=0}^{\infty}\big((1-\epsilon)\hat{p}(i)
+\frac{\epsilon}{m^k_l} \bar{Q}^{k}_{l}(i)-\tilde{p}(i)\big)^2-\sum_{i=0}^{\infty}(\hat{p}(i)-\tilde{p}(i))^2\right)\geqslant 0.
\end{eqnarray*}
Therefore we have the following inequalities:
\begin{eqnarray*}
\liminf_{\epsilon\searrow 0}\frac{1}{\epsilon}\left(\sum_{i=0}^{\infty}\big[(\hat{p}(i)-\tilde{p}(i))^2
+2\epsilon\big(\hat{p}(i)-\tilde{p}(i)\big)\big(\frac{\bar{Q}^{k}_{l}(i)}{m^k_l}-\hat{p}(i)\big)+\epsilon^2\big(\frac{\bar{Q}^{k}_{l}(i)}{m^k_l}-\hat{p}(i)\big)^2\big]\right.\\\left.-\sum_{i=0}^{\infty}(\hat{p}(i)-\tilde{p}(i))^2\right)\geqslant 0,
\end{eqnarray*}

\begin{eqnarray*}
\liminf_{\epsilon\searrow 0}\frac{1}{\epsilon}\left(\epsilon^2\sum_{i=0}^{\infty}(\frac{1}{m^k_l}\bar{Q}^{k}_{l}(i)-\hat{p}(i))^2+2\epsilon\sum_{i=0}^{\infty}(\hat{p}(i)-\tilde{p}(i))(\frac{\bar{Q}^{k}_{l}(i)}{m^k_l}-\hat{p}(i))\right)\geqslant 0,
\end{eqnarray*}

\begin{eqnarray*}
\sum_{i=0}^{\infty}(\hat{p}(i)-\tilde{p}(i))\frac{\bar{Q}^{k}_{l}(i)}{m^k_l}- \sum_{i=0}^{\infty}\hat{p}(i)(\hat{p}(i)-\tilde{p}(i))\geqslant 0.
\end{eqnarray*}

Now, using Lemma \ref{rec.lemme} (see
  Section~\ref{prooflemmas.st}) we have that for all $k\geqslant 2$
and for all positive discrete measure $f$:
\begin{eqnarray*}
 \forall l\in\mathbb{N}^*, \sum_{i=0}^{\infty}f(i) \bar{Q}^{k}_{l}(i)=\sum_{i=0}^{l}f(i) \bar{Q}^{k}_{l}(i)=F^k_f(l).
\end{eqnarray*}
We choose $f=\hat{p}$ and we obtain exactly (\ref{aux83}).

Let us now show that $\hat{p}$ satisfies \textit{2.}. Let $l$ be a $k-$knot
of  $\hat{p}$, we need to show that Inequality (\ref{aux83}) is an
equality. As before we consider $q_{\epsilon l}$ defined at Equation (\ref{defqeps}) and show that $q_{\epsilon l}$ is a $k-$monotone probability for $\epsilon$ nonpositif small enough. Thanks to the following equality:
\begin{eqnarray*}
(-1)^k\Delta^k \bar{Q}^k_l(i)=\left\{
    \begin{array}{l}
    \vspace{0.25cm}
     0 \text{ if $i\neq l$}\\ 
     1 \text{ if $i=l$}
\end{array}
\right.
\end{eqnarray*}
we get :
\begin{eqnarray*}
(-1)^k\Delta^k q_{\epsilon l}(i)=\left\{\begin{array}{l}
\vspace{0.25cm}
     (1-\epsilon)(-1)^k\Delta^k \hat{p}(l)+\epsilon/m_l \text{ if $i=l$}\\
        
     (1-\epsilon)(-1)^k\Delta^k \hat{p}(i) \text{ if $i\neq l$}.
     \end{array}
\right.
\end{eqnarray*}
Because $\hat{p}$ is $k-$monotone, $(-1)^k\Delta^k q_{\epsilon l}(i)\geqslant 0$ for $\epsilon$ nonpositive small enough and $i\neq l$. As $l$ is a $k-$knot of $\hat{p}$, $(-1)^k\Delta^k \hat{p}(l)>0$   then $(-1)^k\Delta^k q_{\epsilon l}(i)\geqslant 0$. Therefore we get:
\begin{eqnarray*}
\left\{
    \begin{array}{l}
    \vspace{0.25cm}
     \liminf_{\epsilon\searrow 0}\frac{1}{\epsilon}\big(\mathcal{Q}(q_{\epsilon l})-\mathcal{Q}(\hat{p})\big)\geqslant 0\\ 
    \limsup_{\epsilon\nearrow 0}\frac{1}{\epsilon}\big(\mathcal{Q}(q_{\epsilon l})-\mathcal{Q}(\hat{p})\big)\leqslant 0,
\end{array}
\right.
\end{eqnarray*}
which leads us to:
\begin{eqnarray*}
\sum_{i=0}^{\infty}(\hat{p}(i)-\tilde{p}(i))\frac{\bar{Q}^{k}_{l}(i)}{m^k_l}- \sum_{i=0}^{\infty}\hat{p}(i)(\hat{p}(i)-\tilde{p}(i))= 0,
\end{eqnarray*}
that is exactly an equality in (\ref{aux83}).

Conversely assuming that $f$ is a $k-$monotone probability that satisfies :
\begin{eqnarray}
\label{aide4}
\frac{F^{k}_{f}(l)- F^{k}_{\tilde{p}}(l)}{m^k_l}\geqslant \beta(f),
\end{eqnarray}
with equality if $l$ is a $k-$knot of $f$, we have to show that $f=\hat{p}$. By definition of $\hat{p}$ we need to show that for all $k-$monotone probability $g$ we have $\mathcal{Q}(g)\geqslant \mathcal{Q}(f)$.\\
Let $g$ be a $k$-monotone probability. Using Lemma \ref{lemme-beta} (see Section \ref{prooflemmas.st}):
\begin{align*}
\mathcal{Q}(g)-\mathcal{Q}(f)&=\frac{1}{2}\vert\vert g-f\vert\vert^2_2+<f-\tilde{p},g-f>\\
&\geqslant\text{ } <f-\tilde{p},g-f>\\
&= \sum_{i=0}^{\infty}(g(i)-f(i))(f(i)-\tilde{p}(i))\\
&=\sum_{i=0}^{\infty}(-1)^k\Delta^{k}(g-f)_i(F_f^k(i)-F_{\tilde{p}}^k(i)).
\end{align*}
The function $g$ is $k-$monotone then for all $i$, $(-1)^k\Delta^{k}g(i)\geqslant 0$ and using \eref{aide4} and lemma \ref{lemme-beta} (see Section \ref{prooflemmas.st}):
\begin{align*}
\mathcal{Q}(g)-\mathcal{Q}(f)&\geqslant \sum_{i=0}^{\infty}(-1)^k\Delta^{k}g(i)(F_f^k(i)-F_{\tilde{p}}^k(i))-\sum_{i=0}^{\infty}(-1)^k\Delta^{k}f(i)(F_f^k(i)-F_{\tilde{p}}^k(i))\\
&\geqslant \sum_{i=0}^{\infty}(-1)^k\Delta^{k}g(i)\beta(f)m^k_i-\sum_{i=0}^{\infty}(-1)^k\Delta^{k}f(i)(F_f^k(i)-F_{\tilde{p}}^k(i))\\
&\geqslant \beta(f)\sum_{i=0}^{\infty}(-1)^k\Delta^{k}g(i)m^k_i-\beta(f).\\
\end{align*}
Moreover $g$ being a $k-$monotone probability, according to Property \ref{decomposition}, we have the decomposition on the spline basis:
\begin{eqnarray*}
\sum_{i=0}^{\infty}(-1)^k\Delta^{k}g(i)m^k_i=1.
\end{eqnarray*}
Finally for all $k-$monotone probability $g$ we find :
\begin{eqnarray*}
\mathcal{Q}(g)-\mathcal{Q}(f)\geqslant \beta(f)-\beta(f)=0.
\end{eqnarray*}
By unicity of the projection we have $f=\hat{p}$.
\\

\subsubsection{Proof of Theorem \ref{supportfini.prop} : Support of $\hat{p}$}

\paragraph{The support of $\hat{p}$ is finite.}
Let us first consider the case where $\beta(\hat{p})=0$. According to Property \ref{beta} this is equivalent to $\hat{p}=\hat{p}^{*}$.\\

The result is proved by contradiction. Let us assume that $\hat{p}$ has an infinite support then we can build a probability $\bar{p}$ satisfying the following properties:
\begin{align*}
\text{i) } &\bar{p}\leqslant \hat{p}.\\
\text{ii) } &\text{ for all } i\leqslant \tilde{s}, \bar{p}(i)=\hat{p}(i).\\
\text{iii) } &\text{ there exists $i$ such as $\bar{p}(i)<\hat{p}(i)$}.\\
\text{iv) } &\bar{p} \text{ is $k$-monotone and non-negative},
\end{align*}
with $\tilde{s}$ the maximum of the support of $\tilde{p}$.\\

For this $\bar{p}$ we have the inequality $\vert\vert\bar{p}-\tilde{p}\vert\vert_{2}<\vert\vert\hat{p}-\tilde{p}\vert\vert_{2}$
which contradicts the definition of $\hat{p}$.\\

The probability $\bar{p}$ is constructed as follows.\\
We define for all $j\in\{1,\ldots,k-1\}$ and for all $i\in\mathbb{N}$,  the $j^{\text{th}}$ derivative function $q_j$ of $\hat{p}$ : 
\begin{equation}
\label{defder}
\left.\begin{array}{rcl}
q_{1}(i) &=& -\hat{p}(i+1)+\hat{p}(i)=-\Delta^1\hat{p}(i),\\
q_{2}(i) &=& -q_{1}(i+1)+q_{1}(i)=\Delta^2\hat{p}(i),\\
         &\vdots&\\
q_{k-1}(i) &=& -q_{k-2}(i+1)+q_{k-2}(i)=(-1)^{k-1}\Delta^{k-1}\hat{p}(i)
\end{array}
\right.         
\end{equation}

We have $q_{j+1}(i)=\Delta^1 q_{j}(i)$ so for all $i\in\mathbb{N}$:
\begin{eqnarray*} (-1)^{k}\Delta^{k}\hat{p}(i)=(-1)^{k-1}\Delta^{k-1}(q_{1}(i))=\ldots=\Delta^1 q_{k-1}(i).
\end{eqnarray*}
Then $\hat{p}$ is $k$-monotone (and non-negative) if and only if $q_{k-1}$ is non-increasing (and non-negative).\\
Because $\hat{p}$ has an infinite support, all the functions $q_j$ have infinite suppport too. \\
Moreover for all $i\in\mathbb{N}$ we have the following inequalities:
\begin{eqnarray*}
\hat{p}(i)=-\sum_{h=0}^{i-1}q_{1}(h)+\hat{p}(0),\\
\forall j\in\{1,\ldots,k-2\}, q_{j}(i)=-\sum_{h=0}^{i-1}q_{j+1}(h)+q_{j}(0).
\end{eqnarray*}
The next step is to modify $q_{k-1}$ to $\bar{q}_{k-1}$ such that if $\bar{q}_j$ is defined as:
\begin{eqnarray*}
\bar{q}_{j}(i)=-\sum_{h=0}^{i-1}\bar{q}_{j+1}(h)+q_{j}(0), \forall j\in\{1,\hdots,k-2\},\\
\end{eqnarray*}
and if $\bar{p}$ is defined as:
\begin{eqnarray}
\label{aux73}
\bar{p}(i)=-\sum_{h=0}^{i-1}\bar{q}_{1}(h)+\hat{p}(0), \forall j\in\{1,\ldots,k-2\},
\end{eqnarray}
then $\bar{p}$ satisfies i)ii)iii)iv).\\

\begin{figure}[!h]
\centering
\includegraphics[width=0.45\textwidth]{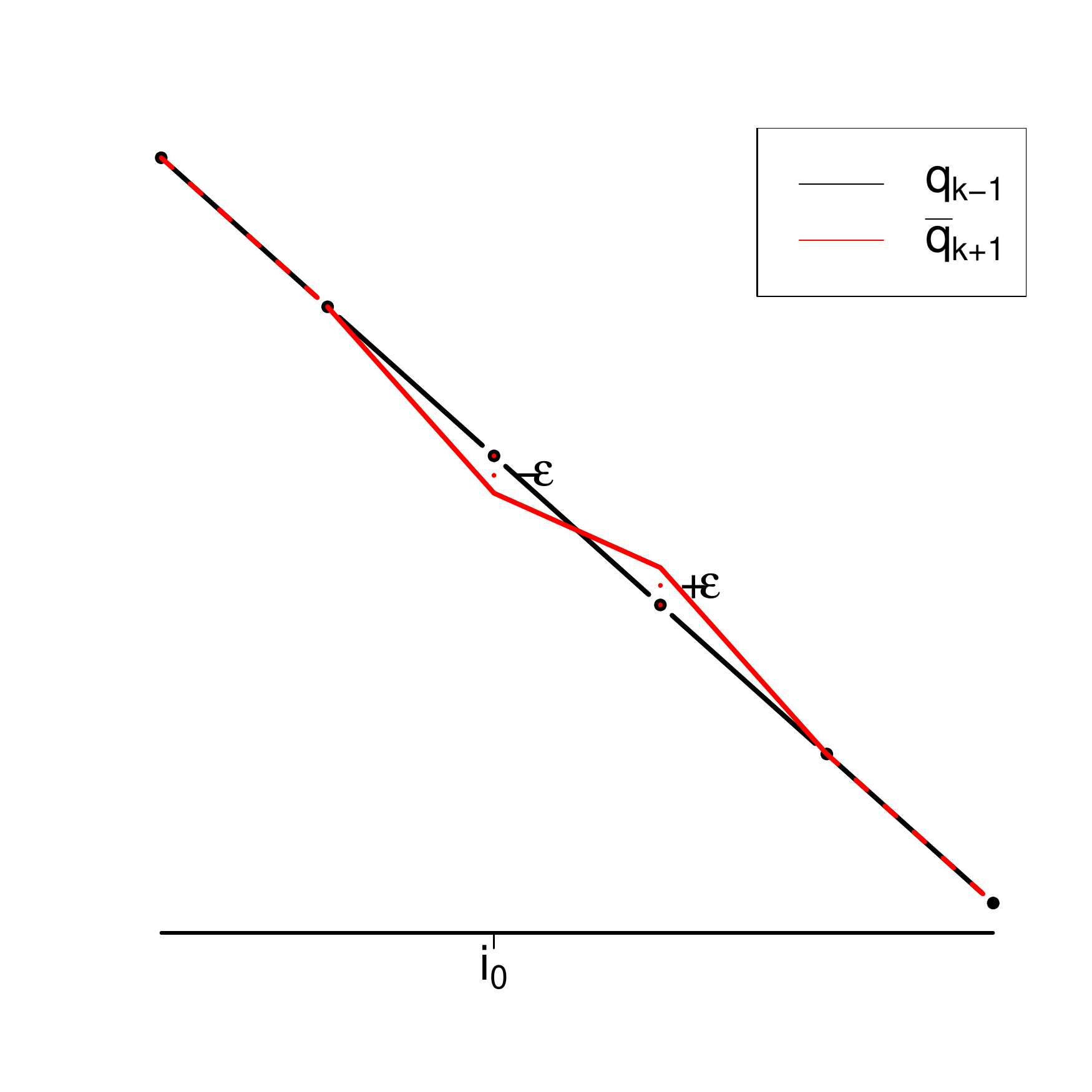}
\caption{Functions $q_{k-1}(i)$ and $\bar{q}_{k-1}(i)$ versus i ($k$ odd).}
\label{fig1}
\end{figure}

The function $q_{k-1}$ has an infinite support and is non-increasing, therefore it has an infinity of $1-$knots (points where $q_{k-1}$ is strictly non-increasing). \\
Assume first that $k$ is odd ($k\geqslant 3$). Let $i_0$ be a $1-$knot of $q_{k-1}$ such that $i_0>\tilde{s}$. We define $\bar{q}_{k-1}$ as follows: 
\begin{eqnarray}
\label{aux74}
\left\{
    \begin{array}{l}
    \vspace{0.25cm}
    \bar{q}_{k-1}(i)=q_{k-1}(i) \text{   if $i\neq i_0, i_0+1$}\\
    \vspace{0.25cm}
    \bar{q}_{k-1}(i_0)=q_{k-1}(i_0)-\epsilon\\
    \vspace{0.25cm}
    \bar{q}_{k-1}(i_0+1)=q_{k-1}(i_0+1)+\epsilon.
    \end{array}
  \right.
\end{eqnarray}
where $\epsilon$ is some positive real number chosen such that $\bar{q}_{k-1}$ is still non-increasing. For example take $\epsilon=(\bar{q}_{k-1}(i_0)-\bar{q}_{k-1}(i_0+1))/2$. The function $\bar{q}_{k-1}$ is shown at Figure \ref{fig1}.

Then the distribution $\bar{p}$ defined at Equation (\ref{aux73}) satisfies $iv)$.
\\
To show the properties $i)$ to $iii)$, we will use the following equality whose proof is straightforward and omitted:
\begin{eqnarray}
\label{lemme7.eq}
\forall i\in\mathbb{N}, \hat{p}(i)-\bar{p}(i)=(-1)^{k-1}\sum_{h_1=0}^{i-1}\sum_{h_2=0}^{h_1-1}\ldots\sum_{h_{k-1}=0}^{h_{k-2}-1}\big(q_{k-1}(h_{k-1})-\bar{q}_{k-1}(h_{k-1})\big),
\end{eqnarray}
where the indice $h_{k-1}$ is in the set $\{0,\ldots,i-k+1\}$ which is empty if $i\leqslant k-1$.

Let $i\leqslant \tilde{s}$. According to Equation (\ref{lemme7.eq}) we get $\hat{p}(i)-\bar{p}(i)\geq 0$ because $\bar{q}_{k-1}(h_{k-1})=q_{k-1}(h_{k-1})$ for all $h_{k-1}\in\{0,\ldots,\tilde{s}-k+1\}$. Then the point $ii)$ is true.
\\

Let $i=i_0+k-1$. Noting that $q_{k-1}(h_{k-1})=\bar{q}_{k-1}(h_{k-1})$ except in $h_{k-1}=i_0$ we get
\begin{eqnarray*}
\hat{p}(i)-\bar{p}(i)=(-1)^{k-1}(q_{k-1}(i_0)-\bar{q}_{k-1}(i_0))=+\epsilon \text{ (because $k$ is odd).}
\end{eqnarray*}
and point $iii)$ is shown.
\\

It remains to show that $\bar{p}\leqslant \hat{p}$. By construction of $\bar{q}_{k-1}$, the primitive of $q_{k-1}$ is greater than the primitive of $\bar{q}_{k-1}$, and because $\hat{p}(i)-\bar{p}(i)$ is nonnegative and the following equality:
\begin{eqnarray*}
\hat{p}(i)-\bar{p}(i)=\sum_{h_1=0}^{i-1}\ldots\sum_{h_{k-2}=0}^{h_{k-3}-1}\left(\sum_{h_{k-1}=0}^{h_{k-2}-1}q_{k-1}(h_{k-1})-\sum_{h_{k-1}=0}^{h_{k-2}-1}\bar{q}_{k-1}(h_{k-1})\right).
\end{eqnarray*}
we get point $i)$.
\\

If $k$ is even the proof is based on another construction of $\bar{q}_{k-1}$. Let us first recall that $i$ is a $1-$knot of $q_{k-1}$ if $\Delta^1q_{k-1}(i)=q_{k-1}(i+1)-q_{k-1}(i)$ is strictly negative (because $k$ is even). We have two cases:\\
Case 1 : There exists $(i_0,i_1)$ such that $\tilde{s}\leq i_0<i_1$, $i_1-i_0\geq 2$, $\Delta^1q_{k-1}(i_0)$ and $\Delta^1 q_{k-1}(i_1)$ are strictly negative and $\Delta^1q_{k-1}(i)=0$. The probability $\bar{q}_{k-1}$ is defined as follows:
\begin{eqnarray*}
\left\{
    \begin{array}{l}
    \vspace{0.25cm}
    \bar{q}_{k-1}(i)=q_{k-1}(i) \text{   if $i<i_0+1,i> i_1$}\\
    \vspace{0.25cm}
    \bar{q}_{k-1}(i_0+1)=q_{k-1}(i_0)+\epsilon\\
    \vspace{0.25cm}
    \bar{q}_{k-1}(i_1)=q_{k-1}(i_1)-\epsilon\\
    \vspace{0.25cm}
    \text{$\bar{q}$ is an affine function on $[i_0+1,i_1]$}.
    \end{array}
  \right.
\end{eqnarray*}

Case 2 : For all $i\geq \tilde{s}+1$, $\Delta^1q_{k-1}(i)<0$. let $i_0=\tilde{s}+1$, then the probability $\bar{q}_{k-1}$ is defined as follows:
\begin{eqnarray*}
\left\{
    \begin{array}{l}
    \vspace{0.25cm}
    \bar{q}_{k-1}(i)=q_{k-1}(i) \text{   if $i<i_0+1,i> i_0+2$}\\
    \vspace{0.25cm}
    \bar{q}_{k-1}(i_0+1)=q_{k-1}(i_0)+\epsilon\\
    \vspace{0.25cm}
    \bar{q}_{k-1}(i_0+2)=q_{k-1}(i_0+2)-\epsilon.
    \end{array}
  \right.
\end{eqnarray*}

\begin{figure}[!h]
\centering
\begin{subfigure}[b]{0.45\textwidth}
\includegraphics[width=0.95\textwidth]{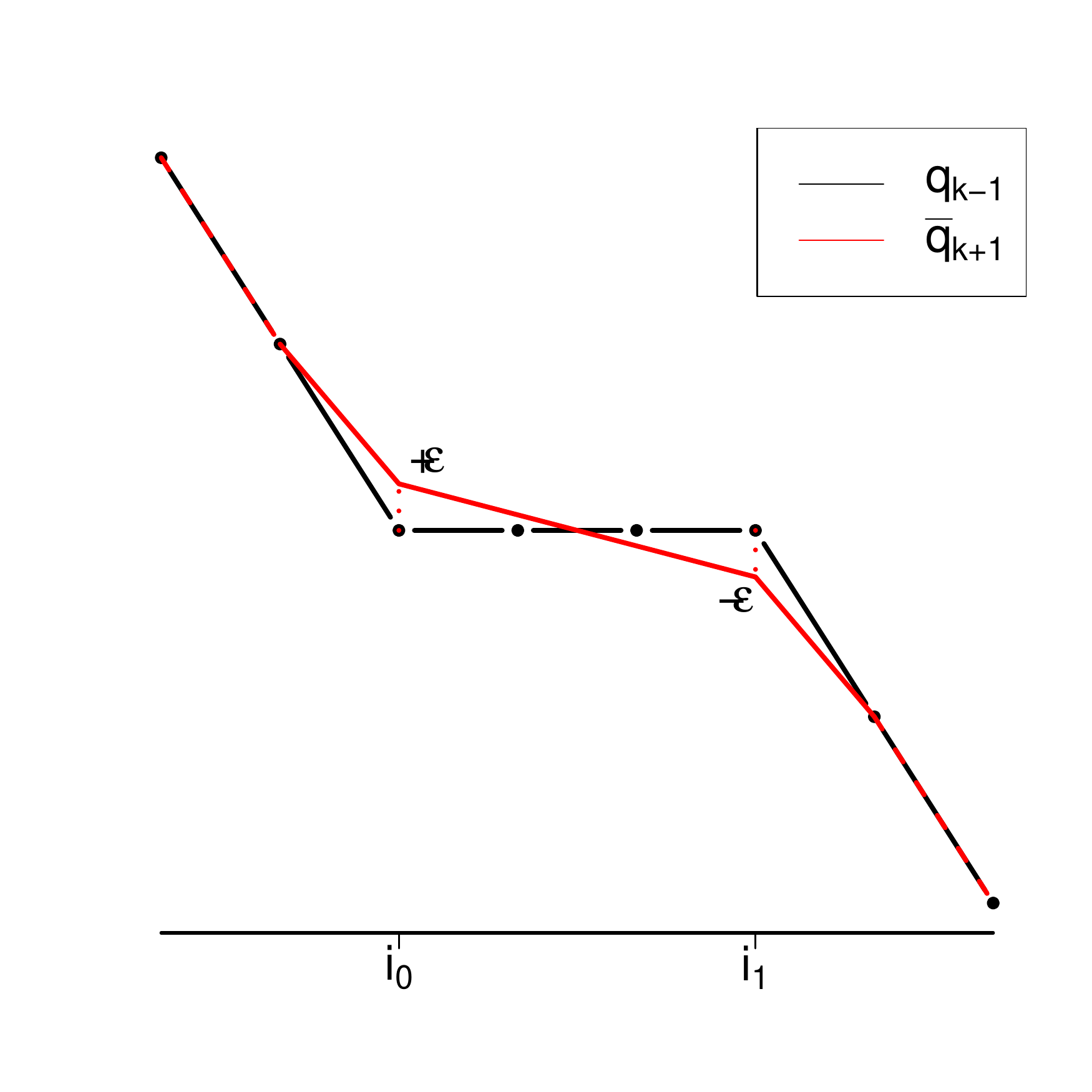}
\caption{Case 1.}
\label{fig2a}
\end{subfigure}
\hfill
\begin{subfigure}[b]{0.45\textwidth}
\includegraphics[width=0.95\textwidth]{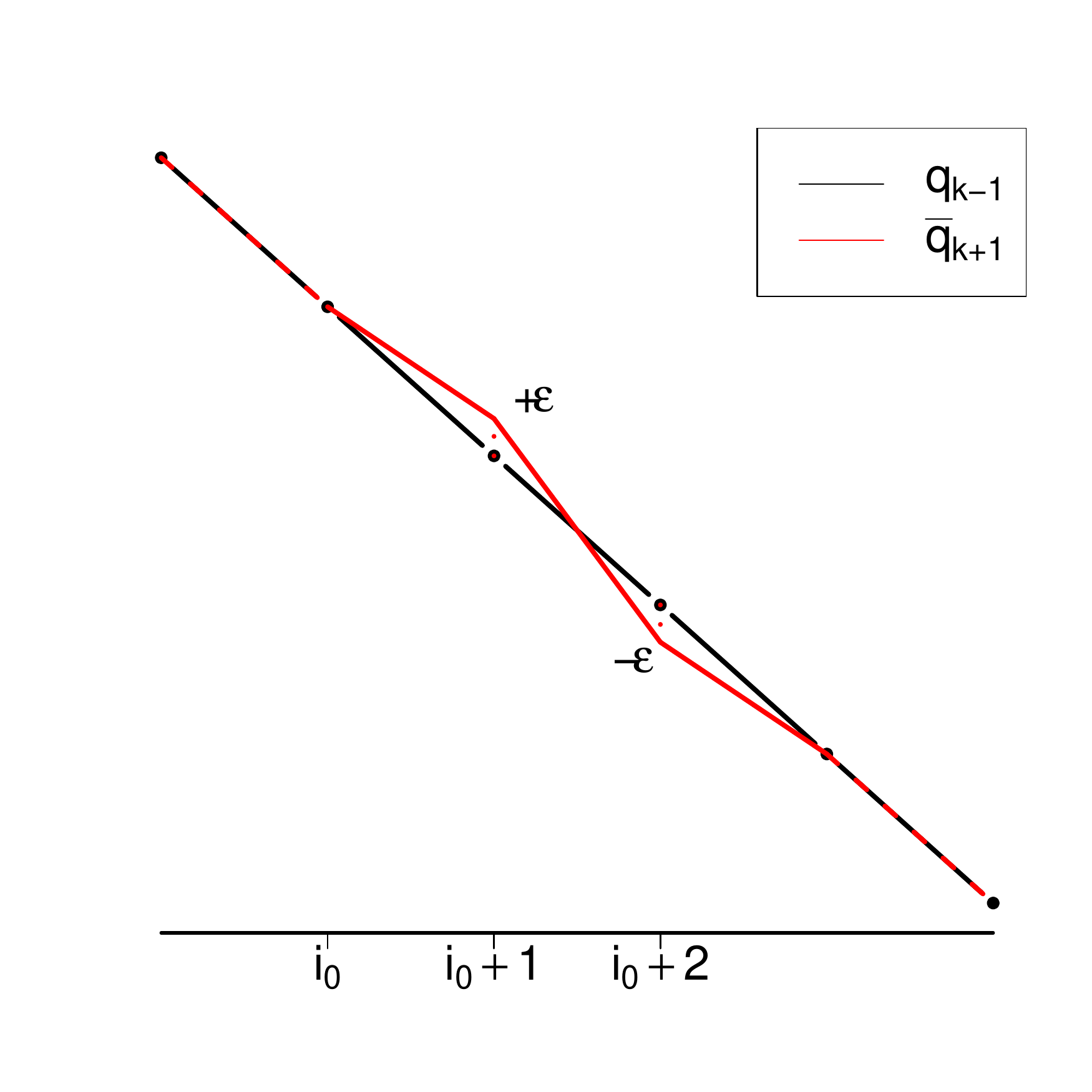}
\caption{Case 2.}
\label{fig2b}
\end{subfigure}
\caption{Functions $q_{k-1}(i)$ and $\bar{q}_{k-1}(i)$ versus i ($k$ even).}
\label{fig2}
\end{figure}

The functions $\bar{q}_{k-1}$ are presented in Figure \ref{fig2}. The rest of the proof is similar to the one when $k$ is odd.

Therefore, Theorem \ref{supportfini.prop} is proved in the case $\beta(\hat{p})=0$. Assume now that $\beta(\hat{p})\neq 0$. By Theorem \ref{theo2} we know that if $l$ is a $k-$knot of $\hat{p}$, $\beta(\hat{p})$ is written as follows:
\begin{eqnarray}
\label{aide33}
\frac{F^{k}_{\hat{p}}(l)- F^{k}_{\tilde{p}}(l)}{m^k_l}= \beta(\hat{p}).
\end{eqnarray} 
Let us proved that if the support of $\bar{p}$ is infinite then $\beta(\hat{p})=0$. Indeed if the support of $\bar{p}$ is infinite, $\bar{p}$ has an infinite numbers of $k-$knots and Equation (\ref{aide33}) is true for an infinite numbers of integers $l$.\\
Moreover by Equation~(\ref{massespline}) the term $m^k_l$ is a polynomial function in the variable $l$ with degree $k$ and by Lemma \ref{lemme16} (see Section \ref{prooflemmas.st}) the term $F^{k}_{\hat{p}}(l)- F^{k}_{\tilde{p}}(l)$ is a polynomial function with degree less than $k-1$. Therefore  the left side in Equation (\ref{aide33}) tends to zero when $l$ tends to infinity, showing that $\beta(\hat{p})= 0$.

\paragraph{$k$-knots' repartition beyond $\widetilde{s} -k +2$.}

Let us assume that $k$ is odd and prove that if $\hat{s}\geq \tilde{s}+1$ then $\Delta^k\hat{p}(r)=0$ for all $r\in[\tilde{s}-k+2,\hat{s}-2]$. We consider $q_1,\hdots,q_{k-1}$ the derivative functions of $\hat{p}$ defined as before in Equation (\ref{defder}).

As $\hat{s}$ is a $k-$knot of $\hat{p}$ there exist two consecutive $1-$knots between $r$ and $\tilde{s}$. This allows to define the function $\bar{q}_{k-1}$ and $\bar{p}$ as before in Equation (\ref{aux73}) and Equation (\ref{aux74}).

By construction $\bar{q}_{k-1}$ is non-increasing (and nonnegative) and therefore $\bar{p}$ is $k-$monotone (and nonnegative). Moreover $\bar{p}$ is lower than $\hat{p}$, equal to $\hat{p}$ on $\{0,\ldots,\tilde{s}\}$ and for $i=r+k-1$ we have $\bar{p}(i)<\hat{p}(i)$. Moreover $r+k-1>\tilde{s}$ because $r\in\{\tilde{s}-k+2$,\ldots,$\hat{s}-1\}$. It follows that $\vert\vert\bar{p}-\tilde{p}\vert\vert_2<\vert\vert\hat{p}
-\tilde{p}\vert\vert_2$ which contradicts the definition of $\hat{p}$. Therefore $\hat{p}$ does not have any $k-$knot on $\{\tilde{s}-k+2$,\ldots,$\hat{s}-1\}$.
\\

The proof is similar when $k$ is even.

\begin{Rem}The second case requires $q_{k-1}(r+2)>0$ for $\bar{q}$ to be nonnegative. That is to say we need that $r+2\leqslant \hat{s}$. This is the reason why the two sets $\{\tilde{s}-k+2$,\ldots,$\hat{s}-2\}$ and $\{\tilde{s}-k+2$,\ldots,$\hat{s}-1\}$ are different if $k$ is odd or $k$ is even.
\end{Rem}

\subsubsection{Proof of Theorem \ref{moments} : Comparison between the moments of $\tilde{p}$ and the moments of $\hat{p}$.}

Let $q$ a sequence and let $\epsilon$ a real number such that $(1-\epsilon)\hat{p}+\epsilon q$ is a $k-$monotone probability. Because $\hat{p}$ minimizes $\mathcal{Q}$ over the set of $k-$monotone probabilities we get:
\begin{eqnarray*}
\liminf_{\epsilon\searrow 0}\frac{1}{\epsilon}(\mathcal{Q}((1-\epsilon)\hat{p}+\epsilon q)-\mathcal{Q}(\hat{p}))\geqslant 0,
\end{eqnarray*}
that is equivalent to:
\begin{eqnarray*}
\liminf_{\epsilon\searrow 0}\frac{1}{\epsilon}\left(\sum_{i=0}^{\infty}\big((1-\epsilon)\hat{p}(i)
+\epsilon q(i)-\tilde{p}(i)\big)^2-\sum_{i=0}^{\infty}(\hat{p}(i)-\tilde{p}(i))^2\right)\geqslant 0,
\end{eqnarray*}
which leads after simplifications (see the proof of Theorem \ref{theo2} for more explanations) to:
\begin{eqnarray*}
\sum_{i=0}^{\infty}(\hat{p}(i)-\tilde{p}(i))q(i)\geqslant \beta(\hat{p}).
\end{eqnarray*}
For $q(i)=\vert i-a\vert^u_+/m(a,u)$ we get the result. Moreover for $q=\delta_0$ we find $\hat{p}(0)-\tilde{p}(0)\geqslant \beta(\hat{p})$.

\subsubsection{Proof of Theorem \ref{vitesse} : Rate of convergence of $\hat{p}$.}

The proof is based on Lemma 6.2 of Jankowski and Wellner (2009) \cite{jankowski2009estimation}. First we assume that $r=2$. Banach's Theorem for projection on a closed convex set says that the projection on the set of $k-$monotone probabilities is $1-$lipschitzienne. Then if $p_{\mathcal{S}_k}$ is the projection of $p$ on the set $\mathcal{S}_k$ we have:
\begin{eqnarray*}
\sqrt{n}\vert\vert p_{\mathcal{S}_k}-\hat{p}\vert\vert_2\leqslant \sqrt{n}\vert\vert p-\tilde{p}\vert\vert_2
\end{eqnarray*}
We need to show that $\sqrt{n}\vert\vert p-\tilde{p}\vert\vert_2=O_\mathbb{P}(1)$, or equivalently that the series $W_n=\sqrt{n}(p-\tilde{p})$ is tight in $L_2(\mathbb{N})$. Using Lemma 6.2 of Jankowski and Wellner (2009), we have to show that:
\begin{eqnarray*}
\left\{
    \begin{array}{l}
    \vspace{0.25cm}
     \underset{n\in\mathbb{N}}{\sup }\text{ }\mathbb{ E}\big[\vert\vert W_n\vert\vert^2_2\big]<\infty,\\ 
    \underset{m\rightarrow\infty}{\lim}\underset{n\in\mathbb{N}}{\sup }\sum_{j\geqslant m}\mathbb{E}\big[\vert\vert W_n\vert\vert_2\big]=0.
\end{array}
\right.
\end{eqnarray*}
This is easily verified by noting that $\text{var}(\tilde{p}(j))=p(j)(1-p(j))/n$. Then for all $r\in[2,\infty]$, $\sqrt{n}\vert\vert p_{S_k}-\hat{p}\vert\vert_r\leqslant \sqrt{n}\vert\vert p_{S_k}-\hat{p}\vert\vert_2=O_\mathbb{P}(1)$.

\subsubsection{Proof of Theorem \ref{estimnoeud.prop} : The case of a finite support for $p$. \label{estimnoeud.proof}}

\paragraph{First part}

For all integer $i$, by the strong law of large numbers $\tilde{p}(i)$ tends a.s. to $p(i)$. Because the maximum of the support $s$ of $p$ is finite we have the following result:
\begin{eqnarray*}
a.s. \underset{n\rightarrow\infty}{\lim}\vert\vert\tilde{p}-p\vert\vert^2_2=\underset{n\rightarrow\infty}{\lim}\sum_{i=0}^{s}(\tilde{p}(i)-p(j))^2=0.
\end{eqnarray*}
Then by Theorem \ref{limnorme2.prop} we get that for all integer $i$:
\begin{eqnarray*}
a.s.\underset{n\rightarrow\infty}{\lim}(\hat{p}(i)-p(i))^2\leqslant\underset{n\rightarrow\infty}{\lim} \vert\vert\hat{p}-p\vert\vert^2_2\leqslant \underset{n\rightarrow\infty}{\lim}\vert\vert\tilde{p}-p\vert\vert^2_2=0.
\end{eqnarray*}

It follows that:
\begin{eqnarray*}
a.s.\underset{n\rightarrow\infty}{\lim}\left[(-1)^k\Delta^k\hat{p}(j)-(-1)^k\Delta^kp(j)\right]=0.
\end{eqnarray*}
Because $(-1)^k\Delta^kp(j)>0$, almost surely for $n$ large enough $(-1)^k\Delta^k\hat{p}(j)>0$, which proves that $j$ is a $k-$knot of $\hat{p}$.

\paragraph{Second part}

If $\hat{s}\leqslant s$ the theorem is true. We assume now that $\hat{s}> s$.\\
Let us first consider the case where $k$ is odd. 
Thanks to the second point of Theorem~\ref{supportfini.prop}, if we note $\tilde{s}$ the maximum of the support of $\tilde{p}$ then $\hat{p}$ has no $k-$knot on $\{\tilde{s}-k+2,\ldots,\hat{s}-2\}$.\\
Moreover as $\tilde{s}\leqslant s$, $\hat{p}$ has no $k-$knot in $\{s-k+2,\ldots,\hat{s}-2\}$ (this set may be empty).\\
The function $p$ is $k-$monotone and $s$ is a $k-$knot of $p$, then by Theorem \ref{estimnoeud.prop} almost surely there exists $n_0$ such as for all $n\geqslant n_0$, $s$ is a $k-$knot of $\hat{p}$.\\
It follows that (almost surely) $s$ is not in the set $\{s-k+2,\ldots,\hat{s}-2\}$ and therefore $s\geqslant \hat{s}-1$ or $\hat{s}\leqslant s+1$.\\
The proof of the result in the case where $k$ is even is similar.

\subsubsection{Proof of Theorem \ref{critere_k=34}\label{critere_k=34.st} : Stopping criterion when $k\in\{3,4\}$}

We first show that $\hat{p}$ satisfies the four properties stated in \textit{1}. We know by Theorem \ref{theo2} that it satisfies \textit{1.(a)} and \textit{1.(b)}. Moreover by Lemma \ref{signebeta} (see Section \ref{prooflemmas.st}) it satisfies \textit{1.(d)}. It remains to show 1.(c).\\

The proof is similar to the proof of Theorem \ref{supportfini.prop}. For $\epsilon$ a real number, and for any $j\in\{2,k-1\}$, the function $q_{\epsilon}$ is defined as follows:
\begin{equation*}
 q_{\epsilon}(i)=(1-\epsilon)\hat{p}(i)+\epsilon \frac{\bar{Q}^{j}_{\hat{s}+1}(i)}{m^j_{\hat{s}+1}}
\end{equation*}
where $\bar{Q}_{s+1}^j$ is defined at Equation (\ref{splinedef}).\\
The function $q_{\epsilon}$ is a $k-$monotone probability for $\epsilon$ small enough. Indeed $(-1)^k\Delta^k \bar{Q}^j_{\hat{s}+1}(i)$ is strictly nonpositive only for $i=\hat{s}$. Moreover $(-1)^k\Delta^k \hat{p}(\hat{s})=\hat{p}(\hat{s})>0$ then for $\epsilon$ smaller enough, $(-1)^k\Delta^k q_{\epsilon}(\hat{s})=(1-\epsilon)\hat{p}(\hat{s})-\epsilon Q^{j}_{\hat{s}+1}(i)/m^j_{\hat{s}+1}$ is nonnegative.\\
On the other hand as $\hat{p}$ minimizes $\mathcal{Q}$, we get:
\begin{eqnarray*}
\liminf_{\epsilon \searrow 0}\frac{1}{\epsilon}(\mathcal{Q}(q_{\epsilon})-\mathcal{Q}(\hat{p}))\geqslant 0,
\end{eqnarray*}
which is equivalent to:
\begin{eqnarray*}
\liminf_{\epsilon\searrow 0}\frac{1}{\epsilon}\left[\epsilon^2\sum_{i=0}^{\infty}\left(\frac{1}{m^j_{\hat{s}+1}}Q^{j}_{\hat{s}+1}(i)-\hat{p}(i)\right)^2+2\epsilon\sum_{i=0}^{\infty}\left(\hat{p}(i)-\tilde{p}(i)\right)\left(\frac{Q^{j}_{\hat{s}+1}(i)}{m^j_{\hat{s}+1}}-\hat{p}(i)\right)\right]\geqslant 0.
\end{eqnarray*}
This leads to the following inequality:
\begin{eqnarray*}
\sum_{i=0}^{\hat{s}+1}(\hat{p}(i)-\tilde{p}(i))\frac{Q^{j}_{\hat{s}+1}(i)}{m^j_{\hat{s}+1}}- \sum_{i=0}^{\infty}\hat{p}(i)(\hat{p}(i)-\tilde{p}(i))\geqslant 0.
\end{eqnarray*}
By Lemma \ref{rec.lemme} (see Section \ref{prooflemmas.st}) we deduce that:
\begin{eqnarray*}
\frac{F^j_{\hat{p}}(\hat{s}+1)-F^j_{\tilde{p}}(\hat{s}+1)}{m^j_{\hat{s}+1}}\geqslant \beta(\hat{p}),
\end{eqnarray*}
which is exactly \textit{2.(c)}.
\\

Reciprocally we assume now that $f$ satisfies \textit{1}. and we show that $f=\hat{p}$, which by Theorem \ref{theo2}, is equivalent to show that
\begin{eqnarray*}
F^k_{f}(l)-F^k_{\tilde{p}}(l)\geqslant\beta(f)m^k_l,
\end{eqnarray*}
for all $l\in\mathbb{N}$ with equality if $l$ is a $k-$knot of $f$.\\
This is true for $l\leq s+1$ because $f$ satisfies \textit{2.(a)} and \textit{2.(b)}. Because $f$ has no $k-$knot after $s$ it remains to show that the inequality is true for $l\geqslant s+2$. \\
We begin with the case $k=3$. Because $s\geqslant \tilde{s}$ and $f$ and $\tilde{p}$ are probabilities, applying Theorem \ref{theoclé}, we obtain that for all $l\geqslant s+1$,
\begin{align*}
F^3_{f}(l)-F^3_{\tilde{p}}(l)&=\sum_{j=2}^3Q_{l-1}^{3-j+1}(s+1)\big(F^j_{f}(s+1)-F^j_{\tilde{p}}(s+1)\big)\\
&=\big(F^3_{f}(s+1)-F^3_{\tilde{p}}(s+1)\big)+(l-s-1)\big(F^2_{f}(s+1)-F^2_{\tilde{p}}(s+1)\big).
\end{align*}
As $f$ satisfies \textit{1.(a)} we have:
\begin{eqnarray*}
F^3_{f}(l)-F^3_{\tilde{p}}(l)\geqslant\beta(f)m^3_{s+1}+(l-s-1)\big(F^2_{f}(s+1)-F^2_{\tilde{p}}(s+1)\big).
\end{eqnarray*}
Moreover as $f$ satisfies \textit{1.(c)} we have:
\begin{eqnarray*}
F^3_{f}(l)-F^3_{\tilde{p}}(l)\geqslant\beta(f)m^3_{s+1}+(l-s-1)\beta(f)m^2_{s+1}=\beta(f)\big(m^3_{s+1}+(l-s-1)m^2_{s+1}\big).
\end{eqnarray*}
Finally, because $\beta(f)\leqslant 0$ by \textit{1.(d)}, it remains to show that:
\begin{eqnarray}\label{aux8}
m^3_{s+1}+(l-s-1)m^2_{s+1}\leqslant m^3_{l}.
\end{eqnarray}
By Equation~\eref{massespline}, Equation \eref{aux8} may be written as follows:
\begin{eqnarray*}
\frac{(s+4)(s+3)(s+2)}{6}+(l-s-1)\frac{(s+3)(s+2)}{2}\leqslant \frac{(l+3)(l+2)(l+1)}{6}.
\end{eqnarray*}
After some calculations, we can show that \eref{aux8} is satisfied if and only if $ P_3(l)\geqslant 0$ where $P_3$ is the polynomial function $P_3(l)=(l-s)(l-(s+1))(l+2s+7)$. This is true because $l\geqslant s+1$.\\
Let us now prove the case $k=4$. By Theorem \ref{theoclé} we obtain for all $l\geqslant s+1$:
\begin{align*}
F^4_{f}(l)-F^4_{\tilde{p}}(l)
&=\big(F^4_{f}(s+1)-F^4_{\tilde{p}}(s+1)\big)+(l-s-1)\big(F^3_{f}(s+1)-F^3_{\tilde{p}}(s+1)\big)\\
&+Q_{l-1}^3(s+1)\big(F^2_{f}(s+1)-F^2_{\tilde{p}}(s+1)\big).
\end{align*}
Let $A_I(s)=\Pi_{i\in I}(s+i)$, then using Equation~\eref{massespline} we need to show that:
\begin{eqnarray}
\label{aux35}
A_{[3,5]}(s)+4A_{[3,4]}(s)(l-s-s)+6(l-s-1)(l-s-2)A_{\{3\}}(s)
\leqslant \frac{(l+1)A_{[1,4]}(l)}{s+2}.
\end{eqnarray}
After some calculations, we can show that \eref{aux35} is satisfied if and only if $P_4(l)\geqslant 0$ where $P_4$ is the polynomial function
\begin{eqnarray*}
P_4(l)=A_{[1,4]}(l)-6(l-(s+1))(l-(s+2))A_{[2,3]}(s)
+4(l-(s+1))A_{[2,4]}(s)-A_{[2,5]}(s).
\end{eqnarray*}
We have $P_4(l+1)-P_4(l)=4(P_3(l+1)+3(s+2)(s+3))$ and $P_4(s+2)=12(s+3)(s+4)>0$ then $P_4(l)\geqslant 0$ because $l\geqslant s+2$.

\subsection{\label{algo2.st} Estimating $\pi$ on a finite support}
\begin{Theo}The algorithm described at Table~\ref{algo.tb} returns $\hat{p}_{L}$ in a finite number of steps.\label{algo2}
\end{Theo}

\subsubsection{Proof of Theorem \ref{algo2}}

During step 1 the set $S$ is a subset of $\{0,\ldots,L\}$ and $\pi$ is the minimizer of $\Psi$ on the set $S$.\\ 
The criterion allowing us to determine if $\pi=\hat{\pi}_L$ (and to stop the algorithm) is given by Lemma \ref{lemme3} (see Section \ref{prooflemma}).
\\

In order to show that the algorithm returns $\hat{\pi}_{L}$ in a finite number of steps we need to show the both assertions :
\begin{itemize}
\item \textbf{Assertion 1} : Going from Step 2 to Step 1 is done in a finite number of runs.
\item \textbf{Assertion 2} : If $\pi_m$ denotes the value of $\pi$ at iteration $m$ of the algorithm, then $(\Psi(\pi_m))$ converges to the minimum of $\Psi$ on the set of probabilities with support on $\{0,\ldots,L\}$ that is to say to $\hat{\pi}_L$.
\end{itemize}
At Step 2 the set $S'$ may be reduced up to one element, but it can not be empty because the minimizer of $\Psi$ on a singleton is non-negative. That proves Assertion 1.\\
Let us show Assertion 2 by proving that for all
%During Step 2, each time we enter on Case 2 we remove an element to $S'$. By the lemma \ref{lemme2} the first time we enter on Case 2 this element is not the integer $j$ added on the end of Step 1 (Because 2. assure that $b_j=\pi_{S'}(j)\geqslant 0=a_j$). Then we do not return on Step 1 with the same $S$ and do not curl.
$m\in\mathbb{N}^*$, $\Psi(\pi_{m+1})<\Psi(\pi_m)$. Let $S$ be the support of $\pi_m$ at iteration $m$, and let $j\in\{0,\ldots,L\}$ be an integer such as $D_{\delta_j}\Psi(\pi_m)<0$. We have $S'=S+{j}$ and  $\Psi(\pi_{S'})<\Psi(\pi_{S})$ by Lemma \ref{lemme2} (see Section \ref{prooflemma}).

We consider two cases : 
\\

\textbf{1 }: If $\pi_{S'}$ is a nonnegative measure we go to Step 1 with $\pi=\pi_{S'}$. In other terms $\pi_{m+1}=\pi_{S'}$ and therefore $\Psi(\pi_{m+1})<\Psi(\pi_m)=\Psi(\pi_S)$.
\\

\textbf{2 }: If $\pi_{S'}$ is not a nonnegative measure the algorithm iterates inside Step 2 and $\pi_{S'}$ is updated at each loop. We need to verify that at the end of this iterative procedure:
\begin{eqnarray*}
\Psi(\pi_{S''})<\Psi(\pi_S).
\end{eqnarray*}
Let $r$ be the number of times when we go in Step 2 during the $m$-th loop and let $S''_h$ be the value of the set $S''$ the $h$-th time we go in Step 2. We have $\pi_{m+1}=\pi_{S''_h}$.\\
We show by induction the following property :
\begin{eqnarray*}
HR_h : \Psi(\pi_{S''_h})<\Psi(\pi_{S}).
\end{eqnarray*}

Thanks to the property \textit{2.} in Lemma \ref{lemme2} the property $HR_1$ is true. Assume now that $HR_h$ is true for some $h\leqslant r-1$, $\Psi(\pi_{S''_h})<\Psi(\pi_S)$. Let $l$ and $\epsilon$ be defined as follows:
\begin{eqnarray*}
l = \underset{i\in S'}{\text{argmin}}\{\frac{a_i}{a_i-\pi_{S''_h}(i)}, \text{pour } i, \pi_{S''_h}(i)<a_i\},
\end{eqnarray*}
\begin{eqnarray*}
\epsilon=\frac{a_l}{a_l-\pi_{S''_h}(l)}.
\end{eqnarray*}
Then $(1-\epsilon)\pi_S+\epsilon\pi_{S''_h}$ is a $1$-mass function with support $S''_{h+1}=S''_{h}-\{l\}$. It follows, by convexity of $\Psi$ that:
\begin{align*}
\Psi(\pi_{S''_{h+1}})&\leqslant \Psi((1-\epsilon)\pi_S+\epsilon\pi_{S''_h})\\
&\leqslant(1-\epsilon)\Psi(\pi_S)+\epsilon\Psi(\pi_{S''_h}).
\end{align*}
Thanks to $HR_h$, it follows that:
\begin{eqnarray*}
\Psi(\pi_{S''_{h+1}})<\Psi(\pi_S),
\end{eqnarray*} 
and $HR_{h+1}$ is true.\\
Then $HR_r$ is true, that is to say $\Psi(\pi_{m+1})<\Psi(\pi_m)$ for all integer $m$, and $(\Psi(\pi_m))_{m\in\mathbb{N}})$ converges when $m$ tends to infinity (because it is a nonincreasing and bounded sequence).
The limit is the minimum of $\Psi$ because the nondecreasing is strict.

\subsubsection{Proof of the lemmas \label{prooflemma}}

The proof of Theorem \ref{algo2} is based on the following lemmas whose proofs are given afterwards. All the notations used in this section were defined in Section \ref{algorithm.section}.

\begin{Lemme}\label{lemme6}Let $\pi$ and $\mu$ be two probabilities with support on the set $\{0,\ldots,L\}$. Then we have the following equality:
\begin{eqnarray*}
D_{\mu}\Psi(\pi)=\sum_{j=0}^L\mu(j)D_{\delta_j}\Psi(\pi).
\end{eqnarray*}
\end{Lemme}

\begin{Lemme}\label{lemme3}
There is equivalence between :
\begin{enumerate}
\item $\pi=\hat{\pi}_{L}$.
\item $\forall j\in\{0,\ldots,L\}$, $D_{\delta_j}\Psi(\pi)\geqslant 0$.
\end{enumerate}
Moreover if $\pi=\hat{\pi}_{L}$ then for all $j\in$ supp($\pi$) we have $D_{\delta_j}\Psi(\pi)= 0$.
\end{Lemme}

\begin{Lemme}
\label{lemme5}
Let $\mathcal{M}_S$ be the set of positive measure $\pi$ whose support is included in the set $S$. Let $\pi_S$ and $\pi_{S'}$ be defined as follows:
\begin{eqnarray*}
\pi_S=\underset{\pi\in\mathcal{M}_S}{\underset{\sum_{j\in S}\pi(j)=1}{\text{argmin}}}\big(\Psi(\pi)\big),
\end{eqnarray*}
\begin{eqnarray*}
\pi'_S=\underset{\pi\in\mathcal{M}_S}{\text{argmin}}\left(\Psi(\pi)+\lambda_S(\sum_{j\in S}\pi(j)-1)\right).
\end{eqnarray*}
Then we have $\pi_S=\pi'_S$.
\end{Lemme}

The proof of the following lemma is in Durot and al. \cite{durot_least}:
\begin{Lemme}\label{lemme2}Let $\pi_{S}=\sum_{i\in S}a_i\delta_{i}$ be the minimizer of $\Psi$ over the set of nonnegative sequences with support $S\subset\{0,\ldots,L\}$.\\
Let $j$ an integer such that $j\notin S$ and $D_{\delta_j}\Psi(\pi_L)<0$.\\
Let $\pi^*_{S'}=\sum_{i\in S'}b_i\delta_{i}$ be the minimizer of $\Psi$ over the set of sequences with support $S'=S+\{j\}$.\\
Then, the two following results hold:
\begin{enumerate}
\item $\Psi(\pi_{S'})<\Psi(\pi_S)$.
\item Assume that $b_i$ for some $i\in S$ is strictly nonpositive and let:
\begin{eqnarray*}
l = \underset{1\in S'}{\text{argmin}}\{\frac{a_i}{a_i-b_i}, \text{pour } i, b_i<a_i\}.
\end{eqnarray*}
If $\pi_{S''}$ is the minimizer of $\Psi$ over the set of sequences with support $S''=S'-\{l\}$, then $\Psi(\pi_{S''})<\Psi(\pi_S)$.
\end{enumerate}
\end{Lemme}

\paragraph{Proof of Lemma \ref{lemme6} }
Let $\mu$ be a probability with support included in $\{0,\ldots,L\}$. We write $\mu=\sum_{j=0}^L\mu(j)\delta_j$ then, for $L\leqslant \tilde{s}$:
\begin{align*}
D_{\mu}\Psi(\pi)=&\underset{\epsilon\searrow 0^+}{\lim}\frac{1}{\epsilon}\big(\Psi \big((1-\epsilon)\pi+\epsilon\mu\big)-\Psi(\pi)\big)\\
=&\underset{\epsilon\searrow 0^+}{\lim}\frac{1}{\epsilon}\left(\sum_{i=0}^{L}\left(\sum_{l=0}^L[(1-\epsilon)\pi(l)+\epsilon
\mu(l)]Q_l^k(i))-\tilde{p}(i)\right)^2\right.\\
&-\left.\sum_{i=0}^{L}\left(\sum_{l=0}^L\pi(l)Q_l^k(i))
-\tilde{p}(i)\right)^2\right)\\
=&\underset{\epsilon\searrow 0^+}{\lim}\frac{1}{\epsilon}\sum_{i=0}^{L}\left[2\epsilon\left(\sum_{l=0}^L(\mu(l)-\pi(l))
Q_l^k(i)\right)\left(\sum_{l=0}^L\pi(l)Q_l^k(i)
-\tilde{p}(i)\right)\right]\\
&\left.+\epsilon^2\left(\sum_{l=0}^L(\mu(l)-\pi(l))
Q_l^k(i)\right)^2\right]\\
=&2\sum_{i=0}^{L}\left(\sum_{l=0}^L\mu(l)Q^k_l(i)-\sum_{l=0}^L\pi(l)Q^k_l(i)\right)
\left(\sum_{l=0}^L\pi(l)Q^k_l(i)-\tilde{p}(i)\right).
%=& \sum_{l=0}^L\mu(l)\left(\sum_{i=0}^{L}\big(Q^k_l(i)-
%\sum_{l=i}^{L}\pi(l)Q^k_{l}(i)\big)\right)
%\left(\sum_{j=i}^{L}\pi(j)Q^k_{j}(i)-\tilde{p}(i)\right)\\
%=&\sum_{j=0}^L\mu(j)D_{\delta_j}\Psi(\pi).
\end{align*}

In particular for $\mu=\delta_j$ we find:
\begin{eqnarray*}
D_{\delta_j}\Psi(\pi)=2\sum_{i=0}^{L}\left(Q^k_j(i)-
\sum_{l=0}^L\pi(l)Q^k_l(i)\right)
\left(\sum_{l=i}^L\pi(l)Q^k_l(i)-\tilde{p}(i)\right).
\end{eqnarray*}

Then, by noting that $\sum_j\mu(j)=1$ we have the following equalities:
\begin{align*}
\sum_{j=0}^{L}\mu(j)D_{\delta_j}\Psi(\pi)&=2\sum_{i=0}^{L}
\left(\sum_{j=0}^{L}\mu(j)\left(Q^k_j(i)-
\sum_{l=0}^L\pi(l)Q^k_l(i)\right)\right)
\left(\sum_{l=i}^L\pi(l)Q^k_l(i)-\tilde{p}(i)\right)\\
&=2\sum_{i=0}^{L}
\left(\sum_{j=0}^{L}\mu(j)Q^k_j(i)-
\sum_{l=0}^L\pi(l)Q^k_l(i)\right)
\left(\sum_{l=i}^L\pi(l)Q^k_l(i)-\tilde{p}(i)\right)
\end{align*}
and the lemma is proved.

\paragraph{Proof of Lemma \ref{lemme3}}

We first show that $\hat{\pi}_{L}$ satisfies \textit{2.}\\
For all $0<\epsilon<1$ and $j\in\{0,\ldots,L\}$ the function $(1-\epsilon)\hat{\pi}_{L}+\epsilon\delta_j$ is a probability and then by definition of $\hat{\pi}_L$ we have the following inequality:
\begin{eqnarray*}
\underset{\epsilon\searrow 0^+}{\lim}\frac{1}{\epsilon}\big(\Psi ((1-\epsilon)\hat{\pi}_L+\epsilon\delta_j)-\Psi (\hat{\pi}_L)\big)\geqslant 0,
\end{eqnarray*}
which leads to $D_{\delta_j}\Psi (\hat{\pi}_L)\geqslant 0
$, showing the point \textit{2.}
\\

Reciprocally, for $\pi$ a probability that satisfies \textit{2.}, let us show that $\pi=\hat{\pi}_{L}$. Precisely we show that for all probability $\mu$ with support on $\{0,\ldots,L\}$ we have $\Psi(\mu)-\Psi(\pi)\geqslant 0$. Because $\Psi$ is convex we have:
\begin{align*}
D_{\mu}\Psi(\pi)&=\underset{\epsilon\searrow 0^+}{\lim}\frac{1}{\epsilon}\left(\Psi \big((1-\epsilon)\pi+\epsilon\mu\big)-\Psi (\pi)\right)\\
&\leqslant \underset{\epsilon\searrow 0^+}{\lim}\frac{1}{\epsilon}\big((1-\epsilon)\Psi (\pi)+\epsilon\Psi(\mu))-\Psi (\pi)\big)\\
&\leqslant \Psi(\mu)-\Psi(\pi),
\end{align*}
and by Lemma \ref{lemme6} we have:
\begin{eqnarray*}
D_{\mu}\Psi(\pi)=\sum_{j=0}^L\mu(j)D_{\delta_j}\Psi(\pi).
\end{eqnarray*}
%Comme $D^1_{\delta_j}\Psi(\pi)\geqslant 0$ pour tout $j$ on peut donc écrire :
%\begin{align*}
%D^1_{\mu}\Psi(\pi)=&\sum_{j=0}^L\mu(j)^2D^1_{\delta_j}\Psi(\pi)+2\sum_{j_1,j_2=0; j_1\neq j_2}^L\mu_{j_1}\mu_{j_2}\sqrt{D^1_{\delta_j}\Psi(\pi)}\\
%=&\sum_{j=0}^L\mu(j)D^1_{\delta_j}\Psi(\pi).
%\end{align*}
Because $\pi$ satisfies \textit{2.}, $D_{\mu}\Psi(\pi)\geqslant 0$, and finally $\Psi(\mu)-\Psi(\pi)\geqslant 0$ and $\pi=\hat{\pi}_L$.\\
To conclude assume now that $j\in$ supp($\hat{\pi}_{L}$). Then the function $(1+\epsilon)\hat{\pi}_L-\epsilon\delta_j$ is a probability for $\epsilon$ positive small enough, and we have the following inequality:
\begin{eqnarray*}
-D_{\delta_j}\Psi(\pi)=\underset{\epsilon\searrow 0^+}{\lim}\frac{1}{\epsilon}\big(\Psi ((1+\epsilon)\hat{\pi}_L-\epsilon\delta_j)-\Psi (\hat{\pi}_L)\big)\geqslant 0,
\end{eqnarray*}
which concludes the proof of the lemma.

\paragraph{Proof of Lemma \ref{lemme5}}
Let $\pi_S$ be the solution of the first problem of minimization. Let $Q_{S}$ and $H_{S}$ be defined as in Section \ref{algorithm.section}.
The KKT's conditions give us that $\pi_S$ is the unique sequence that satisfies:
\begin{equation}
\exists \lambda_S\in\mathbb{R}, 
\left\{
    \begin{array}{l}
    \vspace{0.25cm}
    \sum_{j\in S}\pi_{S}(j)=1\\
    \frac{\partial}{\partial \pi} \mathcal{L}(\pi_S,\lambda_S)=0
    \end{array}
 \right.
\end{equation}
where $\mathcal{L}$ is the Lagrange's function:
\begin{eqnarray*}
\mathcal{L}(\pi,\lambda)=\Psi(\pi)+\lambda(\sum_{j\in S}\pi(j)-1).
\end{eqnarray*}

The partial derivative function of $\mathcal{L}$ is:
\begin{eqnarray*}
\frac{\partial}{\partial \pi} \mathcal{L}(\pi,\lambda)=-Q_S^T(\tilde{p}-Q_S\pi)+\lambda Q_S^T\mathbb{I},
\end{eqnarray*}
where $\mathbb{I}$ is the vector with $L+1$ components equal to $1$. We have
\begin{eqnarray*}
\pi_S=(Q_S^TQ_S)^{-1}Q_S^T(\tilde{p}-\lambda_S\mathbb{I}) \text{ and } Q_S\pi_S=H_S(\tilde{p}-\lambda_S\mathbb{I}),
\end{eqnarray*}
leading to:
\begin{eqnarray*}
1=<Q_S\pi_S,\mathbb{I}>=<H_S\tilde{p},\mathbb{I}>-\lambda_S<H\mathbb{I},\mathbb{I}>.
\end{eqnarray*}
Finally we obtain:
\begin{eqnarray*}
\lambda_S=\frac{<H\tilde{p},\mathbb{I}>-1}{<H\mathbb{I},\mathbb{I}>}.
\end{eqnarray*}

Then for all $\pi$ with support included on $S$ we have $\mathcal{L}(\pi_S,\lambda_S)\leqslant\mathcal{L}(\pi,\lambda_S)$ and $\pi_S$ is solution for the second problem:
\begin{eqnarray*}
\pi_S=\underset{\text{supp}(\pi)\subset S}{\text{argmin}}\big(\mathcal{L}(\pi,\lambda_S)\big).
\end{eqnarray*}
Because we are considering strictly convex minimization problems, we get $\pi_{S}=\pi^{'}_{S}$.

\subsection{Proof of properties}

\subsubsection{Proof of the link between $k-$monotony and $(k-1)$-monotony\label{monotony.sec}}

We will prove the following property about $k-$monotone discrete functions:
\begin{Prop}
For all $k \geqslant 2$, if $p\in L^1(\mathbb{N})$ is a $k-$monotone
discrete function then $p$ is $j-$monotone and strictly $j-$monotone
on its support for all $j<k$. 
\end{Prop}

We show this result by iteration. First a convex (or $2$-monotone) discrete function on $L^1(\mathbb{N})$ is nonincreasing (see \cite{theseKoladjo}).

Let now $k\geqslant 3$. Let $p\in L^1(\mathbb{N})$ be a $k-$monotone function. We denote $q$ the following discrete function: 
\begin{eqnarray*}
\forall i\in\mathbb{N}, q(i)=(-1)^{k-2}\Delta^{k-2}p(i).
\end{eqnarray*}
The function $q$ is in $L^1(\mathbb{N})$ and $\Delta^{2}q(i)=(-1)^{k}\Delta^{k}p(i) \geqslant 0$ for all $i\in\mathbb{N}$. Therefore $q$ is convex and nonincreasing.\\
It follows that for all $i\in\mathbb{N}$, $(-\Delta^1)((-1)^{k-2} \Delta^{k-2}p(i))=q(i)-q(i+1)\geqslant 0$ i.e. $(-1)^{k-1} \Delta^{k-1}p(i)\geqslant 0$ and
$p$ is $(k-1)-$monotone.

\subsubsection{Proof of Property \ref{theoclé}\label{theoclé.st}}

We prove this property by induction. First for $k=2$, we have the following equalities:
\begin{align*}
F^{2}_{f}(l)-F^{2}_{\tilde{p}}(l)=&\sum_{h_1=0}^{l}\sum_{h_2=0}^{h_1}(f(h_2)-\tilde{p}(h_2))\\
   =&\sum_{h_1=0}^{s}\sum_{h_2=0}^{h_1}(f(h_2)-\tilde{p}(h_2))+\sum_{h_1=s+1}^{l}\sum_{h_2=0}^{h_1}(f(h_2)-\tilde{p}(h_2))\\
   =&\;F^{2}_{f}(s)-F^{2}_{\tilde{p}}(s)+\sum_{h_1=s+1}^{l}\sum_{h_2=0}^{s}(f(h_2)-\tilde{p}(h_2))\\
   =&\;F^{2}_{f}(s)-F^{2}_{\tilde{p}}(s)+(F^{1}_{f}(s)-F^{1}_{\tilde{p}}(s))(l-s)_+.
\end{align*}
Because $\bar{Q}_{l-1}^2(s)=(l-s)_+$ the property is true for $k=2$.
\\
Assume now that the property is true until the rank $k-1$. We have the following properties:
\begin{align*}
F^{k}_{f}(l)-F^{k}_{\tilde{p}}(l)=& \sum_{h=0}^{l}(F^{k-1}_{f}(h)-F^{k-1}_{\tilde{p}}(h))\\
=& \sum_{h=0}^{s}(F^{k-1}_{f}(h)-F^{k-1}_{\tilde{p}}(h))+\sum_{h=s+1}^{l}(F^{k-1}_{f}(l)-F^{k-1}_{\tilde{p}}(l))\\
=& \;F^{k}_{f}(s)-F^{k}_{\tilde{p}}(s)+\sum_{h=s+1}^{l}\big(\sum_{j=1}^{k-1}\bar{Q}_{h-1}^{k-1-j+1}(s)\big(F^j_{f}(s)-F^{j}_{\tilde{p}}(s)\big)\big).
\end{align*}
The last equality is obtained by iteration. Using the definition of the $Q_j^k$ we get:
\begin{align*}
\bar{Q}_{h-1}^{k-j}(l)=&\; \bar{Q}_{l-1}^{k-j}(l-h+s)\\
             =& \;\bar{Q}_{l-1}^{k-j+1}(l-h+s)-\bar{Q}_{l-1}^{k-j+1}(l-h+s+1),
\end{align*}
and the following equalities:
\begin{align*}
F^{k}_{f}(l)-F^{k}_{\tilde{p}}(l)=& \;F^{k}_{f}(s)-F^{k}_{\tilde{p}}(s)\\
&\;+\sum_{j=1}^{k-1}\big(F^j_{f}(s)-F^{j}_{\tilde{p}}(s)\big)\sum_{h=s+1}^{1}\big(\bar{Q}_{l-1}^{k-j+1}(l-h+s)-\bar{Q}_{l-1}^{k-j+1}(l-h+s+1)\big)\\
=& \;F^{k}_{f}(s)-F^{k}_{\tilde{p}}(s)+ \sum_{j=1}^{k-1}\big(F^j_{f}(s)-F^{j}_{\tilde{p}}(s)\big)\big(\bar{Q}_{l-1}^{k-j+1}(s)-\bar{Q}_{l-1}^{k-j+1}(l)\big).
\end{align*}
Because $\bar{Q}_{l-1}^{k-j+1}(l)=0$ and $\bar{Q}^{k-k+1}_{l-1}(s)=1$, we finally obtain:
\begin{eqnarray*}
F^{k}_{f}(l)-F^{k}_{\tilde{p}}(l)=\sum_{j=1}^{k}\bar{Q}_{l-1}^{k-j+1}(s)\big(F^j_{f}(s)-F^{j}_{\tilde{p}}(s)\big).
\end{eqnarray*}

\subsubsection{Proof of Properties \ref{beta} \label{beta.st}\label{theo1.st}}

The following property gives a characterization of the estimator $\hat{p}^*$:
\begin{Prop}\label{theo1}
Let $f\in L^1(\mathbb{N})$. There is an equivalence between :
\begin{enumerate}
\item \begin{itemize}
\item For all $l\in\mathbb{N}$ we have $F^{k}_{f}(l)\geqslant F^{k}_{\tilde{p}}(l).$
\item If $l$ is a $k-$knot of $f$, then the previous inequality is an equality. 
\end{itemize}
\item $f=\hat{p}^{*}$.
\end{enumerate}
\end{Prop}
The proof is similar to the proof of Theorem \ref{theo2} and is omitted. Property \ref{beta} is deduced from Property \ref{theo1}.

\subsubsection{The mass of $\widehat{p}^{*}$ is greater than 1\label{masse.prop.st}}

Let $s_{\max}$ the maximum of $\hat{s}^*$ and $\tilde{s}$ (the maxima of the supports of $\hat{p}^{*}$ and $\tilde{p}$ respectively), then using Property \ref{theoclé}, for all $l\geqslant s_{\max}$ we have:
\begin{eqnarray*}
F^{k}_{\hat{p}^{*}}(l)-F^{k}_{\tilde{p}}(l)=\sum_{j=1}^{k}\bar{Q}_{l-1}^{k-j+1}(s_{\max})\big(F^j_{\hat{p}^{*}}(s_{\max})-F^{j}_{\tilde{p}}(s_{\max})\big).
\end{eqnarray*}
Because the quantities $\bar{Q}_{l-1}^{j}(s_{\max})$ are polynomial functions of $l-s_{\max}$ with degree $j-1$, we get:
\begin{align*}
F^{k}_{\hat{p}^{*}}(l)-F^{k}_{\tilde{p}}(l) =&F^{1}_{\hat{p}^{*}}(s_{\max})-F^{1}_{\tilde{p}}(s_{\max})\frac{(l-s_{\max})^{k-1}}{(k-1)!}+o(l^{k-1})\\
=&(\sum_{j=0}^{s_{\max}}\hat{p}^{*}(j)-1)\frac{(l-s_{\max})^{k-1}}{(k-1)!}+o(l^{k-1}).
\end{align*}

If $\sum_{j=0}^{s_{\max}}\hat{p}^{*}(j)<1$ then, when $l$ tends to infinity, the right-hand term tends to $-\infty$ and the left-hand term is non-negative by Theorem \ref{theo1}. Therefore $\sum_{j=0}^{s_{\max}}\hat{p}^{*}(j)\geqslant 1$.
\\

\subsubsection{Proof of Property \ref{prop5} \label{prop5.st}}

By Theorem \ref{maxsupport.prop} almost surely it exists $n_0\in\mathbb{N}$ such as for all $n\geqslant n_0$ we have $\hat{s}^*\leqslant s+1$, where $s$ is the support of $p$ and $\hat{s}^*$ the support of $\hat{p}^*$. Then almost surely:
\begin{eqnarray*}
\forall n\geqslant n_0, \big\vert\sum_{i=0}^{\infty}\hat{p}^*(i)-1\big\vert = \big\vert\sum_{i=0}^{s+1}(\hat{p}^*(i)-p(i))\big\vert\leqslant \sum_{i=0}^{s+1}\vert\hat{p}^*(i)-p(i)\vert.
\end{eqnarray*}
Moreover, almost surely by Theorem \ref{limnorme2.prop} $\vert\vert p-\hat{p}^*\vert\vert_{2}\leqslant \vert\vert p-\tilde{p}_{n}\vert\vert_{2}$ and then $\underset{n\rightarrow\infty}{\lim}\vert\vert\hat{p}^*(i)-p(i)\vert\vert_2^2=0$ i.e.:
\begin{eqnarray*}
\underset{n\rightarrow\infty}{\lim}\sum_{i=0}^{\infty}(\hat{p}^*(i)-p(i)
)^2=0.
\end{eqnarray*}
Then almost surely for all $i\in\mathbb{N}$ we have $\underset{n\rightarrow\infty}{\lim}(\hat{p}^*(i)-p(i)
)=0$ and $\underset{n\rightarrow\infty}{\lim}\sum_{i=0}^{s+1}\vert\hat{p}^*(i)-p(i)\vert=0$. Finally almost surely we have:
\begin{eqnarray*}
\underset{n\rightarrow\infty}{\lim}\big\vert\sum_{i=0}^{\infty}\hat{p}^*(i)-1\big\vert=0.
\end{eqnarray*}

\subsubsection{Proof of Property \ref{Poiss.prop}\label{Poiss.st}}

We prove that Poisson distribution $\mathcal{P}(\lambda)$ is $l-$monotone if and only if $\lambda\leq\lambda_l$. The distribution $q=\mathcal{P}(\lambda)$ is $l-$monotone if and only if for all $i\in\mathbb{N}$ we have $(-1)^k\Delta^kq(i)\geq 0$. We have for all $l\in\mathbb{N}$ the following equalities:
\begin{eqnarray*}
(-1)^k\Delta^kq(i)=\sum_{h=0}^{l}(-1)^h(_h^l)\frac{\lambda^{h+i}e^{-\lambda}}{(h+i)!}=\frac{\lambda^{i}e^{-\lambda}}{(h+l)!}R_l(\lambda,i)
\end{eqnarray*}
where $R_l$ is the polynomial function defined as follows:
\begin{eqnarray*}
R_l(\lambda,i)=\sum_{h=0}^{l}(-1)^h(_h^l)\lambda^h(h+l)\ldots(h+i+1).
\end{eqnarray*}
Therefore a necessary condition for $\mathcal{P}(\lambda)$ to be $l-$monotone is $R_l(\lambda,0)$ nonnegative which is equivalent to $P_l(\lambda)$ nonnegative where $P_l(\lambda)$ is defined as follows:
\begin{eqnarray*}
P_\ell(\lambda)=\sum_{h=0}^\ell (-1)^h \frac{(\ell!)^2}{h!((\ell-h)!)^2}\lambda^h.
\end{eqnarray*} 
Conversely, because $i\hookrightarrow R_l(\lambda,i)$ is an increasing function for $\lambda\in[0,1]$, the condition is sufficient.\\
When $\lambda$ tends to infinity, $P_l(\lambda,0)$ tends to $+\infty$ then 
$P_l(\lambda)$ is nonnegative until the smallest root of $P_l$ which is nonnegative. In other terms the previous condition is true in particular for $\lambda\leq \lambda_l$.

\subsubsection{Projection of $\delta_{1}$ onto ${\mathcal S}^{3}$\label{cex_masse.st}}

Our purpose is to show that the projection of $\delta_1$ on the cone $\mathcal{S}_3$ has a mass strictly greater than one. After some calculations, we know that this projection is written as $f=\alpha\bar{Q}_5^3+\beta\bar{Q}_6^3$. To calculate the coefficients $(\alpha,\beta)$, we need to establish a necessary and sufficient condition which makes sure that $f$ is $\hat{p}^{*3}$. This condition is given in Property \ref{theo1} (see Section \ref{theo1.st}).

We search $\alpha$ et $\beta$ such as $f=\alpha \bar{Q}^3_5+\beta \bar{Q}^3_6$ satisfies the stopping criterion. For this $p$ we have $f=(21\alpha+28\beta,15\alpha+21\beta,10\alpha+15\beta,6\alpha+10\beta,
3\alpha+6\beta,\alpha+3\beta,\beta,0\ldots)$.
With elementary calculations we obtain the following necessary conditions for $\alpha$ et $\beta$:
\begin{eqnarray*}
S1=\left\{
    \begin{array}{l}
    \vspace{0.25cm}
     F^3_p(0)=21\alpha+28\beta\geqslant 0=F^3_{\delta_1}(0)\\
    \vspace{0.25cm}
     F^3_p(1)=78\alpha+105\beta\geqslant 1=F^3_{\delta_1}(1)\\ 
    \vspace{0.25cm}
    F^3_p(2)=181\alpha+246\beta\geqslant 3=F^3_{\delta_1}(2)\\
    \vspace{0.25cm}
    F^3_p(3)=336\alpha+461\beta\geqslant 6=F^3_{\delta_1}(3)\\
    \vspace{0.25cm}
    F^3_p(4)=546\alpha+756\beta\geqslant 10=F^3_{\delta_1}(4)\\
    \vspace{0.25cm}
    F^3_p(5)=812\alpha+1134\beta\geqslant 15=F^3_{\delta_1}(5)\\
    \vspace{0.25cm}
    F^3_p(6)=1134\alpha+1596\beta\geqslant 21=F^3_{\delta_1}(6)\\
    
    F^3_p(7)=1512\alpha+2142\beta\geqslant 28=F^3_{\delta_1}(7)\\
    \end{array}
 \right.
\end{eqnarray*}
and
\begin{eqnarray*}
S2=\left\{
    \begin{array}{l}
    \vspace{0.25cm}
     F^1_p(7)=61\alpha+89\beta\geqslant 1=F^1_{\delta_1}(7)\\ 
    
    F^2_p(2)=378\alpha+546\beta\geqslant 7=F^2_{\delta_1}(7)\\
     \end{array}
 \right.
\end{eqnarray*}
and
\begin{eqnarray*}
S3=\left\{
    \begin{array}{l}
    \vspace{0.25cm}
    (-1)^3\Delta^3p(i)=0 \text{   si $i>8$ ou si $i<5$}\\ 
    \vspace{0.25cm}
    (-1)^3\Delta^3p(6)=2\beta\\

    (-1)^3\Delta^3p(5)=2\alpha.
    \end{array}
 \right.
\end{eqnarray*}
The condition $S1$ assure that $f$ satisfies 1., $S2$ that $p$ satisfies 2.(c) and $S3$ that $p$ satisfies 2.(b).\\
If we assume that $\alpha$ and $\beta$ are strictly nonnegative we find the more restrictive necessary condition:
\begin{eqnarray*}
S4=\left\{
    \begin{array}{l}
    \vspace{0.25cm}
    812\alpha+1134\beta= 15\\ 
    
    1134\alpha+1596\beta= 21
    \end{array}
 \right.
  =\left\{
    \begin{array}{l}
    \vspace{0.25cm}
    812\alpha+1134\beta= 15\\ 
    
    54\alpha+76\beta=1,
    \end{array}
 \right.
\end{eqnarray*}
whose unique solution is
\begin{eqnarray*}
\left\{
    \begin{array}{l}
    \vspace{0.25cm}
    \alpha=\frac{3}{238}\\ 
    
    \beta=\frac{1}{238}.
    \end{array}
 \right.
 \end{eqnarray*}
 
Reciprocally if we take $\alpha$ and $\beta$ like before then $f$ satisfies the conditions $S_1$, $S_2$ and $S3$. Using Property \ref{theo1} it follows that $f$ is the projection of $\delta_1$ on the set of $3-$monotone sequences. 
\\

\subsection{Proofs of the technical lemmas \label{prooflemmas.st}}

Let us first state technical lemmas used in the proofs given before. Their proofs are given afterwards.

\begin{Lemme}For all integer $k\geqslant 2$, for all $l\in\mathbb{N}$ and for all $f\in\mathcal{P}$, the following assumption is true:
\label{rec.lemme}
\begin{eqnarray}
\label{pterec}
\sum_{i=0}^{l}f(i)\bar{Q}^{k}_{l}(i)=F^k_f(l).
\end{eqnarray}
\end{Lemme}

\begin{Lemme}\label{lemme-beta}For all $k\geq 0$, for all $f\in\mathcal{S}_k$, for all $g\in L^1(\mathbb{N})$:
\begin{eqnarray*}
\sum_{i=0}^{\infty}f(i)g(i)=\sum_{i=0}^{\infty}(-1)^k\Delta^kf(i)F^k_{g}(i).
\end{eqnarray*}
In particular for all $f\in \mathcal{S}_k$ the coefficient $\beta(f)$ defined at Equation (\ref{beta.eq}) satisfies:
\begin{eqnarray*}\beta(f)=\sum_{i=0}^{\infty}f(i)(f(i)-\tilde{p}(i))=\sum_{i=0}^{\infty}(-1)^{k}\Delta^{k}f(i)\big(
F^{k}_{f}(i)-F^{k}_{\tilde{p}}(i)\big).
\end{eqnarray*}
\end{Lemme}

\begin{Lemme}\label{signebeta}The coefficient $\beta(\hat{p})$ defined at Equation (\ref{beta.eq}) is always non-positive
\end{Lemme}

\begin{Lemme}\label{lemme16}Let $k\geq 2$. Let $f\in L^1(\mathbb{N})$, $s\in\mathbb{N}$ and $l\geqslant s$. The following equality is true:
\begin{eqnarray*}
F^k_f(l)-F^k_{\tilde{p}}(l)=\frac{(l-s)^{k-1}}{(k-1)!}(F^1_f(s)-F^1_{\tilde{p}}(s))+o(l^{k-1}).
\end{eqnarray*}
\end{Lemme}

\paragraph{Proof of Lemma \ref{rec.lemme}}
The lemma is proved by induction. Let us first consider $k=2$. Let $f$ be a positive sequence and $l\in\mathbb{N}$. We have:
\begin{eqnarray*}
F_f^2(l)=\sum_{h=0}^{l}\sum_{i=0}^{h}f(i)
        =\sum_{i=0}^{l}\sum_{h=i}^{l}f(i)
        =\sum_{i=0}^{l}f(i)(l+1-i)
        =\sum_{i=0}^{l}f(i) \bar{Q}^2_l(i),
\end{eqnarray*}
and Equation (\ref{pterec}) is shown.
Assume that Equation (\ref{pterec}) is true for $k-1\geqslant 2$. We have the following equalities:
\begin{eqnarray*}
F_f^k(l)=\sum_{h=0}^{l}F_f^{k-1}(h)
        =\sum_{h=0}^{l}\sum_{i=0}^{h}f(i) \bar{Q}^{k-1}_h(i)
        =\sum_{i=0}^{l}f(i)\sum_{h=i}^{l}\bar{Q}^{k-1}_h(i).
\end{eqnarray*}
Using Pascal's Triangle and the definition of $\bar{Q}^k_j$, we get:
\begin{align*}       
F_f^k(l)&=\sum_{i=0}^{l}f(i)\sum_{h=i}^{l}
        \big(\bar{Q}^{k}_h(i)-\bar{Q}
        ^{k}_h(i+1)\big) \\
        &=\sum_{i=0}^{l}f(i)\big(\sum_{h=i}^{l}\bar{Q}^{k}_h(i)-
        \sum_{h=i}^{l}\bar{Q}^{k}_{h-1}(i)\big)
\end{align*}
where the last equality comes from $\bar{Q}^{k}_{h}(i+1)=\bar{Q}^{k}_{h-1}(i)$ with the convention $\bar{Q}^{k}_{0}=0$. Finally we obtain:
\begin{eqnarray*}
F_f^k(l)=\sum_{i=0}^{l}f(i)\big(\bar{Q}^{k}_l(i)\big),
\end{eqnarray*}
and the lemma is shown.

\paragraph{Proof of Lemma \ref{lemme-beta}}

The lemma is proved by induction. First it is true for $k=0$ with the convention $\Delta^0 f(i)=f(i)=F^0_f(i)$. Assume now that the result if true for some $k-1\geqslant 0$. We have the following inequalities:
\begin{align*}
\sum_{i=0}^{\infty}\Delta^k f(i)F^k_{g}(i)=&\sum_{i=0}^{\infty}(\Delta^{k-1} f(i+1)-\Delta^{k-1} f(i))F^k_{g}(i)\\
=&\sum_{i=1}^{\infty}\Delta^{k-1} f(i)F^k_{g}(i-1)
-\sum_{i=0}^{\infty}\Delta^{k-1} f(i)F^k_{g}(i)\\
=&\sum_{i=1}^{\infty}\Delta^{k-1} f(i)[F^k_{g}(i-1)-F^k_{g}(i)]
-\Delta^{k-1} f(0)F^k_{g}(0)\\
=&-\sum_{i=1}^{\infty}\Delta^{k-1} f(i)F^{k-1}_{g}(i))
-\Delta^{k-1} f(0)F^{k-1}_{g}(0)
\end{align*}
because $F^k_{g}(0)=\sum_{h_1=0}^{0}\ldots\sum_{h_k=0}^0 g(h_k)=f(0)=F^{k-1}_g(0).$
\begin{Rem}This sums are well-defined because thanks to Lemma \ref{rec.lemme} ((see Section \ref{prooflemmas.st}) we have:
\begin{align*}
\sum_{l=0}^{\infty}\vert F^k_g(l)\Delta^kf(l) \vert=&
\sum_{l=0}^{\infty}\sum_{i=0}^lg(i)\bar{Q}_l^k(i)(-1)^k\Delta^kf(l)\\
=&\sum_{i=0}^{\infty}\left(\sum_{l=i}^{\infty}(-1)^k\Delta^kf(l)\bar{Q}_l^k(i)\right)g(i).
\end{align*}
By Property \ref{decomposition2} (see Section \ref{characterization.section}) we have the equality:
\begin{eqnarray*}
f(i)=\sum_{l=0}^{\infty}(-1)^k\Delta^kf(l)\bar{Q}_l^k(i).
\end{eqnarray*}
Then $\sum_{l=i}^{\infty}(-1)^k\Delta^kf(l)\bar{Q}_l^k(i)\leqslant 1$ and finally :
\begin{eqnarray*}
\sum_{l=0}^{\infty}\vert F^k_f(l)\Delta^kf(l)\vert \leqslant \sum_{i=0}^{\infty}f(i)<\infty.
\end{eqnarray*}
\end{Rem}
It follows that $\sum_{i=0}^{\infty}\Delta^k f(i)F^k_{g}(i)=-\sum_{i=0}^{\infty}\Delta^{k-1} f(i)F^{k-1}_{g}(i)$ and the lemma is proved.

\paragraph{Proof of Lemma \ref{signebeta}}

We note $\hat{s}$ and $\tilde{s}$ the maxima of the supports of $\hat{p}$ and $\tilde{p}$ respectively. We note $s_{\max}=\max(\hat{s},\tilde{s})$. We use Property \ref{theoclé} with $f=\hat{p}$ and we obtain that for all $l\geqslant s_{\max}+1$:
\begin{align*}
F^k_{\hat{p}}(l)-F^k_{\tilde{p}}(l)&=\sum_{j=1}^k\bar{Q}_{l-1}^{k-j+1}(s_{\max})\big(F^j_{\hat{p}}(s_{\max})-F^j_{\tilde{p}}(s_{\max})\big)\\
&=\sum_{j=2}^k\bar{Q}_{l-1}^{k-j+1}(s_{\max})\big(F^j_{\hat{p}}(s_{\max})-F^j_{\tilde{p}}(s_{\max})\big).
\end{align*}
The last equality comes from $F^1_{\hat{p}}(s_{\max})=F^1_{\tilde{p}}(s_{\max})=1$ because $\hat{p}$ and $\tilde{p}$ are probabilities and $s_{\max}$ is greater than $\hat{p}$ and $\tilde{p}$. \\
Because the quantities $\bar{Q}^j_{l-1}(s_{\max})$ are polynomial functions with degree $j-1$ in the variable $l-s_{\max}$ we write $F^k_{\hat{p}}(l)-F^k_{\tilde{p}}(l)$ in the following form:
\begin{eqnarray*}
F^k_{\hat{p}}(l)-F^k_{\tilde{p}}(l)=\frac{(F^2_{\hat{p}}(s_{\max})-F^2_{\tilde{p}}(s_{\max}))}{(k-2)!}(l-s)^{k-2}+o(l^{k-2}).
\end{eqnarray*}
Thanks to Equation~\eref{massespline}, $m^k_l$ is a polynomial function with degree $k$ and we have the following limit:
\begin{eqnarray*}
\underset{l\rightarrow\infty}{\lim}\frac{F^k_{\hat{p}}(l)-F^k_{\tilde{p}}(l)}{m^k_l}=0.
\end{eqnarray*}
Moreover for all $l\in\mathbb{N}$ the characterization of  $\hat{p}$ gives us:
\begin{eqnarray*}
\frac{F^k_{\hat{p}}(l)-F^k_{\tilde{p}}(l)}{m^k_l}\geqslant\beta(\hat{p}).
\end{eqnarray*}
Necessarily $\beta(\hat{p})\leqslant 0$.

\paragraph{Proof of Lemma \ref{lemme16}}

We show this result by induction. For $k=2$ the result is shown in \cite{durot_least}. Assume that the result is true for some $k-1\geqslant 2$. We have the following equalities:
\begin{align*}
F^k_f(l)-F^k_{\tilde{p}}(l)&=\sum_{h=0}^{l}\left(F^{k-1}_f(h)-F^{k-1}_{\tilde{p}}(h)\right)\\
&=\sum_{h=0}^{s}\left(F^{k-1}_f(h)-F^{k-1}_{\tilde{p}}(h)\right)+\sum_{h=s+1}^{l}\left(F^{k-1}_f(h)-F^{k-1}_{\tilde{p}}(h)\right)\\
&=\big(F^k_f(s)-F^k_{\tilde{p}}(s)\big)+\sum_{h=s+1}^{l}\left(\frac{(h-s)^{k-2}}{(k-2)!}\big(F^{1}_f(s)-F^{1}_{\tilde{p}}(s)\big)+o(h^{k-2})\right)\\
&=\frac{(l-s)^{k-1}}{(k-1)!}(F^1_f(s)-F^1_{\tilde{p}}(s))+o(l^{k-1}).
\end{align*}
The last equality is due to a result of Bernoulli for Faulhaber's sum : the $k-$th sum of Faulhaber is denoted by $S_k$ and defined as follows :
\begin{eqnarray*}
S_k(m)=\sum_{i=1}^{m}i^k.
\end{eqnarray*}
It is shown that:
\begin{eqnarray*}
S_k(m)=\frac{1}{k+1}\sum_{j=0}^kC^{j}_{k+1}B_jm^{k+1-j}=\frac{1}{k+1}\big(m^{k+1}+\frac{k+1}{2}m^k+\ldots\big)
\end{eqnarray*}
where the $B_j$ are Bernoulli's numbers (with the convention $B_1=\frac{1}{2}$). A proof of this result can be found in \cite{faulhaber}.

\bibliographystyle{plain}
\bibliography{bibliothese}

\begin{thebibliography}{10}

\bibitem{balabdaoui2004nonparametric}
Fadoua Balabdaoui.
\newblock {\em Nonparametric estimation of a k-monotone density: A new
  asymptotic distribution theory.}
\newblock PhD thesis, University of Washington, 2004.

\bibitem{balabdaoui2014asymptotics}
Fadoua Balabdaoui, C{\'e}cile Durot, and Fran{\c{c}}ois Koladjo.
\newblock On asymptotics of the discrete convex lse of a pmf.
\newblock {\em arXiv preprint arXiv:1404.3094}, 2014.

\bibitem{balabdaoui2013asymptotics}
Fadoua Balabdaoui, Hanna Jankowski, Kaspar Rufibach, and Marios Pavlides.
\newblock Asymptotics of the discrete log-concave maximum likelihood estimator
  and related applications.
\newblock {\em Journal of the Royal Statistical Society: Series B (Statistical
  Methodology)}, 75(4):769--790, 2013.

\bibitem{balabdaoui2009limit}
Fadoua Balabdaoui, Kaspar Rufibach, and Jon~A Wellner.
\newblock Limit distribution theory for maximum likelihood estimation of a
  log-concave density.
\newblock {\em Annals of statistics}, 37(3):1299, 2009.

\bibitem{balabdaoui2007estimation}
Fadoua Balabdaoui and Jon~A Wellner.
\newblock Estimation of a k-monotone density: limit distribution theory and the
  spline connection.
\newblock {\em The Annals of Statistics}, 35:2536--2564, 2007.

\bibitem{balabdaoui2010estimation}
Fadoua Balabdaoui and Jon~A Wellner.
\newblock Estimation of a k-monotone density: characterizations, consistency
  and minimax lower bounds.
\newblock {\em Statistica Neerlandica}, 64(1):45--70, 2010.

\bibitem{boyd2004convex}
Stephen Boyd and Lieven Vandenberghe.
\newblock {\em Convex optimization}.
\newblock Cambridge university press, 2004.

\bibitem{bunge2014estimating}
John Bunge, Amy Willis, and Fiona Walsh.
\newblock Estimating the number of species in microbial diversity studies.
\newblock {\em Annual Review of Statistics and Its Application}, 1:427--445,
  2014.

\bibitem{faulhaber}
J.~H. Conway and R.~K. Guy.
\newblock {\em The Book of Numbers}.
\newblock Springer-Verlag, 1996.

\bibitem{dumbgen2010logcondens}
Lutz D{\"u}mbgen and Kaspar Rufibach.
\newblock logcondens: Computations related to univariate log-concave density
  estimation.
\newblock {\em Journal of Statistical Software}, 39:1--28, 2011.

\bibitem{dumbgen2009maximum}
Lutz D{\"u}mbgen, Kaspar Rufibach, et~al.
\newblock Maximum likelihood estimation of a log-concave density and its
  distribution function: Basic properties and uniform consistency.
\newblock {\em Bernoulli}, 15(1):40--68, 2009.

\bibitem{durot_least}
C{\'e}cile Durot, Sylvie Huet, Fran{\c{c}}ois Koladjo, and St{\'e}phane Robin.
\newblock Least-squares estimation of a convex discrete distribution.
\newblock {\em Computational Statistics \& Data Analysis}, 67:282--298, 2013.

\bibitem{durot_nonparametric}
C{\'e}cile Durot, Sylvie Huet, Fran{\c{c}}ois Koladjo, and St{\'e}phane Robin.
\newblock Nonparametric species richness estimation under convexity constraint.
\newblock {\em Environmetrics}, 26(7):502--513, 2015.

\bibitem{fejer1936trigonometrische}
Leopold Fej{\'e}r.
\newblock Trigonometrische reihen und potenzreihen mit mehrfach monotoner
  koeffizientenfolge.
\newblock {\em Transactions of the American Mathematical Society},
  39(1):18--59, 1936.

\bibitem{feller1939completely}
Willy Feller et~al.
\newblock Completely monotone functions and sequences.
\newblock {\em Duke Math. J}, 5:661--674, 1939.

\bibitem{fisher1943relation}
Ronald~Aylmer Fisher, A~Steven Corbet, and Carrington~B Williams.
\newblock The relation between the number of species and the number of
  individuals in a random sample of an animal population.
\newblock {\em The Journal of Animal Ecology}, 12:42--58, 1943.

\bibitem{grenander1956theory}
Ulf Grenander.
\newblock On the theory of mortality measurement: part ii.
\newblock {\em Scandinavian Actuarial Journal}, 1956(2):125--153, 1956.

\bibitem{groeneboom2001estimation}
Piet Groeneboom, Geurt Jongbloed, and Jon~A Wellner.
\newblock Estimation of a convex function: Characterizations and asymptotic
  theory.
\newblock {\em Annals of Statistics}, 29:1653--1698, 2001.

\bibitem{groeneboom2008support}
Piet Groeneboom, Geurt Jongbloed, and Jon~A Wellner.
\newblock The support reduction algorithm for computing non-parametric function
  estimates in mixture models.
\newblock {\em Scandinavian Journal of Statistics}, 35(3):385--399, 2008.

\bibitem{jankowski2009estimation}
Hanna~K Jankowski and Jon~A Wellner.
\newblock Estimation of a discrete monotone distribution.
\newblock {\em Electronic journal of statistics}, 3:1567, 2009.

\bibitem{jewell1982mixtures}
Nicholas~P Jewell.
\newblock Mixtures of exponential distributions.
\newblock {\em The Annals of Statistics}, 10:479--484, 1982.

\bibitem{knopp1925mehrfach}
Konrad Knopp.
\newblock Mehrfach monotone zahlenfolgen.
\newblock {\em Mathematische Zeitschrift}, 22(1):75--85, 1925.

\bibitem{theseKoladjo}
Fran{\c{c}}ois Koladjo.
\newblock {\em Estimation d'une distribution discrete sous contrainte de
  convexit{\'e} : Application a l'estimation du nombre d'especes de la faune
  ichtyologique du bassin du fleuve Ou{\'e}m{\'e}.}
\newblock PhD thesis, Universit{\'e} Paris-Sud XI et d'Abomey-Calavi, 2013.

\bibitem{lefevre2013multiply}
Claude Lefevre, St{\'e}phane Loisel, et~al.
\newblock On multiply monotone distributions, continuous or discrete, with
  applications.
\newblock {\em Journal of Applied Probability}, 50(3):827--847, 2013.

\bibitem{levy1962extensions}
Paul L{\'e}vy.
\newblock Extensions d'un th{\'e}or{\`e}me de d. dugu{\'e} et m. girault.
\newblock {\em Probability Theory and Related Fields}, 1(2):159--173, 1962.

\bibitem{williamson1955multiply}
Richard~Edmund Williamson et~al.
\newblock {\em On Multiply Monotone Functions and Their Laplace Transforms.}
\newblock PhD thesis, Graduate School of Arts and Sciences, University of
  Pennsylvania, 1955.

\end{thebibliography}
\nocite{*}

\end{document}